\documentclass[reqno]{amsart}

\usepackage{amsmath}
\usepackage{amscd}
\usepackage{amsthm}
\usepackage{amsfonts}
\usepackage{amssymb}

\newcommand{\divides}{\mid}
\newcommand{\notdivides}{\nmid}
\DeclareMathOperator{\Hl}{{\mathit{Hl}}}
\DeclareMathOperator{\K}{{\mathsf{K}}}
\DeclareMathOperator{\AG}{AG}

\renewcommand\ge\geqslant
\renewcommand\geq\geqslant
\renewcommand\le\leqslant
\renewcommand\leq\leqslant

\theoremstyle{plain}
\newtheorem{theorem}{Theorem}[section]
\newtheorem{proposition}[theorem]{Proposition}
\newtheorem{lemma}[theorem]{Lemma}
\newtheorem{corollary}{Corollary}[theorem]

\theoremstyle{definition}
\newtheorem{definition}{Definition}[section]

\theoremstyle{remark}
\newtheorem*{remark*}{Remark}

\newtheorem*{altproof}{Alternative Proof}

\numberwithin{equation}{section}

\begin{document}

\title{On Rationally Parametrized Modular Equations}
\author{Robert S. Maier}
\date{}

\begin{abstract}
Many rationally parametrized elliptic modular equations are derived.  Each
comes from a family of elliptic curves attached to a genus-zero congruence
subgroup~$\Gamma_0(N),$ as an algebraic transformation of elliptic curve
periods, parametrized by a Hauptmodul (function field generator).  The
periods satisfy a Picard--Fuchs equation, of hypergeometric, Heun, or more
general type; so the new modular equations are algebraic transformations of
special functions.  When $N=4,3,2,$ they are modular transformations of
Ramanujan's elliptic integrals of signatures $2,3,4$.  This gives a modern
interpretation to his theories of integrals to alternative bases: they
are attached to certain families of elliptic curves.  His anomalous theory
of signature~$6$ turns~out to fit into a general Gauss--Manin rather than a
Picard--Fuchs framework.
\end{abstract}


\maketitle

\section{Introduction}

\subsection{Context and overview}

The theory of elliptic modular equations is classical, and predates Gauss's
introduction of the homogeneous modular group ${\it PSL}(2,\mathbb{Z})$.
It can be traced to Landen's tranformation of the first complete elliptic
integral $\K=\K(\alpha),$ where the independent
variable~$\alpha$ is often denoted~$k^2,$ after Jacobi.  The
function~$\K(\alpha)$ is the $r=2$ (i.e., `signature~$2$'
or~`base~$2$') case of Ramanujan's elliptic integral
\begin{equation}
\K_r(\alpha) := \tfrac{\sin(\pi/r)}2 \int_0^1 x^{-1/r}
(1-x)^{-1+1/r} (1-\alpha x)^{-1/r}\,dx,
\end{equation}
which is defined on $0< \alpha< 1,$ and extends to a single-valued analytic
function on the Riemann sphere $\mathbb{P}^1(\mathbb{C})_\alpha,$ slit
between $\alpha=1$ and $\alpha=\infty$.  (Without the slit, it would be
multivalued.)  Landen's transformation is
\begin{subequations}
\label{eq:unum}
\begin{equation}
  \label{eq:one}
  \K(\alpha) = (2/\alpha)(1-\sqrt{1-\alpha})\, \K(\beta),
\end{equation}
where $\alpha,\beta$ are constrained by the modular relation
\begin{equation}
  \label{eq:two}
  \alpha^2(1-\beta)^2 - 16(1-\alpha)\beta=0.
\end{equation}
\end{subequations}
A uniformized version of~(\ref{eq:unum}) is
\begin{equation}
\label{eq:duo}
\K\biggl(\frac{t(t+8)}{(t+4)^2}\biggr)=2\left[\frac{t+4}{t+8}\right]\,\K\biggl(\frac{t^2}{(t+8)^2}\biggr),
\end{equation}
where $t$~is an auxiliary parameter.  Landen's transformation is a special
function identity; in~particular, a quadratic hypergeometric
transformation, as $\K_r(\cdot)$~equals $2/\pi$ times the Gauss
hypergeometric function ${}_2F_1(\frac1r,1-\frac1r;1;\cdot)$.

The transformation theory of $\K:=\K_2$ was of intense
interest to nineteenth-century mathematicians, and led to the modern theory
of elliptic curves.  Ramanujan's alternative integrals~$\K_r$ (for
$r=3,4,6,$ in~addition to~$r=2$), which he introduced while deriving
rapidly convergent series for~$\pi,$ have until recently been much less
well understood.  For most of the twentieth century, elliptic integrals
were best known to applied mathematicians.

In recent years, interest in elliptic integrals and modular equations among
pure mathematicians has revived~\cite{Borwein87,Borwein91,McKean99}.
In~part, it has been stimulated by a desire to understand Ramanujan's
modular equations: to~derive them
algorithmically~\cite{Berndt95,Berndt2001}, and also place them in a modern
conceptual framework.  The proper setting of identities such as
(\ref{eq:unum}) or~(\ref{eq:duo}) is now felt to be a general one, a
Gauss--Manin connection over the base curve~$X$ of an elliptic surface
$\mathfrak{E}\stackrel{\pi}{\to}X$.  The generic fibre of such a surface is
an elliptic curve over~$\mathbb{C},$ so the surface can be viewed as an
elliptic family, parametrized by the base curve, which may not be of genus
zero.  This setting subsumes the classical (Picard--Fuchs) situation where
$\mathfrak{E}\stackrel{\pi}{\to}X$ is an elliptic {\em modular\/} surface,
i.e., (i)~$X=\Gamma\setminus\mathcal{H}^*,$ the quotient of the
compactified upper half plane $\mathcal{H}^*\ni\tau=\tau_1/\tau_2$ by a
subgroup $\Gamma<{\it PSL}(2,\mathbb{Z})$ of finite index, and
(ii)~$\mathfrak{E}={\mathfrak{E}}_\Gamma,$ the family of elliptic curves
attached to~$\Gamma$.  The Picard--Fuchs equation
for~$\mathfrak{E}_\Gamma\to\Gamma\setminus\mathcal{H}^*,$ which has
elliptic curve periods as its solutions, defines a Gauss--Manin connection
on a two-dimensional period bundle over~$X$.

The classical Picard--Fuchs theory is the natural setting of Landen's
transformation in the form~(\ref{eq:unum}).  The parameter~$\alpha$ can be
viewed as a Hauptmodul (rational parameter) for the genus-zero modular
curve $X_0(4)\cong\mathbb{P}^1(\mathbb{C}),$ the quotient
of~$\mathcal{H}^*$ by the Hecke subgroup $\Gamma_0(4)$.  Likewise,
$\K(\alpha)$~is a weight-$1$ modular form (with character)
for~$\Gamma_0(4)$.  The second-order differential equation satisfied
by~$\K$ as a function of~$\alpha,$ i.e., the Gauss\ hypergeometric
equation, is the Picard--Fuchs equation attached to~$\Gamma_0(4)$.  It
defines a flat connection on a $2$-plane period bundle over~$X_0(4)$.
Landen's transformation is a relation on
${\mathfrak{E}_{\Gamma_0(4)}\stackrel{\pi}{\to}X_0(4)}\cong\mathbb{P}^1(\mathbb{C})_\alpha,$
which ties together fibres (elliptic curves) over distinct points
$\alpha,\beta\in X_0(4)$ if and only if they are related by a $2$-isogeny.
The `multiplier' $\K(\alpha)/\K(\beta)$ is a quotient of period ratios.  It
is automorphic of weight~$0$ and must be rational in~$\alpha$; or more
accurately (taking characters into account) finite-valued, i.e., algebraic
in~$\alpha,$ as one sees from~(\ref{eq:one}).

The {\em computational\/} theory of modular equations for Gauss--Manin
connections is still in an incomplete state, even for elliptic-modular
families ${\mathfrak{E}}_\Gamma\to\Gamma\setminus\mathcal{H}^*$; and even
in the genus-zero case, when the differential equation defining the flat
connection on the $2$-plane bundle over the base is, in~effect, a Fuchsian
equation on~$\mathbb{P}^1(\mathbb{C})$.  We~recently began the systematic
derivation of such modular equations, viewed as special function
identities, i.e., algebraic transformations of ${}_2F_1$~\cite{Maier11}.
We began with certain genus-$0$ elliptic families, attached to subgroups
$\Gamma<{\it PSL}(2,\mathbb{Z})$; and more generally to their extensions by
Atkin--Lehner involutions, which are subgroups not of ${\it
PSL}(2,\mathbb{Z})$ but of ${\it PSL}(2,\mathbb{R})$.  It soon became clear
that the family of transformations of~${}_2F_1$ of modular origin is larger
than previous treatments had revealed.  {\em Rational\/} transformation
of~${}_2F_1$ were investigated in the nineteenth century, most
systematically by Goursat; though Vid\=unas~\cite{Vidunas2005} has recently
shown that Goursat's classification was incomplete.  {\em Algebraic\/}
tranformations of~${}_2F_1$ are much more numerous, and our treatment in
Ref.~\cite{Maier11} only scratched the surface.

Our ultimate goal is determining which transformations of~${}_2F_1$ come
`from geometry'.  This includes known ${}_2F_1$ and special function
identities, such as Ramanujan's.  As a first step, in this article we go
beyond Landen's transformation by systematically working~out all {\em
rationally parametrized\/} modular equations associated with the
$14$~modular curves $X_0(M)$ that are of genus zero.  The parameter in the
degree-$N$ modular equation is a parameter for~$X_0(NM)$; so such an
equation exists if and only if $\Gamma_0(NM)$ like~$\Gamma_0(M)$ is of
genus zero.  Our equations include (i)~correspondences between elliptic
curves related by $N$-isogenies, i.e., between points on the base of an
elliptic family $\mathfrak{E}_M\stackrel{\pi_M}{\to}X_0(M),$ which are
analogous to the $\alpha$--$\beta$ relation~(\ref{eq:two}); and (ii)~full
modular equations, analogous to~(\ref{eq:one}).  We express these in~terms
of a canonical Hauptmodul $t_M=t_M(\tau)$ and a weight-$1$ modular
form~$\mathfrak{h}_M=\mathfrak{h}_M(\tau),$ treated initially as a
(multivalued) function of~$t_M,$ i.e., $\mathfrak{h}_M=h_M(t_M)$.  Our key
Theorem~\ref{thm:key} says essentially the following.  If
$\nabla_M,\nabla_{NM}$ are the Gauss--Manin connections coming from the
Picard--Fuchs equations for~$\mathfrak{E}_M,\mathfrak{E}_{NM},$ then
pulling back~$\nabla_M$ along the maps $t_{NM}(\tau)\mapsto t_M(\tau)$ and
$t_{NM}(\tau)\mapsto t_M(N\tau)$ yields the {\em same\/} connection;
namely,~$\nabla_{NM}$.  The equality between the two yields a rationally
parametrized degree-$N$ modular equation for the multivalued
function~$h_M,$ analogous to~(\ref{eq:duo}).

To facilitate the interpretation of our modular equations as special
function identities, we first express each weight-$1$ modular
form~$\mathfrak{h}_M$ (regarded as a function~$h_M$ of the corresponding
Hauptmodul~$t_M$) in~terms of~${}_2F_1,$ and give the Picard--Fuchs
equation that $h_M$ satisfies.  To place the identities in context, we also
express each Hauptmodul~$t_M$ and modular form~$\mathfrak{h}_M$ in~terms of
the Dedekind eta function, and give explicit $q$-expansions where
appropriate.  These $q$-expansions yield combinatorial identities
resembling those of Fine~\cite{Fine88}.  The $q$-expansions of the modular
forms tend to be simple, but those of the Hauptmoduln are complicated, and
we mostly omit them.  In~fact, in deriving Picard--Fuchs and modular
equations, we do~not rely on $q$-series at~all.  This contrasts with recent
work of Lian and Wiczer~\cite{Lian2006}, who derived Picard--Fuchs
equations for $175$ genus-zero subgroups of ${\it PSL}(2,\mathbb{R})$ by
$q$-series manipulations.

The culmination of this article is
\S\S\ref{sec:ramanujan}--\ref{sec:final}, where we succeed in placing
Ramanujan's theories of elliptic integrals to alternative bases, which have
been developed by Berndt, Bhargava, and Garvan~\cite{Berndt95} among
others, in a modern setting.  We show that Ramanujan's modular equations
for his elliptic integral~$\K_r,$ where the signature~$r$ equals
$2,3,4,$ come from elliptic families parametrized by $X_0(4),\allowbreak
X_0(3),\allowbreak X_0(2),$ respectively.  In~fact, the (multivalued)
functions $\K_2,\K_3,\K_4$ can be viewed as
defining weight-$1$ modular forms
$\mathcal{A}_2,\mathcal{A}_3,\mathcal{A}_4$ for
$\Gamma_0(4),\Gamma_0(3),\Gamma_0(2),$ which are modified versions
of~$\mathfrak{h}_4,\mathfrak{h}_3,\mathfrak{h}_2$.

Modular interpretations of Ramanujan's theories were pioneered by the
Borweins~\cite{Borwein87,Borwein91}, but this new interpretation is very
fruitful.  It leads to new parametrized modular equations for
$\K_3,\K_4$.  (See our Table~\ref{tab:ramanujan}, which is
likely to be of broad interest.)  We also find that his somewhat mysterious
theory of signature~$6$ is associated to a nonclassical elliptic family
$\mathfrak{E}\stackrel{\pi}{\to}X$ that is not an elliptic {\em modular\/}
family: the base curve~$X$ is not of the
form~$\Gamma\setminus\mathcal{H}^*$.  Hence, his theory of signature~$6$
fits into a general Gauss--Manin rather than a Picard--Fuchs framework.

\subsection{Detailed overview}
\label{subsec:12}

Any elliptic curve over~$\mathbb{C}$ is necessarily of the form
$\mathbb{C}/(\mathbb{Z}\tau_1\oplus\mathbb{Z}\tau_2)\cong\mathbb{C}/(\mathbb{Z}\tau\oplus\mathbb{Z}),$
where the period ratio $\tau=\tau_1/\tau_2$ lies in the upper half
plane~$\mathcal{H}$.  Any
$\left(\begin{smallmatrix}a&b\\c&d\end{smallmatrix}\right)
\in\overline{\Gamma}(1):={\it SL}(2,\mathbb{Z}),$ the inhomogeneous modular
group, acts projectively by $\tau\mapsto\frac{a\tau+b}{c\tau+d},$ giving an
action of $\Gamma(1):={\it PSL}(2,\mathbb{Z})$ on~$\mathcal{H}$.  The space
of isomorphism classes of elliptic curves over~$\mathbb{C}$ is the quotient
$\Gamma(1)\setminus\mathcal{H}$.  Its one-point compactification is the
modular curve $X(1):=\Gamma(1)\setminus\mathcal{H}^*,$ where
$\mathcal{H}^*:=\mathcal{H}\cup\mathbb{P}^1(\mathbb{Q})=\mathcal{H}\cup\mathbb{Q}\cup\{\rm
i\infty\}$ includes cusps.  The function field of~$X(1)$ is generated by
the $j$-invariant, which as traditionally defined has $q$-expansion
$q^{-1}+744 + 196884q + O(q^2)$ about the infinite cusp $\tau={\rm
i}\infty$.  Here $q:=e^{2\pi{\rm i}\tau}$.

For any integer $N>1,$ the algebraic relation ${\Phi_N(j,j')=0},$ with
$j'(\tau):=j(N\tau),$ is called the {\em classical modular equation of
degree~$N$\/}.  The polynomial~$\Phi_N$ is symmetric and
in~$\mathbb{Z}[j,j']$.  Its degree
is
\begin{equation}
\psi(N):=N\prod_{\substack{p\divides{N}\\\text{$p$ prime}}}\Bigl(1+\frac1p\Bigr),
\end{equation}
a multiplicative function of~$N$~\cite{Cox89}.  If $N$~is prime, the $N+1$
roots $j'\in\mathbb{C}$ of the modular equation correspond to
$\tau'=N\tau,$ and to
$\tau/N,(\tau+\nobreak1)/N,\dots,\allowbreak(\tau+\nobreak{N}-\nobreak1)/N$.
In~general the $\psi(N)$ roots correspond to values
$\tau'=\frac{a\tau+b}{c\tau+d},$ where the matrices
$\left(\begin{smallmatrix}a&b\\c&d\end{smallmatrix}\right)\in{\it
GL}(2,\mathbb{Z})$ are reduced level-$N$ modular correspondences, i.e.,
they satisfy $ad=N,$ $c=0,$ $0\le b< d,$ $(a,b,d)=1$.  They are bijective
with the $\psi(N)$ (isomorphism classes~of) unramified $N$-sheeted
coverings of a general elliptic curve
$E=\mathbb{C}/(\mathbb{Z}\tau\oplus\mathbb{Z})$ by an elliptic curve
$\mathbb{C}/(\mathbb{Z}\tau'\oplus\mathbb{Z})$; or equivalently, with the
order-$N$ cyclic subgroups of its group of $N$-division points $E_N\cong
C_N\times C_N$.

The coefficients of~$\Phi_N$ are large, even for small~$N$.  E.g.,
$\Phi_2(j,j')=0$ is
\begin{multline}
\label{eq:modular1}
(j^3+{j'}^3) - j^2{j'}^2 + 2^4\,3\cdot31(j^2j'+j{j'}^2)- 2^4\,3^4\,5^3(j^2+{j'}^2)\\
{}+3^4\,5^3\,4027\,jj' +2^8\,3^7\,5^6(j+j') -2^{12}\,3^9\,5^9=0.
\end{multline}
Due~to the difficulty of computing~$\Phi_N,$ modular equations of
alternative forms have long been of~interest.  One approach begins by
viewing $\Phi_N(j,j')=0$ as a singular plane model of an algebraic curve
over~$\mathbb{C},$ the function field of which is~$\mathbb{C}(j,j')$.  The
curve is $X_0(N),$ the quotient of~$\mathcal{H}^*$ by the level-$N$ Hecke
subgroup $\Gamma_0(N)$ of~$\Gamma(1),$ comprising all
$\pm\left(\begin{smallmatrix}a&b\\c&d\end{smallmatrix}\right)\in\Gamma(1)$
with $c\equiv0\pmod N$.  This is because the reduced level-$N$ modular
correspondences are bijective with the $\psi(N)$ cosets of $\Gamma_0(N)$.
If $X_0(N)$ is of genus zero, like~$X(1),$ then its function field too will
be singly generated, i.e., will be generated by some univalent function
(`Hauptmodul') $t_N\in\mathbb{C}(j,j')$; and $j,j'$ will be rational
functions of~$t_N$.  An example is the case $N=2$.  Expressions for $j,j'$
in~terms of an appropriate~$t_2$ turn~out to be
\begin{equation}
\label{eq:modular2}
j=\frac{(t_2+16)^3}{t_2},\qquad
j'=\frac{(t_2+256)^3}{t_2^2}.
\end{equation}
These constitute a degree-$2$ {\em rationally parametrized\/} modular
equation, more pleasant and understandable than~(\ref{eq:modular1}).

In \S\S\ref{sec:prelims} and~\ref{sec:parametrized1} we begin by tabulating
parametrized modular equations for~$j$ that are of degrees
$N=2,3,\allowbreak4,\allowbreak5,\allowbreak6,\allowbreak7,\allowbreak8,\allowbreak9,\allowbreak10\allowbreak,1\allowbreak2,\allowbreak13,\allowbreak16,\allowbreak18,25$.
(See Tables \ref{tab:coverings} and~\ref{tab:coverings2}.)  These~$N$ are
the ones for which $X_0(N)$~is of genus zero, and the parameters are
Hauptmoduln~$t_N,$ normalized in accordance with a convention that we
introduce and follow consistently.  The rational formulas $j=j(t_N),$ or
certain more general relations $\Psi_N(t_N,j)=0$ that can be derived when
$X_0(N)$~is of positive genus, have been called `canonical modular
equations'~\cite{Morain95}, since the formulas $j'=j'(t_N)$ can be deduced
from them.  The cases $N=2,3,5,7,13,$ where $N$~is a prime (and
$N-\nobreak1\divides12$) were treated by Klein~\cite{Klein1879}, and his
formulas $j=j(t_N)$ are frequently reproduced.  (See,
e.g.,~\cite[\S4]{Elkies98}.)  Composite values of~$N$ were treated by
Gierster~\cite{Gierster1879}, and magisterially, by
Fricke~\cite[II.~Abschnitt, 4.~Kap.]{Fricke22}, but their findings are not
reproduced in recent references.  Our Tables \ref{tab:coverings}
and~\ref{tab:coverings2} are based on their results, but are more
consistently presented.

Better known than modular equations for~$j$ are modular equations for the
$\lambda$-invariant, or equivalently for the $\alpha$-invariant, which were
extensively investigated in the nineteenth century.  Any elliptic curve
$E/\mathbb{C}$ has a Legendre model $y^2=x(x-1)(x-\nobreak\lambda),$ the
parameter $\lambda\in\mathbb{C}\setminus\{0,1\}$ of~which is the
$\lambda$-invariant.  Like the $j$-invariant, it can be chosen to be a
single-valued function of~$\tau\in\mathcal{H}$.  Its $q$-expansion about
the infinite cusp is $2^4\cdot[q_2-8q_2^2+44q_2^3-192q_2^4+\cdots],$ where
the bracketed series has integer coefficients, and
$q_2:=\sqrt{q}=e^{\pi{\rm i}\tau}$.  The $j$-invariant can be expressed
in~terms of~$\lambda$ by
\begin{equation}
j(\tau)=2^8\,\frac{(\lambda^2-\lambda+1)^3}{\lambda^2(\lambda-1)^2}(\tau) =
2^4\,\frac{(\lambda^2+14\lambda+1)^3}{\lambda(\lambda-1)^4}(2\tau).
\label{eq:mformula}
\end{equation}
Any elliptic curve $E/\mathbb{C}$ also has a quartic Jacobi model
$y^2=(1-\nobreak{x}^2)(1-\nobreak{\alpha}x^2),$ with parameter
$\alpha\in\mathbb{C}\setminus\{0,1\}$.  From its birational equivalence to
the Legendre model one can deduce that $\lambda=\frac{4k}{(1+k)^2},$ where
$\alpha=:k^2$.  In fact, one can choose $\alpha(\tau)=\lambda(2\tau)$; so
the $\alpha$-invariant equals $2^4\cdot[q-8q^2+44q^3-192q^4+\cdots]$.

Many large-$N$ modular equations for~$j$ assume a simple form when written
as an (unparametrized) algebraic relation between $\lambda:=\lambda(\tau)$
and~$\mu:=\lambda(N\tau),$ or equivalently $\alpha:=\alpha(\tau)$
and~$\beta:=\alpha(N\tau)$; or alternatively, between the square roots
$k:=k(\tau)$ and~$l:=k(N\tau)$.  Classical work on modular equations
focused on deriving modular equations of the $k$\textendash\nobreak$l$
type, or the related $u$\textendash\nobreak$v$
type~\cite{Borwein87,Hanna28}.

In \S\ref{sec:parametrized2} we follow a different path: we compute
rationally parametrized modular equations for the Hauptmoduln~$t_M$.  (See
Table~\ref{tab:intermediate}, and for related `factored' modular equations
for~$j,$ see Table~\ref{tab:jfactored}.)  A rational parametrization of the
equation of degree~$N$ at level~$M$ is possible iff $NM,$ as~well as~$M,$
is one of the numbers
$2,3,\allowbreak4,\allowbreak5,\allowbreak6,\allowbreak7,\allowbreak8,\allowbreak9,\allowbreak10\allowbreak,1\allowbreak2,\allowbreak13,\allowbreak16,\allowbreak18,25$.
The connection between parametrized modular equations of this type, and the
$\alpha$\textendash\nobreak$\beta$ or $k$\textendash\nobreak$l$ type, is
not distant.  The invariant~$\lambda$ is a Hauptmodul of the modular
curve~$X(2),$ the quotient of~$\mathcal{H}^*$ by the level-$2$ principal
congruence subgroup $\Gamma(2)$.  The subgroup relation
$\Gamma(2)<\Gamma(1)$ induces an injection of the function field of~$X(1)$
into that of~$X(2),$ so $j$~must be a rational function of~$\lambda$; which
is where (\ref{eq:mformula}) comes from.  But the groups
$\Gamma(2),\Gamma_0(4)<\Gamma(1)$ are conjugated to each other in~${\it
PSL}(2,\mathbb{R})$ by a $2$-isogeny, so the quotients $X(2),X_0(4)$ are
isomorphic.  Their canonical Hauptmoduln $\lambda,t_4,$ the latter to be
defined below, are closely related.  In~fact,
$\lambda(\tau)=[t_4/(t_4+\nobreak16)](\tau/2)$; so
$\alpha=t_4/(t_4+\nobreak16),$ revealing that the Jacobi model is more
closely associated with $\Gamma_0(4)$ than with~$\Gamma(2)$.
In~consequence, any degree-$N$ modular equation of the
$\alpha$\textendash\nobreak$\beta$ type, etc., is equivalent to an
algebraic relation between $t_4(\tau)$ and~$t_4(N\tau)$.

The connection between the traditional approach and ours is illustrated by
Landen's transformation, which is based on a
$\alpha$\textendash\nobreak$\beta$ modular equation of degree~$2,$ namely
Eq.~(\ref{eq:two}), with $\alpha:=\alpha(\tau),$ $\beta:=\alpha(2\tau)$.
Its more familiar $k$\textendash\nobreak$l$ counterpart is
\begin{equation}
l=(1-\nobreak{k}')/(1+\nobreak{k}'),\qquad k':=\sqrt{1-k^2}\,.
\end{equation}
If Eq.~(\ref{eq:two}) is converted to an
$t_4(\tau)$\textendash\nobreak$t_4(2\tau)$ relation, it becomes the
degree-$2$ modular equation at level~$4$ that will be derived
in~\S\ref{sec:parametrized2}, with rational parameter~$t_8$; namely,
\begin{equation}
\label{eq:modular2new}
t_4 = t_8(t_8+8),\qquad t'_4=\frac{t_8^2}{t_8+4},
\end{equation}
where $t_4,t'_4$ signify $t_4(\tau),t_4(2\tau)$.  This is the level-$4$
counterpart of~(\ref{eq:modular2}).  The parameter $t_8$ will be
interpreted as the canonical Hauptmodul of~$X_0(8)$.

In \S\S\ref{sec:modular}, \ref{sec:modular2}, \ref{sec:modular3}
and~\ref{sec:families}, we go beyond Hauptmodul relations, and derive
modular equations for families of elliptic curves.  Along with each
degree-$N$ modular equation for a Hauptmodul of $\Gamma<\Gamma(1),$ e.g.,
(\ref{eq:modular2}) and~(\ref{eq:modular2new}), there is a degree-$N$
modular equation for the elliptic family ${\mathfrak{E}_\Gamma}
\stackrel{\pi}{\to}\Gamma\setminus\mathcal{H}^*$.  It is a functional
equation satisfied by a certain canonical weight-$1$ modular form
for~$\Gamma,$ viewed as a function of the Hauptmodul.  It connects the
periods $\tau_1,\tau_2$ of elliptic curves, i.e., fibres, over related
points on the base, and is really a modular equation for a Gauss--Manin
connection.

To see all this, consider the case $\Gamma=\Gamma(1),$ where the Hauptmodul
is~$j$ and the family is the universal family of elliptic curves.  Since
$j$~is a Fuchsian function of the first kind on $\mathcal{H}\ni\tau,$ it
follows from a theorem on conformal mapping~\cite{Ford51} that any branch
of the multivalued function~$\tau$ on the curve~$X(1)$ can be written
locally as the ratio of two solutions of a second-order differential
equation, with independent variable~$j$.  This equation, used for
constructing the uniformizing variable~$\tau,$ is the classical
Picard--Fuchs equation.  As Stiller~\cite{Stiller88} notes, it is best to
use $\hat J:=1/J=12^3/j$ as the independent variable.  The Picard--Fuchs
equation will then be hypergeometric, with solution space $\hat h_1(\hat
J)(\mathbb{C}\tau+\nobreak\mathbb{C}),$ where $\hat h_1(\hat J)$ is the
Gauss\ hypergeometric function ${}_2F_1(\frac1{12},\frac5{12};1;\hat J)$.
The associated hypergeometric differential operator will be a Gauss--Manin
connection for the universal family.  Near the cusp $\tau={\rm i}\infty,$
where $\hat J=0,$ $\hat h_1(\hat J(\tau))$ is the fourth root of the
Eisenstein sum~$E_4(\tau),$ a weight-$4$ modular form for~$\Gamma(1)$.  So
$\hat{\mathfrak{h}}_1(\tau):=\hat h_1(\hat J(\tau))$ is {\em formally\/} a
weight-$1$ modular form for~$\Gamma(1)$ (formally only, because globally
on~$\mathcal{H},$ it is not single-valued).

By pulling back (i.e., lifting) the Picard--Fuchs equation and its
solutions from $X(1)$ to each genus-zero curve
$X_0(N)\cong\mathbb{P}^1(\mathbb{C})_{t_N},$ we derive a rationally
parametrized degree-$N$ modular equation for~$\hat h_1$.  For instance,
when $N=2$ this equation is {\small
\begin{equation*}
\hat h_1(12^3 t_2/(t_2+16)^3) = 2[(t_2+256)/(t_2+16)]^{-1/4}\,\hat h_1(12^3
t_2^2/(t_2+256)^3),
\end{equation*}
}%
where the arguments of~$\hat h_1$ on the two sides are $\hat J=12^3/j$ and
${\hat J}'=12^3/j',$ written in~terms of the Hauptmodul~$t_2$ of~$X_0(2),$
as in~(\ref{eq:modular2}).  This functional equation is the degree-$2$
modular equation for the universal family of elliptic curves parametrized
by~$j$ (or~$\hat J$).  It relates the periods of two elliptic curves with
$j$-invariants related by~(\ref{eq:modular1}), i.e.,
by~(\ref{eq:modular2}).  Because $\hat
h_1(\cdot)={}_2F_1(\frac1{12},\frac5{12};1;\cdot),$ each such modular
equation is an {\em algebraic hypergeometric transformation}: it relates
a~${}_2F_1$ to another~${}_2F_1$ with an algebraically transformed
argument.

At each level $M>1$ with a genus-zero $X_0(M),$ there is a rational
elliptic family
$\mathfrak{E}_M\stackrel{\pi}{\to}X_0(M)\cong\mathbb{P}^1(\mathbb{C})_{t_M}$.
By pulling back along $X_0(M)\to X(1),$ one can derive a Gauss--Manin
connection for~it.  We make this concrete by pulling back the function
$\hat h_1=\hat h_1(\hat J)$ to a function $h_M=h_M(t_M),$ and working~out
eta product and $q$-series representations for the weight-$1$ modular form
$\mathfrak{h}_M(\tau):=h_M(t_M(\tau))$ on~$\mathcal{H}$.  (See Tables
\ref{tab:xNhN} and~\ref{tab:windup}.)  The Picard--Fuchs equation for the
level-$M$ family has solution space
$h_M(t_M)(\mathbb{C}\tau+\nobreak\mathbb{C})$.  Our key theorem on {\em
further\/} pullings back of the modular form $h_M(t_M(\tau)),$
Theorem~\ref{thm:key}, efficiently produces a degree-$N$ modular equation
for~$h_M$ if $NM,$ as~well as~$M,$ is one of
$2,3,\allowbreak4,\allowbreak5,\allowbreak6,\allowbreak7,\allowbreak8,\allowbreak9,
\allowbreak10,\allowbreak12,\allowbreak13,\allowbreak16,\allowbreak18,25$.
These equations are listed in Table~\ref{tab:keymodular}, the most
important single table of this article (though its corollary
Table~\ref{tab:ramanujan} is more immediately understandable).  They may be
viewed as (i)~algebraic transformations of Picard--Fuchs equations, or
(ii)~transformations of the special functions that satisfy them (i.e., the
functions $h_M=h_M(t_M)$), or (iii)~transformations of certain integer
sequences defined by recurrences (i.e., the coefficient sequences of the
functions~$h_M,$ expanded in~$t_M$).

The cases $M=2,3,4$ are especially interesting, partly because in these
cases the multiplier system of the modular form
$\mathfrak{h}_M(\tau)=h_M(t_M(\tau))$ is nontrivial: only if $M\geq5$ is it
given by a Dirichlet character mod~$M$.  When $M=2,3,4,$ the Picard--Fuchs
equation for~$\Gamma_0(M)$ is of {\em hypergeometric type\/}, with three
regular singular points on~$X_0(M)\cong\mathbb{P}^1(\mathbb{C})_{t_M}$;
namely, the fixed points of~$\Gamma_0(M)$.  For instance, the level-$4$
function $h_4=h_4(t_4),$ pulled back from~$\hat h_1=\hat h_1(\hat J),$ is
${}_2F_1(\frac12,\frac12;1;-t_4/16)$.  The degree-$2$ modular equation for
the elliptic family parametrized by~$t_4$ is
\begin{equation}
\label{eq:preLanden}
h_4(t_8(t_8+8)) = 2(t_8+4)^{-1/2} h_4(t_8^2/(t_8+4)),
\end{equation}
where the arguments of~$h_4$ on the two sides are $t_4$ and~$t_4'$ written
in~terms of a Hauptmodul~$t_8$ for~$\Gamma_0(8),$ as
in~(\ref{eq:modular2new}).  This is another algebraic hypergeometric
transformation, as are all our modular equations at levels $M=2,3,4$.

If $M=5,6,7,8,9,$ then the Picard--Fuchs equation is of {\em Heun type\/},
with four singular points on~$X_0(M)\cong\mathbb{P}^1(\mathbb{C})_{t_M}$;
as before, the fixed points of~$\Gamma_0(M)$.  Hence, its solutions can be
expressed in terms of the canonical `local Heun'
function~$\Hl$~\cite{Ronveaux95}.  (For basic facts on equations of
hypergeometric and Heun type, see the Appendix.)  Our elliptic-family
modular equations on levels $5,6,7,8,9$ are thus {\em algebraic Heun
transformations\/}: functional equations satisfied by~${\it Hl}$.  By
expanding ${\it Hl}$ in its argument and focusing on the resulting
coefficient sequence, one can view them as identities satisfied by integer
sequences defined by three-term recurrences.  The functional equations of
Proposition~\ref{prop:franel}, which are satisfied by the generating
function for the Franel numbers $\sum_{k=0}^n \binom{n}{k}^3,$ $n\geq0,$
are an example.  The theory of such Heun transformations, whether viewed
combinatorially or not, is entirely undeveloped.

Our modular equations at levels $4,3,2$ fit into Ramanujan's theories of
elliptic integrals in signatures $r=2,3,4,$ respectively.  In special
function terms, this is because his complete elliptic integral in
signature~$r,$ denoted $\K_r(\alpha_r)$ here, is proportional to
${}_2F_1(\frac1r,1-\nobreak\frac1r;1;\alpha_r)$.  By employing Pfaff's
transformation of~${}_2F_1,$ given in the Appendix as~(\ref{eq:Pfaff}), one
can convert identities involving $h_4,h_3,h_2$ to ones involving
$\K_2,\K_3,\K_4$.  For instance, the modular
equation~(\ref{eq:preLanden}) becomes
\begin{equation}
\label{eq:K2modular}
\K_2(4p/(p+1)^2) = (1+p)\,\K_2(p^2).
\end{equation}
This is a parametric form of Landen's transformation, since $r=2$ is the
classical base and $\K_2$ equals~$\K,$ the traditional complete elliptic
integral.  It is equivalent to~(\ref{eq:duo}).  Like~$h_4,$ ${\K}_2$~can be
viewed as defining a weight-$1$ modular form for~$\Gamma_0(4),$ and the
modular equation~(\ref{eq:K2modular}) expresses $\K_2$ at a general point
$\tau\in\mathcal{H}$ in~terms of its value at $\tau'=2\tau$.

In \S\ref{sec:ramanujan}, by systematically converting our modular
equations for the elliptic families parametrized by $t_4,t_3,t_2,$ we
compute parametrized modular equations of degrees $N=2,3,4$ in Ramanujan's
theory of signature~$2$; of~degrees $N=2,3,\allowbreak4,{\underline6}$ in
that of signature~$3$; and of degrees
$N=2,3,\allowbreak4,\allowbreak{\underline5},\allowbreak{\underline6},\allowbreak{\underline8},{\underline9}$
in that of signature~$4$.  (See Table~\ref{tab:ramanujan}; the underlined
ones are new to the literature.)  We~also derive equations relating
$\K_2,\K_3,\K_4,$ arising from commensurability of
the subgroups $\Gamma_0(4),\allowbreak\Gamma_0(3),\allowbreak\Gamma_0(2)$
of~$\Gamma(1)$.  Finally, in~\S\ref{sec:final} we discuss the modular
underpinnings of Ramanujan's theories, and explain why his underdeveloped
but very interesting theory of signature~$6$ is essentially nonclassical.

\section{Preliminaries: Fixed Points}
\label{sec:prelims}

The modular curve $X_0(N)=\Gamma_0(N)\setminus\mathcal{H}^*$ classifies
each cyclic $N$-isogeny $\phi:E\to E/C$ up~to isomorphism, where an
isomorphism $(E_1,C_1)\cong(E_2,C_2)$ is an isomorphism of elliptic
curves~$E_1,E_2$ that takes $C_1$ to~$C_2$.  Since $X_0(N)$~is defined
over~$\mathbb{Q},$ it classifies up~to isomorphism over~$\overline{\mathbb
Q}$.  As~noted, its function field is $\mathbb{C}(j,j')$ where
$j'(\tau)=j(N\tau)$.  The Fricke involution $W_N:(E,C)\to(E/C,E_N/C),$
which on~$\mathcal{H}$ is the map $\tau\mapsto-1/N\tau,$
interchanges~$j,j'$.

The cusps of~$X_0(N)$ are as~follows~\cite{Ogg73}.  The set of cusps
of~$\mathcal{H}^*\ni\tau,$ i.e., $\mathbb{P}^1(\mathbb{Q}),$ is partitioned
into equivalence classes under~$\Gamma_0(N)$.  A~system of representatives,
i.e., a choice of one from each, may be taken to comprise certain fractions
$\tau=\frac{a}d$ for each~$d\divides{N},$ with $1\le a<d$ and~$(a,d)=1$.
Here $a$~is reduced modulo $(d,N/d)$ while remaining coprime to~$d,$ so
there are $\varphi((d,N/d))$ possible values of~$a,$ and hence
$\varphi((d,N/d))$ inequivalent cusps in~$\mathbb{P}^1(\mathbb{Q})$
associated to~$d$.  (In this statement $(\cdot,\cdot)$ is the greatest
common divisor function, and $\varphi(\cdot)$~is the Euler totient
function.)  So in~all, $X_0(N)$ has
\begin{equation}
\sigma_\infty(N):=\sum_{d\divides{N}}\varphi((d,N/d))
\end{equation}
cusps, which are involuted by~$W_N$.  The {\em rational\/}
cusps~$\frac{a}d$ are those for which $\varphi((d,N/d))=1,$ i.e, the ones
with $(d,N/d)=1$ or~$2$.

The covering $j:X_0(N)\to\mathbb{P}^1(\mathbb{C})\cong X(1)$ is
$\psi(N)$-sheeted, since the index $[\Gamma(1):\Gamma_0(N)]$
equals~$\psi(N)$.  It is ramified only above the cusp $j=\infty$ and the
elliptic fixed points $j=0,12^3,$ corresponding to equianharmonic and
lemniscatic elliptic curves; i.e., only above the three vertices $\tau={\rm
i}\infty,\allowbreak{\zeta_3:=e^{2\pi{\rm i}/3}},{\rm i}$ of the
fundamental half-domain of~$\Gamma(1)$.  The fibre above ${j=\infty}$
includes (the equivalence class~of) each cusp $\tau=\frac{a}d$ with
multiplicity equal to its width $e_{d,N}:=N/d(d,N/d)$.  To~indicate that a
cusp is an equivalence class, the notation $\bigl[\frac{a}d\bigr]$
or~$\bigl[\frac{a}d\bigr]_N$ will be used as~appropriate.

\begin{table}
\caption{Basic data on each genus-$0$ group
$\Gamma_0(N)<\Gamma(1)$.}
\begin{center}
{\small
\begin{tabular}{cccllcc}
\hline\noalign{\smallskip}
$N$ & $\psi(N)$ & $\sigma_\infty(N)$ & \hfil Cusps ($\tau$ values) & \hfil Cusp widths & $\epsilon_2(N)$ & $\epsilon_3(N)$\\
\noalign{\smallskip}\hline\noalign{\smallskip}
$2$ & $3$ & $2$ & $\frac11;\frac12$ & $2;1$ & $1$ & $0$ \\ 
$3$ & $4$ & $2$ & $\frac11;\frac13$ & $3;1$ & $0$ & $1$ \\ 
$4$ & $6$ & $3$ & $\frac11;\frac12,\frac14$ & $4;1,1$ & $0$ & $0$ \\ 
$5$ & $6$ & $2$ & $\frac11;\frac15$ & $5;1$ & $2$ & $0$ \\ 
$6$ & $12$ & $4$ & $\frac11;\frac12;\frac13;\frac16$ & $6;3;2;1$ & $0$ & $0$ \\ 
$7$ & $8$ & $2$ & $\frac11;\frac17$  & $7;1$ & $0$ & $2$ \\ 
$8$ & $12$ & $4$ & $\frac11;\frac12;\frac14,\frac18$ & $8;2;1,1$ & $0$ & $0$ \\ 
$9$ & $12$ & $4$ & $\frac11;\frac13,\frac23,\frac19$ & $9;1,1,1$ & $0$ & $0$ \\ 
$10$ & $18$ & $4$ & $\frac11;\frac12;\frac15;\frac1{10}$  & $10;5;2;1$ & $2$ & $0$ \\ 
$12$ & $24$ & $6$ & $\frac11;\frac13;\frac12,\frac14;\frac16,\frac1{12}$ & $12;4;3,3;1,1$ & $0$ & $0$ \\ 
$13$ & $14$ & $2$ & $\frac11;\frac1{13}$ & $13;1$ & $2$ & $2$ \\ 
$16$ & $24$ & $6$ & $\frac11;\frac12;\frac14,\frac34,\frac18,\frac1{16}$ & $16;4;1,1,1,1$ & $0$ & $0$ \\ 
$18$ & $36$ & $8$ & $\frac11;\frac12;\frac13,\frac23,\frac19;$ & $18;9;2,2,2;$ & &  \\ 
     &      &     & ${}\qquad\quad\;\;\frac16,\frac56,\frac1{18}$ & ${}\qquad\quad\;\;1,1,1$ & $0$ & $0$ \\ 
$25$ & $30$ & $6$ & $\frac11;\frac15,\frac25,\frac35,\frac45,\frac1{25}$ & $25;1,1,1,1,1$ & $2$ & $0$ \\ 
\noalign{\smallskip}\hline
\end{tabular}
}%
\end{center}
\label{tab:basic}
\end{table}

The fibre above $j=0$ (resp.~$12^3$) includes $\epsilon_3(N)$ cubic
elliptic points (resp.\ $\epsilon_2(N)$ quadratic ones), each with unit
multiplicity; other points on the fibre appear with cubic (resp.\
quadratic) multiplicity~\cite{Schoeneberg74}.  Here
\begin{subequations}
\begin{align}
\epsilon_2(N)&:=
\begin{cases}
\prod_{p\divides{N},\text{ $p$ prime}}\left(1+\left(\frac{-1}{p}\right)\right), &
4\notdivides{N}, \\ 0, & 4\divides{N},
\end{cases}\\*
\epsilon_3(N)&:=
\begin{cases}
\prod_{p\divides{N},\text{ $p$ prime}}\left(1+\left(\frac{-3}{p}\right)\right), &
9\notdivides{N}, \\ 0, & 9\divides{N},
\end{cases}
\end{align}
\end{subequations}
with $\left(\frac{\cdot}{\cdot}\right)$ the Legendre symbol.  Any elliptic
curve $\mathbb{C}/(\mathbb{Z}\tau_1\oplus\mathbb{Z}\tau_2)$ that is
lemniscatic has $\epsilon_2(N)$~non-isomorphic self-$N$-isogenies,
performed by complex multiplication of the period lattice by a Gaussian
integer; so $\epsilon_2(N)$~simply counts the ways of representing~$N$ as
the sum of two squares.  $\epsilon_3(N)$~has a similar interpretation,
in~terms of Eisenstein integers.  It~follows from the Hurwitz formula that
$X_0(N)$~has genus
\begin{equation}
\label{eq:genus}
g = 1+\frac{\psi(N)}{12} - \frac{\sigma_\infty(N)}2 
-\frac{\epsilon_2(N)}{4}-\frac{\epsilon_3(N)}{3}.
\end{equation}
In Table~\ref{tab:basic}, the basic data (number of fixed points, cusp
locations, etc.)\ are listed for each of the 14~genus-zero curves $X_0(N)$.

\section{Hauptmoduln and Parametrized Modular Equations}
\label{sec:parametrized1}

To make the covering $j:X_0(N)\to\mathbb{P}^1(\mathbb{C})\cong X(1)$ more
concrete, one needs (i)~an explicit formula for a Hauptmodul~$t_N$
for~$\Gamma_0(N),$ so that $X_0(N)$ like~$X(1)$ can be identified with
$\mathbb{P}^1(\mathbb{C}),$ and (ii)~an expression for the covering map, as
a degree-$\psi(N)$ rational function of~$t_N$.  For most of the above
14~values of~$N,$ items (i) and~(ii) were worked~out by Klein, Gierster and
Fricke.  But because the classical derivation was somewhat unsystematic,
ours is {\em de~novo\/}.

\subsection{Canonical Hauptmoduln}
\label{subsec:parametrized1a}

For each~$N,$ we can specify the Hauptmodul~$t_N$ uniquely by requiring
that (I) $t_N$~have a simple zero at the cusp $\bigl[\frac1N\bigr]\ni{\rm
i}\infty$ and a simple pole at the cusp $\bigl[\frac11\bigr]\ni0,$ and
(II)~the function $(t_N|_{W_N})(\tau):=t_N(-1/N\tau)$ have a Fourier
expansion on~$\mathcal{H}^*$ that begins $1\cdot q^{-1}+O(q^0),$ where
$q:=e^{2\pi{\rm i}\tau}$.  Since the Fricke involution
$W_N:\tau\mapsto-1/N\tau$ interchanges $\bigl[\frac1N\bigr]$
and~$\bigl[\frac11\bigr],$ the product $t_N(\tau)t_N(-1/N\tau)$ necessarily
equals some constant function of~$\tau,$ to be denoted~$\kappa_N,$ and
condition~(II) fixes this constant.  Imposing this normalization condition
will simplify the modular equations to be deduced, e.g., by forcing many
polynomial factors to be monic.

A $q$-product, such as the one used to define the Dedekind eta function, is
the natural way of defining each~$t_N$.  Though Hauptmoduln of general
genus-zero congruence subgroups cannot be expressed in~terms of the eta
function, it turns~out to be possible for each
genus-zero~$\Gamma_0(N),$ as we explain.  Recall that
$\eta(\tau):=q^{1/24}\prod_{n=1}^\infty (1-q^n)$ satisfies
$\eta(\tau+1)=\zeta_{24}\eta(\tau)$ and $\eta(-1/\tau)=(-{\rm
i}\tau)^{1/2}\eta(\tau)$.  One first defines an unnormalized
Hauptmodul~$\hat t_N$ by
\begin{equation}
  \hat t_N = \prod_{\delta\divides{N}} \eta(\delta\tau)^{r_\delta},
\end{equation}
where $\delta\mapsto r_\delta\in\mathbb{Z}$ is an appropriate `generalized
permutation'.  An obvious requirement is that the Hauptmodul have the
correct order of vanishing at each cusp~\cite{Ligozat75}.  Since
$\eta(\delta\tau)$ has order $\frac1{24}(\delta,d)^2/\delta$ at any cusp
$\tau=\frac{a}d$ on~$\mathcal{H}^*,$ this function $\hat t_N(\tau)$ has
order $\frac1{24}\sum_{\delta\divides{N}}r_\delta (\delta,d)^2/\delta$
there.  To~compute its order at $\bigl[\frac{a}d\bigr]\in X_0(N),$ one must
multiply by the cusp width $e_{d,N}=N/d(d,N/d)$.  The resulting order of
vanishing should be $+1$ at~$\bigl[\frac{a}d\bigr]=\frac1N,$
$-1$~at~$\bigl[\frac11\bigr],$ and $0$~at each of the other
$\sigma_\infty(N)-2$ cusps.

The generalized permutation $\delta\mapsto r_\delta$ must also be such that
$\hat t_N$~is a single-valued function on~$X_0(N)$.  A~sufficient condition
for this was given by Newman~\cite{Newman58}.  If it is the case that
$\sum_{\delta\divides{N}}r_\delta\equiv0\pmod{24},$
$\sum_{\delta\divides{N}}\delta r_\delta\equiv0\pmod{24},$ and
$\sum_{\delta\divides{N}} (N/\delta)r_\delta\equiv0\pmod{24},$ and also
$\prod_{\delta\divides{N}}\delta^{r_\delta}={\kappa_N}^2$ for some natural
number~$\kappa_N,$ then $\hat t_N$~will be single-valued on~$X_0(N)$.
Moreover, $\hat t_N|_{W_N}$ will equal
${\kappa_N}^{-1}\prod_{\delta\divides{N}} \eta(\delta\tau)^{r_{N/\delta}}$.
It~follows that to satisfy condition~(II), one should let
$t_N:=\kappa_N\cdot\hat t_N$.  If that choice is made,
$t_N|_{W_N}=\prod_{\delta\divides{N}} \eta(\delta\tau)^{r_{N/\delta}}$ will
have a Fourier expansion beginning $1\cdot q^{-1}+O(q^0)$.

For each of the 14~values of~$N,$ a unique map $\delta\mapsto r_\delta$
satisfying the preceding conditions can be found by inspection.  The
resulting normalized eta-product expressions,
$t_N:=\kappa_N\cdot\prod_{\delta\divides{N}} \eta(\delta\tau)^{r_\delta},$
are listed in the second column of Table~\ref{tab:hauptmoduln}.  Fine's
compact notation $[\delta]$ for the function~$\eta(\delta\tau)$
on~$\mathcal{H}^*\ni\tau$ is used.  These eta products were worked~out by
Fricke, with the exception of those for $N=4,16,18,$ and also with the
exception of the one for $N=2$.  (In~an unfortunate confusion that can be
traced to Klein~\cite[p.~143]{Klein1879}, Fricke's analogue of~$t_2(\tau)$
is proportional to~$t_2(2\tau)$.)  It is well known
that~\cite[Ch.~4]{Apostol90}
\begin{equation}
t_N=N^{12/(N-1)}\cdot[N]^{24/(N-1)}/\,[1]^{24/(N-1)},\qquad\text{if}\quad
N-\nobreak1\divides24.
\end{equation}
But in the six cases when $N-1\notdivides{24},$ the prefactors~$\kappa_N$
seem not to have been tabulated before.  Several of those given by
Klein--Gierster--Fricke differ from ours, since they did not consistently
impose condition~(II)\null.  Shih~\cite{Shih74} reproduces and uses
their~$\kappa_N$'s, each of which is less than ours by a factor equal to a
rational square.  Fortunately, for Shih's application inconsistently
scaled~$\kappa_N$'s are adequate.

For each~$N,$ if $\mathfrak{c}_0,\mathfrak{c}_\infty$ are distinct rational
cusps for the genus-zero subgroup $\Gamma_0(N),$ there is a
Hauptmodul~$\tilde t_N$ on~$X_0(N)$ with a simple zero at~$\mathfrak{c}_0$
and a simple pole at~${c}_\infty,$ which is a rational (over~$\mathbb{Q}$)
degree-$1$ function of the canonical Hauptmodul~$t_N$.  In
Table~\ref{tab:alt} the alternative Hauptmoduln~$ \tilde t_N$ for which
$\mathfrak{c}_0=\bigl[\frac1N\bigr]\ni{\rm i}\infty$ are listed.  They were
extracted from a long list of Hauptmoduln of genus-zero groups obtained by
Ford et~al.~\cite{Ford94} by manipulating $q$-series, but can also be
verified individually by examining behavior at each cusp.  The eta-product
representation $2^{4}\cdot[1]^{8}[4]^{16}/\,[2]^{24}$ of the alternative
Hauptmodul $\alpha=t_4/(t_4+16)$ for $X_0(4)$ should be noted.  Combined
with the known expression $2^{4}\cdot[\frac12]^{8}[2]^{16}/\,[1]^{24}$ for
the $\lambda$-invariant, it confirms that
$\lambda(\tau)=[t_4/(t_4+\nobreak16)](\tau/2)=\alpha(\tau/2),$ as mentioned
in~\S\ref{subsec:12}.

\begin{table}
\caption{The canonical Hauptmodul~$t_N=\kappa_N\cdot\hat t_N$
  for~$\Gamma_0(N),$ and the fixed points of~$\Gamma_0(N)$ on~$X_0(N)$.
  (The polynomial $p_{25}(t)$ equals $t^4+5t^3+15t^2+25t+25$.)}
\begin{center}
{\footnotesize
\begin{tabular}{clll}
\hline\noalign{\smallskip}
$N$ & $t_N(\tau)=\kappa_N\cdot\hat t_N(\tau)$ & Cusps ($t_N$ values) & Elliptic
points ($t_N$ values)\\
\noalign{\smallskip}\hline\noalign{\smallskip}
$2$ & $2^{12}\cdot[2]^{24}/\,[1]^{24}$ & $\infty;0$ & $-64$ [quadratic]\\
$3$ & $3^{6}\cdot[3]^{12}/\,[1]^{12}$ & $\infty;0$ & $-27$ [cubic] \\
$4$ & $2^{8}\cdot[4]^{8}/\,[1]^{8}$ & $\infty;-16,0$ & \\
$5$ & $5^{3}\cdot[5]^{6}/\,[1]^{6}$ & $\infty;0$ & $-11\pm2\sqrt{-1}$
[quadratic] \\
$6$ & $2^{3} 3^{2}\cdot [2][6]^{5}/\,[1]^{5}[3]$ & $\infty;-8;-9;0$ & \\
$7$ & $7^{2}\cdot[7]^{4}/\,[1]^{4}$ & $\infty;0$ & $\frac{-13\pm3\sqrt{-3}}2$ [cubic] \\
$8$ & $2^{5}\cdot[2]^{2}[8]^{4}/\,[1]^{4}[4]^{2}$ & $\infty;-4;-8,0$ & \\
$9$ & $3^{3}\cdot[9]^{3}/\,[1]^{3}$ &
$\infty;\frac{-9\pm3\sqrt{-3}}2,0$ & \\
$10$ & $2^{2}5\cdot[2]\,[10]^{3}/\,[1]^{3}[5]$ & $\infty;-4;-5;0$ &
$-4\pm2\sqrt{-1}$ [quadratic]\\
$12$ & $2^{2}3\cdot[2]^2[3]\,[12]^3/\,[1]^{3}[4]\,[6]^2$ &
$\infty;-3;-2,-4;-6,0$ & \\
$13$ & $13\cdot[13]^{2}/\,[1]^{2}$ & $\infty;0$ & $-3\pm2\sqrt{-1}$ [quadratic],\\
 &  &  & $\frac{-5\pm3\sqrt{-3}}2$ [cubic] \\
$16$ & $2^{3}\cdot[2]\,[16]^{2}/\,[1]^{2}[8]$ & $\infty;-2;$ & \\
     & & $-2\pm2\sqrt{-1},-4,0$ & \\
$18$ & $2\cdot3\cdot[2]\,[3]\,[18]^2/\,[1]^{2}[6]\,[9]$ & $\infty;-2;\frac{-3\pm\sqrt{-3}}2,-3;$ & \\
    &  & $-3\pm\sqrt{-3},0$ & \\
$25$ & $5\cdot[25]\,/\,[1]$ & $\infty$; roots of $p_{25}(t),$ $0$ & $-1\pm2\sqrt{-1}$ [quadratic] \\
\noalign{\smallskip}\hline
\end{tabular}
}%
\end{center}
\label{tab:hauptmoduln}
\end{table}

\begin{table}
\caption{Alternative Hauptmoduln $\tilde t_N$ for $\Gamma_0(N)$.}
\begin{center}
{\small
\begin{tabular}{clcl}
\hline\noalign{\smallskip}
$N$ & \hfil$\tilde t_N$ & pole ($\tau$~class) & eta product representation \\
\noalign{\smallskip}\hline\noalign{\vskip1pt}\hline\noalign{\vskip1.5pt}
$4$ & $t_4/(t_4+16)$ & $\bigl[\frac12\bigr]$ & $2^{4}\cdot[1]^{8}[4]^{16}/\,[2]^{24}$ \\
\noalign{\vskip 1pt}\hline\noalign{\vskip 1pt}
$6$ & $t_6/(t_6+8)$ & $\bigl[\frac12\bigr]$ & $3^{2}\cdot[1]^{4}[6]^{8}/\,[2]^{8}[3]^{4}$ \\
$6$ & $t_6/(t_6+9)$ & $\bigl[\frac13\bigr]$ & $2^{3}\cdot[1]^{3}[6]^{9}/\,[2]^{3}[3]^{9}$ \\
\noalign{\vskip 1pt}\hline\noalign{\vskip 1pt}
$8$ & $t_8/(t_8+4)$ & $\bigl[\frac12\bigr]$ & $2^{3}\cdot[1]^{4}[4]^{2}[8]^{4}/\,[2]^{10}$ \\
$8$ & $t_8/(t_8+8)$ & $\bigl[\frac14\bigr]$ & $2^{2}\cdot[2]^{4}[8]^{8}/\,[4]^{12}$ \\
\noalign{\vskip 1pt}\hline\noalign{\vskip 1pt}
$10$ & $t_{10}/(t_{10}+4)$ & $\bigl[\frac12\bigr]$ & $5\cdot[1]^{2}[10]^{4}/\,[2]^{4}[5]^{2}$ \\
$10$ & $t_{10}/(t_{10}+5)$ & $\bigl[\frac15\bigr]$ & $2^{2}\cdot[1][10]^{5}/\,[2][5]^{5}$ \\
\noalign{\vskip 1pt}\hline\noalign{\vskip 1pt}
$12$ & $t_{12}/(t_{12}+2)$ & $\bigl[\frac12\bigr]$ & $2\cdot3\cdot[1]^{3}[4]^{2}[6][12]^{2}/\,[2]^{7}[3]$ \\
$12$ & $t_{12}/(t_{12}+3)$ & $\bigl[\frac13\bigr]$ & $2^{2}\cdot[1][12]^{3}/\,[3]^3[4]$ \\
$12$ & $t_{12}/(t_{12}+4)$ & $\bigl[\frac14\bigr]$ & $3\cdot[2]^{2}[12]^{4}/\,[4]^{4}[6]^{2}$ \\
$12$ & $t_{12}/(t_{12}+6)$ & $\bigl[\frac16\bigr]$ & $2\cdot[2]^{3}[3]^{3}[12]^{6}/\,[1][4]^{2}[6]^{9}$ \\
\noalign{\vskip 1pt}\hline\noalign{\vskip 1pt}
$16$ & $t_{16}/(t_{16}+2)$ & $\bigl[\frac12\bigr]$ & $2^{2}\cdot[1]^{2}[4]^{2}[16]^{2}/\,[2]^{5}[8]$ \\
$16$ & $t_{16}/(t_{16}+4)$ & $\bigl[\frac18\bigr]$ & $2\cdot[4]^{2}[16]^{4}/\,[8]^{6}$ \\
\noalign{\vskip 1pt}\hline\noalign{\vskip 1pt}
$18$ & $t_{18}/(t_{18}+2)$ & $\bigl[\frac12\bigr]$ & $3\cdot[1][18]^{2}/\,[2]^{2}[9]$ \\
$18$ & $t_{18}/(t_{18}+3)$ & $\bigl[\frac19\bigr]$ & $2\cdot[3][18]^{3}/\,[6][9]^{3}$ \\
\noalign{\smallskip}\hline
\end{tabular}
}%
\end{center}
\label{tab:alt}
\end{table}

One can find $q$-series for most of the Hauptmoduln of Tables
\ref{tab:hauptmoduln} and~\ref{tab:alt} in Sloane's {\em
Encyclopedia\/}~\cite{Sloane2005}.  In a neighborhood of the infinite cusp,
i.e., of the point ${q=0},$ each can be written as ${\rm
prefactor}\cdot\left[q+\sum_{n=2}^\infty a_n q^n\right],$ where the
$a_n$~are integers.  Thus $t_4=2^8\cdot[q+8q^2+44q^3+192q^4+\cdots],$ which
is consistent with the $q$-expansion of $\alpha=t_4/(t_4+16)$ given
in~\S\ref{subsec:12}.  It~is worth remarking that with the exception of
prefactors, the $q$-series for $t_4$ and $\alpha=t_4/(t_4+16)$ are related
by~$q\mapsto-q,$ i.e., $\tau\mapsto\tau+\frac12,$ as are those for $t_8$
and~$t_8/(t_8+4),$ and those for $t_{16}$ and~$t_{16}/(t_{16}+2)$.  The
$q$-series coefficients~$a_n$ for the Hauptmoduln of
Table~\ref{tab:hauptmoduln}, unlike those of Table~\ref{tab:alt}, are
always non-negative.

\begin{remark*}
Many of these Hauptmoduln, or their $q$-series, have cropped~up in the
literature.  For instance, the alternative Hauptmodul
$[3][18]^3/\,[6][9]^3$ for $\Gamma_0(18),$ i.e., ${t_{18}/2(t_{18}+3)},$ was
expanded in an interesting continued fraction by both Ramanujan and
Selberg.  (Duke~\cite[(9.13)]{Duke2005} gives this continued fraction, with
$q\mapsto q^{1/3}$.)  Also, the alternative Hauptmodul
$[4]^{2}[16]^{4}/\,[8]^{6}$ for~$\Gamma_0(16),$ i.e., $t_{16}/2(t_{16}+4),$
is the so-called $\varepsilon$\nobreakdash-invariant, in which Weierstrass
developed the nome~$q$ as a power series.  That is,
\begin{equation}
  q=q(\varepsilon) = \varepsilon+\sum_{k=1}^\infty \delta_k\,\varepsilon^{4k+1}.
\end{equation}
The problem of efficiently calculating the sequence of positive integers
$\{\delta_k\}_{k=1}^\infty$ by reverting the $q$\nobreakdash-series
for~$\varepsilon$ has generated a significant literature~\cite{Ferguson75}.
\end{remark*}

\begin{table}
\caption{The degree-$\psi(N)$ covering map from\hfil\break
  $X_0(N)\cong\mathbb{P}^1(\mathbb{C})_{t_N}$ 
  to~$X(1)\cong\mathbb{P}^1(\mathbb{C})_{j}$.}
\begin{center}
{\small
\begin{tabular}{cl}
\hline\noalign{\smallskip}
$N$ & \hfil $j(\tau)$ as a function of~$t_N(\tau)$ \\
\noalign{\smallskip}\hline\noalign{\smallskip}
$2$ & $\frac{(t+16)^3}t$ \\
    & $\quad = 12^3 + \frac{(t+64)(t-8)^2}t$ \\
$3$ & $\frac{(t+27)(t+3)^3}t$   \\
    & $\quad = 12^3 + \frac{(t^2+18t-27)^2}{t} $ \\
$4$ & $\frac{(t^2+16t+16)^3}{t(t+16)}$  \\
    & $\quad = 12^3 + \frac{(t+8)^2(t^2+16t-8)^2}{t(t+16)}$ \\
$5$ & $\frac{(t^2+10t+5)^3}{t}$  \\
    & $\quad = 12^3 + \frac{(t^2+22t+125)(t^2+4t-1)^2}{t}$ \\
$6$ & $\frac{(t+6)^3(t^3+18t^2+84t+24)^3}{t(t+8)^3(t+9)^2}$  \\
    & $\quad = 12^3 + \frac{(t^2+12t+24)^2(t^4+24t^3+192t^2+504t-72)^2}{t(t+8)^3(t+9)^2}$ \\
$7$ & $\frac{(t^2+13t+49)(t^2+5t+1)^3}{t}$  \\
    & $\quad = 12^3 + \frac{(t^4+14t^3+63t^2+70t-7)^2}{t}$ \\
$8$ & $\frac{(t^4+16t^3+80t^2+128t+16)^3}{t(t+4)^2(t+8)}$  \\
& $\quad = 12^3 + \frac{(t^2+8t+8)^2(t^4+16t^3+80t^2+128t-8)^2}{t(t+4)^2(t+8)}$ \\
$9$ & $\frac{(t+3)^3(t^3+9t^2+27t+3)^3}{t(t^2+9t+27)}$  \\
    & $\quad = 12^3 + \frac{(t^6+18t^5+135t^4+504t^3+891t^2+486t-27)^2}{t(t^2+9t+27)} $ \\
$10$ & $\frac{(t^6+20t^5+160t^4+640t^3+1280t^2+1040t+80)^3}{t(t+4)^5(t+5)^2}$  \\
    & $\quad = 12^3 + \frac{(t^2+8t+20)(t^2+6t+4)^2(t^2+6t+10)^2(t^4+14t^3+66t^2+104t-4)^2}{t(t+4)^5(t+5)^2}$ \\
$12$ & $\frac{(t^2+6t+6)^3(t^6+18t^5+126t^4+432t^3+732t^2+504t+24)^3}{t(t+2)^3(t+3)^4(t+4)^3(t+6)}$  \\
    & {\quad\scriptsize$ = 12^3 + \frac{(t^4+12t^3+48t^2+72t+24)^2(t^8+24t^7+240t^6+1296t^5+4080t^4+7488t^3+7416t^2+3024t-72)^2}{t(t+2)^3(t+3)^4(t+4)^3(t+6)}$} \\
$13$ & $\frac{(t^2+5t+13)(t^4+7t^3+20t^2+19t+1)^3}{t}$  \\
    & $\quad = 12^3 + \frac{(t^2+6t+13)(t^6+10t^5+46t^4+108t^3+122t^2+38t-1)^2}{t} $ \\
$16$ & $\frac{(t^8+16t^7+112t^6+448t^5+1104t^4+1664t^3+1408t^2+512t+16)^3}{t(t+2)^4(t+4)(t^2+4t+8)}$  \\
    & {\quad\scriptsize$ = 12^3 + \frac{(t^4+8t^3+24t^2+32t+8)^2(t^8+16t^7+112t^6+448t^5+1104t^4+1664t^3+1408t^2+512t-8)^2}{t(t+2)^4(t+4)(t^2+4t+8)}$} \\
$18$ & {\scriptsize$\frac{(t^3+6t^2+12t+6)^3(t^9+18t^8+144t^7+666t^6+1944t^5+3672t^4+4404t^3+3096t^2+1008t+24)^3}{t(t+2)^9(t+3)^2(t^2+3t+3)^2(t^2+6t+12)}$}  \\
    & \quad{$ = 12^3 +\, {\rm etc.}$} \\ 
$25$ & $\frac{(t^{10}+10t^9+55t^8+200t^7+525t^6+1010t^5+1425t^4+1400t^3+875t^2+250t+5)^3}{t(t^4+5t^3+15t^2+25t+25)}$  \\
    & \quad{$ = 12^3 + \, {\rm etc.}$} \\ 
\noalign{\smallskip}\hline
\end{tabular}
}%
\end{center}
\label{tab:coverings}
\end{table}

\begin{table}
\caption{The degree-$\psi(N)$ map $t_N\mapsto j'$.}
\begin{center}
{\small
\begin{tabular}{cl}
\hline\noalign{\smallskip}
$N$ & \hfil $j'(\tau):=j(N\tau)$ as a function of~$t_N(\tau)$ \\
\noalign{\smallskip}\hline\noalign{\smallskip}
$2$ & $\frac{(t+256)^3}{t^2}$ \\
    & $\quad = 12^3 + \frac{(t+64)(t-512)^2}{t^2}$ \\
$3$ & $\frac{(t+27)(t+243)^3}{t^3}$ \\
    & $\quad = 12^3 + \frac{(t^2-486t-19683)^2}{t^3}$ \\
$4$ & $\frac{(t^2+256t+4096)^3}{t^4(t+16)}$ \\
    & $\quad = 12^3 + \frac{(t+32)^2(t^2-512t-8192)^2}{t^4(t+16)}$ \\
$5$ & $\frac{(t^2+250t+3125)^3}{t^5}$ \\
    & $\quad = 12^3 + \frac{(t^2+22t+125)(t^2-500t-15625)^2}{t^5}$ \\
$6$ & $\frac{(t+12)^3(t^3+252t^2+3888t+15552)^3}{t^6(t+8)^2(t+9)^3}$ \\
    & $\quad = 12^3 + \frac{(t^2+36t+216)^2(t^4-504t^3-13824t^2-124416t-373248)^2}{t^6(t+8)^2(t+9)^3}$ \\
$7$ & $\frac{(t^2+13t+49)(t^2+245t+2401)^3}{t^7}$ \\
    & $\quad = 12^3 + \frac{(t^4-490t^3-21609t^2-235298t-823543)^2}{t^7}$ \\
$8$ & $\frac{(t^4+256t^3+5120t^2+32768t+65536)^3}{t^8(t+4)(t+8)^2}$ \\
    & $\quad = 12^3 + \frac{(t^2+32t+128)^2(t^4-512t^3-10240t^2-65536t-131072)^2}{t^8(t+4)(t+8)^2}$ \\
$9$ & $\frac{(t+9)^3(t^3+243t^2+2187t+6561)^3}{t^9(t^2+9t+27)}$ \\
    & $\quad = 12^3 + \frac{(t^6-486t^5-24057t^4-367416t^3-2657205t^2-9565938t-14348907)^2}{t^9(t^2+9t+27)}$ \\
$10$ & $\frac{(t^6+260t^5+6400t^4+64000t^3+320000t^2+800000t+800000)^3}{t^{10}(t+4)^2(t+5)^5}$ \\
    & \quad{\scriptsize$= 12^3 + \frac{(t^2+8t+20)(t^2+12t+40)^2(t^2+30t+100)^2(t^4-520t^3-6600t^2-28000t-40000)^2}{t^{10}(t+4)^2(t+5)^5}$} \\
$12$ & $\frac{(t^2+12t+24)^3(t^6+252t^5+4392t^4+31104t^3+108864t^2+186624t+124416)^3}{t^{12}(t+2)(t+3)^3(t+4)^4(t+6)^3}$ \\
    & \quad{$ = 12^3 +\, {\rm etc.}$} \\ 
$13$ & $\frac{(t^2+5t+13)(t^4+247t^3+3380t^2+15379t+28561)^3}{t^{13}}$ \\
    & \quad{\scriptsize $= 12^3 + \frac{(t^2+6t+13)(t^6-494t^5-20618t^4-237276t^3-1313806t^2-3712930t-4826809)^2}{t^{13}}$} \\
$16$ & {\scriptsize$\frac{(t^8+256t^7+5632t^6+53248t^5+282624t^4+917504t^3+1835008t^2+2097152t+1048576)^3}{t^{16}(t+2)(t+4)^4(t^2+4t+8)}$}\\
    & \quad{$ = 12^3 +\, {\rm etc.}$} \\ 
$18$ & {\scriptsize$\frac{{{(t^3+12t^2+36t+36)^3}{(t^9+252t^8+4644t^7+39636t^6+198288t^5+629856t^4+1294704t^3+1679616t^2+1259712t+419904)^3}}}{t^{18}(t+2)^2(t+3)^9(t^2+3t+3)(t^2+6t+12)^2}$} \\
    & $\quad = 12^3 + \,{\rm etc.}$ \\
$25$ & {\scriptsize$\frac{(t^{10}+250t^9+4375t^8+35000t^7+178125t^6+631250t^5+1640625t^4+3125000t^3+4296875t^2+3906250t+1953125)^3}{t^{25}(t^4+5t^3+15t^2+25t+25)}$} \\
    & $\quad = 12^3 + \,{\rm etc.}$ \\
\noalign{\smallskip}\hline
\end{tabular}
}%
\end{center}
\label{tab:coverings2}
\end{table}

\subsection{Coverings and parametrized modular equations}
\label{subsec:parametrized1b}

An explicit formula for each cover $X_0(N)/X(1),$ i.e., for the covering
map $t_N\mapsto j$ of degree $\psi(N),$ is given in
Table~\ref{tab:coverings}.  The rational expressions originated with Klein
and Gierster, but have been modified to agree with our Hauptmodul
normalization convention, which ensures monicity of each polynomial factor.
The factored expressions add much detail to the ramification data of
Table~\ref{tab:basic}.  For any~$N,$ the roots of the denominator are the
$t_N$-values of cusps; and rational and irrational cusps are easy to
distinguish.  The multiplicity of each root is the width of the
corresponding cusp.  Each numerator contains one or more cubic factors, and
possibly a polynomial that is not cubed; if the latter is present, its
roots are the cubic fixed points of~$\Gamma_0(N)$ on~$X_0(N)$.  Similarly,
an unsquared factor (if~any) of the numerator of ${j-12^3}$ has quadratic
fixed points as its roots.  The cusp and elliptic point locations read~off
from Table~\ref{tab:coverings} are listed in the final two columns of
Table~\ref{tab:hauptmoduln}.

The covering maps of Table~\ref{tab:coverings} can be trusted, since they
agree with the known ramification data for $t_N\mapsto j$.  For each~$N,$
the ramification data constitute a schema, such~as the $N=5$ schema, which
is 
\begin{equation}
\psi(5)=6=5+1=3+3={2+2+1+1}
\end{equation}
(in~an obvious notation).  The final three members specify multiplicities
on the fibres over $j=\infty,0,12^3,$ respectively.  Ordinarily a covering
of~$\mathbb{P}^1(\mathbb{C})$ by~$\mathbb{P}^1(\mathbb{C}),$ unramified
over the complement of three points, is not uniquely specified by its
ramification schema.  Rather, by the Grothendieck correspondence it is
specified by its associated \emph{dessin d'enfants}~\cite{Schneps94}.  But
it is the case that for each prime~$N$ in the table, there is only a single
dessin compatible with the schema.  This is a combinatorial statement,
in~fact a graph-theoretic one, which can be verified for each~$N$
individually.  There is also a deeper, Galois-theoretic reason why it is
true.  Within the symmetric group $\mathfrak{S}_{\psi(N)}$ of permutations
of the sheets of the covering, the ramification schema specifies a disjoint
cycle decomposition of the monodromy generators $g_\infty,g_0,g_{12^3}$
associated to loops around $j=\infty,0,12^3$; that~is, it specifies each of
them up~to conjugacy.  For all~$N,$ it can be shown that this triple of
conjugacy classes is {\em rigid\/}: there is essentially only one way of
embedding it in~$\mathfrak{S}_{\psi(N)}$ that is compatible with the
constraint $g_\infty g_0g_{12^3}=1$.  This is the same as saying that only
one dessin is possible.  For a discussion, see
Elkies~\cite[p.~49]{Elkies98} and Malle and
Matzat~\cite[Chap.~I,\,\S7.4]{Malle99}.

Each `canonical modular equation' $j=j(t_N)$ of Table~\ref{tab:coverings}
immediately yields a parametrization of the corresponding classical modular
equation, i.e., of the equation ${\Phi_N(j,j')=0},$ where
$j'(\tau)=j(N\tau)$.  Since $j|_{W_N}=j'$ and $t_N|_{W_N}=\kappa_N/t_N,$
applying $W_N$ to the formula $j=j(t_N)$ yields the other formula
$j'=j'(t_N)$.  The rational expressions for $j'$~as a function of~$t_N,$
and also for~$j'-12^3,$ are given in fully factored form in
Table~\ref{tab:coverings2}.  A~few expressions for $j'-12^3,$ as for
${j-12^3},$ are omitted on account of length.

It should be noted that combined with the formula for the
$\lambda$-invariant, namely $\lambda(\tau)=[t_4/(t_4+\nobreak16)](\tau/2),$
the formulas for $j,j'$ in~terms of~$t_4$ yield Eq.~(\ref{eq:mformula})
of~\S\ref{subsec:12}.  It should also be noted that an uncubed factor
appears in the numerator of any~$j(t_N)$ if~and only~if the {\em same\/}
factor appears in that of~$j'(t_N)$; and that similarly, an unsquared
factor appears in the numerator of $j(t_N)-12^3$ if~and only~if it appears
in that of~$j'(t_N)-12^3$.  Such factors (i.e., their roots) correspond to
non-isomorphic self-$N$-isogenies of equianharmonic and lemniscatic
elliptic curves, respectively.  For instance, the roots of
$t_5^2+22t_5+125,$ the common unsquared factor of the numerators of
$j(t_5)-12^3,\allowbreak j'(t_5)-12^3,$ are bijective with the
representations of $N=5$ as a sum of two squares, i.e., $5=1^2+2^2$
and~$5=2^2+1^2$.

Taken together, Tables \ref{tab:coverings} and~\ref{tab:coverings2}
comprise the 14~rationally parametrized modular equations of level~$1$.
When $N$~is composite, functional decompositions of the coverings
$j=j(t_N),$ $j'=j'(t_N)$ are valuable; these will be given in the next
section.  But for prime~$N,$ the tables are useful as they stand.  In~the
spirit of Mathews~\cite{Mathews1890} and Fricke, one can use them to
compute singular moduli and class invariants.  The case $N=5$ is typical.
Factoring $j(t_5)-j'(t_5)$ yields some factors that one can show are
extraneous, together with the polynomial $t_5^2-125,$ which has roots
$t_5=\pm5\sqrt5$.  There are two because the discriminant $D=-5$ has class
number~$2$.  The associated values of~$j,$ computed from the formula
${j=j(t_5),}$ have minimal polynomial over~$\mathbb{Q}$ equal to
$j^2-2^75^3 79j-880^3,$ which is the class polynomial
of~$\mathbb{Q}(\sqrt{-5})$.

\section{Parametrized Modular Equations at Higher Levels}
\label{sec:parametrized2}

To generate rationally parametrized modular equations at higher levels, one
reasons as~follows.  The invariants $j,j'$ of the last section, where
$j'(\tau):=j(N\tau),$ are Hauptmoduln for $\Gamma(1)$ and the conjugated
subgroup $\Gamma(1)':=w_N^{-1}\Gamma\,w_N$ of ${\it PSL}(2,\mathbb{R}),$
where $w_N=\left(\begin{smallmatrix}0&-1\\N&0\end{smallmatrix}\right)$ is
the Fricke involution.  Since $\Gamma(1)\cap\Gamma(1)'=\Gamma_0(N),$ if
$\Gamma_0(N)$~is of genus zero with Hauptmodul~$t_N$ then one must have
$j,j'\in\mathbb{C}(t_N)$.  In~fact, quotienting~$\mathcal{H}^*$ yields a
pair of covers $X_0(N)/X(1),\allowbreak X_0(N)/X(1)',$ the projections of
which were given in Tables \ref{tab:coverings} and~\ref{tab:coverings2}.
Similarly, for any~$1<d\divides{N}$ one has $\Gamma_0(N)<\Gamma_0(d),$ and
the pair of groups
$\Gamma_0(d),\allowbreak\Gamma_0(d)':=w_N^{-1}\Gamma_0(d)w_N,$ with
Hauptmoduln $t_d,t_d'$ where $t_d'(\tau):=t_d((N/d)\tau),$ have
intersection~$\Gamma_0(N)$.  One must have $t_d,t_d'\in\mathbb{C}(t_N),$
and a pair of covers $X_0(N)/X_0(d),X_0(N)/X_0(d)'$.

For the nine composite~$N$ for which $\Gamma_0(N)$ is of genus zero, the
corresponding rational maps $t_d=t_d(t_N)$ were derived by
Gierster~\cite{Gierster1879}; not all are in Fricke.  They are listed in
the third column of Table~\ref{tab:intermediate}, modified to agree with
our normalization convention.  Like the covers $X_0(N)/X(1),\allowbreak
X_0(N)/X(1)'$ for prime~$N,$ one can show that the covers $X_0(N)/X_0(d)$
are uniquely determined by ramification data.  This uniqueness enabled
Gierster and Fricke to work them~out from `pictorial' ramification data,
i.e., from figures showing how a fundamental region of~$\Gamma_0(N)$
comprises $\psi(N)/\psi(d)$ ones of~$\Gamma_0(d)$.
Knopp~\cite[\S7.6]{Knopp70} gives a linear-algebraic, nonpictorial
derivation of the projection map $t_5=t_5(t_{25})$ of the cover
$X_0(25)/X_0(5),$ but most of the others are not well known.


\begin{proposition}
  For all\/ $1<d\divides{N}$ listed in Table\/~{\rm\ref{tab:intermediate}},
  the Hauptmodul\/~$t_d$ is a polynomial\/ {\rm(}rather than merely
  rational\/{\rm)} function of the Hauptmodul\/~$t_N$ iff all primes that
  divide\/~$N$ also divide\/~$d$.  This occurs if\/ $\Gamma_0(N)$ is a
  normal subgroup of\/~$\Gamma_0(d),$ though the converse does not hold.
\label{prop:short}
\end{proposition}
\begin{proof}
For any $1<d\divides{N},$ the cusps
$\bigl[\frac11\bigr]_d,\bigl[\frac11\bigr]_N$ have widths $d,N,$ so if
$\Gamma_0(d),\allowbreak\Gamma_0(N)$ are of genus zero with Hauptmoduln
$t_d,t_N$ having poles at $\tau=0,$ the rational function $t_d=t_d(t_N)$
must take $t_N=\infty$ to $t_d=\infty$ with multiplicity~$N/d$.  This
function has degree $\psi(N)/\psi(d),$ which equals $N/d$ iff all primes
that divide~$N/d$ also divide~$d$.  The final sentence follows from a
result of Cummins~\cite[Prop.~7.1]{Cummins2004}: $\Gamma_0(N)$ is normal
in~$\Gamma_0(d)$ iff $(N/d) \divides (d,24)$.  The cases $d=2,N=8,16$ show
that the converse does not hold.
\end{proof}

We computed each of the coverings $X_0(N)/X_0(d)'$ shown in
Table~\ref{tab:intermediate}, i.e., each formula $t_d'=t_d'(t_N)$ in the
fourth column, by applying $W_N$ to the corresponding formula
$t_d=t_d(t_N),$ noting that just as $t_N|_{W_N}=\kappa_N/t_N,$ so
$t_d|_{W_N}(\tau)=\kappa_d/t_d((N/d)\tau)$.  The several functional
decompositions appearing in columns $3$ and~$4$ have an intuitive
explanation.  They come from composite covers: e.g., the formula
$[t(t+16)\circ t(t+8)\circ t(t+4)](t_{16})$ for~$t_2$ is associated with
the composite cover $X_0(16)/X_0(8)/X_0(4)/X_0(2)$.  It is the composition
$t_{16}\mapsto t_8\mapsto t_4\mapsto t_2$.  Distinct paths from
$\Gamma_0(d)$ to~$\Gamma_0(N)$ in the subgroup lattice yield distinct
composite covers.  For instance, the two distinct representations given for
the rational function $t_2=t_2(t_{12})$ come from $X_0(12)/X_0(6)/X_0(2)$
and $X_0(12)/X_0(4)/X_0(2),$ respectively.  In the same way, the two
representations given for $t_2(6\tau),$ the Hauptmodul for
$\Gamma_0(2)'=w_{12}^{-1}\Gamma_0(2)w_{12},$ in~terms of~$t_{12},$ the
Hauptmodul for~$\Gamma_0(12),$ come from $t_{12}(\tau)\mapsto
t_4(3\tau)\mapsto t_2(6\tau)$ and $t_{12}(\tau)\mapsto t_6(2\tau)\mapsto
t_2(6\tau)$.

\begin{table}
\caption{The covering maps from $X_0(N)\ni t_N$ to $X_0(d)\ni t_d$ and $X_0(d)'\ni t_d',$ for $1<d\divides{N}$.}
\begin{center}
{\small
\begin{tabular}{ccll}
\noalign{\smallskip}\hline\noalign{\smallskip}
$N$ & $d$ & $t_d(\tau),$  & $t_d'(\tau):=t_d(\frac{N}d\tau),$ \\
    &     & as a function of~$t:= t_N(\tau)$ & as a function of~$t:= t_N(\tau)$\\
\noalign{\smallskip}\hline\noalign{\vskip1pt}\hline\noalign{\vskip1.5pt}
$4$ & $2$ & $t(t+16)$ & $\frac{t^2}{t+16}$ \\
\noalign{\vskip 1pt}\hline\noalign{\vskip 1pt}
$6$ & $2$ & $\frac{t(t+8)^3}{t+9}$ & $\frac{t^3(t+8)}{(t+9)^3}$ \\
$6$ & $3$ & $\frac{t(t+9)^2}{t+8}$ & $\frac{t^2(t+9)}{(t+8)^2}$ \\
\noalign{\vskip 1pt}\hline\noalign{\vskip 1pt}
$8$ & $2$ & $t(t+16)\,\circ\,t(t+8)$ & $\frac{t^2}{t+16} \circ \frac{t^2}{t+4}$ \\
$8$ & $4$ & $t(t+8)$ & $\frac{t^2}{t+4}$ \\
\noalign{\vskip 1pt}\hline\noalign{\vskip 1pt}
$9$ & $3$ & $t(t^2+9t+27)$ & $\frac{t^3}{t^2+9t+27}$ \\
\noalign{\vskip 1pt}\hline\noalign{\vskip 1pt}
$10$ & $2$ & $\frac{t(t+4)^5}{t+5}$ & $\frac{t^5(t+4)}{(t+5)^5}$ \\
$10$ & $5$ & $\frac{t(t+5)^2}{t+4}$ & $\frac{t^2(t+5)}{(t+4)^2}$ \\
\noalign{\vskip 1pt}\hline\noalign{\vskip 1pt}
$12$ & $2$ & $\frac{t(t+8)^3}{t+9}\circ\,t(t+6)$ & $\frac{t^3(t+8)}{(t+9)^3}\circ\frac{t^2}{t+2}$ \\
     &     & $=t(t+16)\, \circ\frac{t(t+4)^3}{t+3}$ & $=\frac{t^2}{t+16}\circ\frac{t^3(t+4)}{(t+3)^3}$ \\
$12$ & $3$ & $\frac{t(t+9)^2}{t+8}\circ\, t(t+6)$ & $\frac{t^2(t+9)}{(t+8)^2}\circ \frac{t^2}{t+2}$ \\
$12$ & $4$ & $\frac{t(t+4)^3}{t+3}$ & $\frac{t^3(t+4)}{(t+3)^3}$ \\
$12$ & $6$ & $t(t+6)$ & $\frac{t^2}{t+2}$ \\
\noalign{\vskip 1pt}\hline\noalign{\vskip 1pt}
$16$ & $2$ & $t(t+16)\,\circ\,t(t+8)\,\circ\,t(t+4)$ & $\frac{t^2}{t+16}\circ\frac{t^2}{t+4}\circ\frac{t^2}{t+2}$ \\
$16$ & $4$ & $t(t+8)\,\circ\,t(t+4)$ & $\frac{t^2}{t+4}\circ\frac{t^2}{t+2}$ \\
$16$ & $8$ & $t(t+4)$ & $\frac{t^2}{t+2}$ \\
\noalign{\vskip 1pt}\hline\noalign{\vskip 1pt}
$18$ & $2$ & $\frac{t(t+8)^3}{t+9}\circ\,t(t^2+6t+12)$ & $\frac{t^3(t+8)}{(t+9)^3}\circ\frac{t^3}{t^2+3t+3}$ \\
$18$ & $3$ & $t(t^2+9t+27)\,\circ\frac{t(t+3)^2}{t+2}$ & $\frac{t^3}{t^2+9t+27}\circ\frac{t^2(t+3)}{(t+2)^2}$ \\
     &     & $=\frac{t(t+9)^2}{t+8}\circ\,t(t^2+6t+12)$  & $=\frac{t^2(t+9)}{(t+8)^2}\circ\frac{t^3}{t^2+3t+3}$ \\
$18$ & $6$ & $t(t^2+6t+12)$ & $\frac{t^3}{t^2+3t+3}$ \\
$18$ & $9$ & $\frac{t(t+3)^2}{t+2}$ & $\frac{t^2(t+3)}{(t+2)^2}$ \\
\noalign{\vskip 1pt}\hline\noalign{\vskip 1pt}
$25$ & $5$ & $t(t^4+5t^3+15t^2+25t+25)$ & $\frac{t^5}{t^4+5t^3+15t^2+25t+25}$ \\
\noalign{\vskip2pt}\hline
\end{tabular}
}%
\end{center}
\label{tab:intermediate}
\end{table}

Taken together, columns $3,4$ of Table~\ref{tab:intermediate} list all
rationally parametrized modular equations of level greater than unity (the
level being~$d,$ and the degree~$N/d$).  As an application, one can compute
parametrized modular equations of the $\alpha$\textendash\nobreak$\beta$
type mentioned in~\S\ref{subsec:12}.  Using the formula
$\alpha=t_4/(t_4+16),$ and defining $\alpha:=\alpha(\tau),\allowbreak
\beta:=\alpha\left((N/d)\tau\right),$ converts the three equations
with~$d=4$ (and degrees $N/d=2,3,4$) respectively to
\begin{subequations}
\begin{align}
\alpha&=\frac{t(t+8)}{(t+4)^2},& \beta&=\frac{t^2}{(t+8)^2}\label{eq:mn2};\\
\alpha&=\frac{t(t+4)^3}{(t+2)^3(t+6)},& \beta&=\frac{t^3(t+4)}{(t+2)(t+6)^3}\label{eq:mn3};\\
\alpha&=\frac{t(t+8)}{(t+4)^2}\,\circ\,t(t+4),& \beta&=\frac{t^2}{(t+8)^2}\,\circ
\,\frac{t^2}{t+2}.\label{eq:mn4}
\end{align}
\end{subequations}
The parameter~$t$ signifies $t_8,t_{12},t_{16},$ respectively.
Equation~(\ref{eq:mn2}), of degree~$2,$ is a parametrization of the
$\alpha$\textendash\nobreak$\beta$ relation~(\ref{eq:two}), and hence of
Landen's transformation; cf.~(\ref{eq:duo}).  Equation~(\ref{eq:mn3}), of
degree~$3,$ is also classical; it was discovered by Legendre and
rediscovered by Jacobi.  The derivation of Cayley \cite[\S265]{Cayley1895}
is perhaps the most accessible.\footnote{Cayley's equation relates the
value of his invariant~$k^2$ at any point $\tau\in\mathcal{H}$ to its value
at $3\tau\in\mathcal{H}$.  The formula he gives for the Klein invariant
$J=j/12^3$ in~terms of~$k^2$ (see his~\S300) makes clear that $k^2$ is to
be interpreted as our $\alpha$-invariant, i.e.,
$\alpha(\tau)=\lambda(2\tau)$. Cf.~Eq.~(\ref{eq:mformula}).}  (His
uniformizing parameter is not our~$t_{12},$ but rather the alternative
Hauptmodul $t_{12}/(t_{12}+6)$ for~$\Gamma_0(12)$; cf.\
Table~\ref{tab:alt}.)  Equation~(\ref{eq:mn4}), of degree~$4,$ is classical
too, though it may not have appeared in this form before; it is the basis
of the little-known quartic (i.e., {\em biquadratic\/})
arithmetic--geometric mean iteration~\cite[p.~17]{Borwein87}.  Cayley's
method of deriving a modular equation of prime degree~$p$ for the
$\alpha$-invariant is difficult to apply when~$p>3,$ and one now sees why:
if~$p>3$ then $\Gamma_0(4p)$ is of positive genus, and no rational
parametrization exists.

\begin{table}
\caption{Rationally parametrized modular equations for the $j$-invariant,
of all degrees~$N$ for which the curve $X_0(N)$~is of genus zero.  Here
$j,j'$ and~$t$ signify $j(\tau),$ $j(N\tau)$ and~$t_N(\tau)$.}
\begin{center}
{\small
\begin{tabular}{cll}
\hline\noalign{\smallskip}
$N$ & $\quad j$ & $\quad j'$\\
\noalign{\smallskip}\hline\noalign{\vskip1pt}\hline\noalign{\vskip1.5pt}
$2$ & $\frac{(t+16)^3}{t}$ & $\frac{(t+256)^3}{t^2}$\\
\noalign{\vskip 1pt}\hline\noalign{\vskip 1pt}
$3$ & $\frac{(t+27)(t+3)^3}t$ & $\frac{(t+27)(t+243)^3}{t^3}$\\
\noalign{\vskip 1pt}\hline\noalign{\vskip 1pt}
$4$ & $\frac{(t+16)^3}t \circ\, t(t+16)$ & $\frac{(t+256)^3}{t^2} \circ \frac{t^2}{t+16}$\\
\noalign{\vskip 1pt}\hline\noalign{\vskip 1pt}
$5$ & $\frac{(t^2+10t+5)^3}{t}$ & $\frac{(t^2+250t+3125)^3}{t^5}$ \\
\noalign{\vskip 1pt}\hline\noalign{\vskip 1pt}
$6$ & $\frac{(t+16)^3}{t} \circ \frac{t(t+8)^3}{t+9}$ & $\frac{(t+256)^3}{t^2} \circ \frac{t^3(t+8)}{(t+9)^3}$ \\
    & $\quad=\frac{(t+27)(t+3)^3}t \circ \frac{t(t+9)^2}{t+8}$ & $\quad=\frac{(t+27)(t+243)^3}{t^3} \circ \frac{t^2(t+9)}{(t+8)^2}$ \\
\noalign{\vskip 1pt}\hline\noalign{\vskip 1pt}
$7$ & $\frac{(t^2+13t+49)(t^2+5t+1)^3}{t}$ & $\frac{(t^2+13t+49)(t^2+245t+2401)^3}{t^7}$ \\
\noalign{\vskip 1pt}\hline\noalign{\vskip 1pt}
$8$ & $\frac{(t+16)^3}{t}\circ\,t(t+16)\,\circ\,t(t+8)$ & $\frac{(t+256)^3}{t^2}\circ\frac{t^2}{t+16} \circ \frac{t^2}{t+4}$ \\
\noalign{\vskip 1pt}\hline\noalign{\vskip 1pt}
$9$ & $\frac{(t+27)(t+3)^3}t\circ\, t(t^2+9t+27)$ & $\frac{(t+27)(t+243)^3}{t^3}\circ\frac{t^3}{t^2+9t+27}$ \\
\noalign{\vskip 1pt}\hline\noalign{\vskip 1pt}
$10$ & $\frac{(t+16)^3}{t}\circ\frac{t(t+4)^5}{t+5}$ & $\frac{(t+256)^3}{t^2}\circ\frac{t^5(t+4)}{(t+5)^5}$ \\
     & $\quad=\frac{(t^2+10t+5)^3}{t}\circ\frac{t(t+5)^2}{t+4}$ & $\quad=\frac{(t^2+250t+3125)^3}{t^5}\circ\frac{t^2(t+5)}{(t+4)^2}$ \\
\noalign{\vskip 1pt}\hline\noalign{\vskip 1pt}
$12$ & $\frac{(t+16)^3}{t}\circ\frac{t(t+8)^3}{t+9}\circ\,t(t+6)$ & $\frac{(t+256)^3}{t^2}\circ\frac{t^3(t+8)}{(t+9)^3}\circ\frac{t^2}{t+2}$ \\
     & $\quad=\frac{(t+27)(t+3)^3}t\circ\frac{t(t+9)^2}{t+8}$ & $\quad=\frac{(t+27)(t+243)^3}{t^3}\circ\frac{t^2(t+9)}{(t+8)^2}\circ \frac{t^2}{t+2}$ \\
     & $\qquad\qquad\qquad\qquad\qquad{}\circ\, t(t+6)$ &  \\
     & $\quad=\frac{(t+16)^3}{t}\circ\,t(t+16)\, \circ\frac{t(t+4)^3}{t+3}$ & $\quad=\frac{(t+256)^3}{t^2}\circ\frac{t^2}{t+16}\circ\frac{t^3(t+4)}{(t+3)^3}$ \\
\noalign{\vskip 1pt}\hline\noalign{\vskip 1pt}
$13$ & {\scriptsize$\frac{(t^2+5t+13)(t^4+7t^3+20t^2+19t+1)^3}{t}$} & {\scriptsize$\frac{(t^{2}+5t+13)(t^{4}+247t^{3}+3380t^{2}+15379t+28561)^{3}}{t^{13}}$} \\
\noalign{\vskip 1pt}\hline\noalign{\vskip 1pt}
$16$ & $\frac{(t+16)^3}{t}\circ\, t(t+16)\,\circ\,t(t+8)$ & $\frac{(t+256)^3}{t^2}\circ\frac{t^2}{t+16}\circ\frac{t^2}{t+4}\circ\frac{t^2}{t+2}$ \\
     & $\qquad\qquad\qquad\qquad\qquad{}\circ\,t(t+4)$ & \\
\noalign{\vskip 1pt}\hline\noalign{\vskip 1pt}
$18$ & {\scriptsize$\frac{(t+16)^3}{t}\circ\frac{t(t+8)^3}{t+9}\circ\,t(t^2+6t+12)$} & {\scriptsize$\frac{(t+256)^3}{t^2}\circ\frac{t^3(t+8)}{(t+9)^3}\circ\frac{t^3}{t^2+3t+3}$} \\
     & {\scriptsize$\quad=\frac{(t+27)(t+3)^3}t\circ\frac{t(t+9)^2}{t+8}$}  & {\scriptsize$\quad=\frac{(t+27)(t+243)^3}{t^3}\circ\frac{t^2(t+9)}{(t+8)^2}\circ\frac{t^3}{t^2+3t+3}$} \\
     & {\scriptsize$\qquad\qquad\qquad\qquad\qquad{}\circ\,t(t^2+6t+12)$} & \\
     & {\scriptsize$\quad=\frac{(t+27)(t+3)^3}t\circ\,t(t^2+9t+27)$} & {\scriptsize$\quad=\frac{(t+27)(t+243)^3}{t^3}\circ\frac{t^3}{t^2+9t+27}\circ\frac{t^2(t+3)}{(t+2)^2}$} \\
     & {\scriptsize$\qquad\qquad\qquad\qquad\qquad{}\circ\frac{t(t+3)^2}{t+2}$} & \\
\noalign{\vskip 1pt}\hline\noalign{\vskip 1pt}
$25$ & $\frac{(t^2+10t+5)^3}{t}$ & $\frac{(t^2+250t+3125)^3}{t^5}$ \\
     & $\quad{}\circ\, t(t^4+5t^3+15t^2+25t+25)$ & {\scriptsize$\qquad\qquad\qquad{}\circ\frac{t^5}{t^4+5t^3+15t^2+25t+25}$} \\
\noalign{\vskip2pt}\hline
\end{tabular}
}%
\end{center}
\label{tab:jfactored}
\end{table}

The parametrizations of the classical modular equations $\Phi_N(j,j')=0$
given in Tables \ref{tab:coverings} and~\ref{tab:coverings2} were quite
complicated, and by exploiting Table~\ref{tab:intermediate} one can write
them in more understandable form.  For each~$d\divides{N},$ $j$~is
rationally expressible in~terms of~$t_d,$ which in~turn is rationally
expressible in~terms of~$t_N$.  This expresses $j=j(t_N)$ as a composition,
and $j'$ (i.e.,~$j(N\tau)$) can be similarly expressed.  Taking into
account the many pairs~$d,N$ of Table~\ref{tab:intermediate}, one obtains
Table~\ref{tab:jfactored}, an improved version of Tables
\ref{tab:coverings} and~\ref{tab:coverings2} that displays each functional
composition.  For example, the three representations given for the
projection $j=j(t_{12})$ of $X_0(12)/X(1)$ come from the composite covers
$X_0(12)/X_0(6)/X_0(2)/X(1),$ $X_0(12)/X_0(6)/X_0(3)/X(1),$ and
$X_0(12)/X_0(4)/X_0(2)/X(1),$ respectively.

No table resembling Table~\ref{tab:jfactored} has appeared in~print before,
and it may prove useful, e.g., in the numerical computation of transformed
invariants~$j'$.  In many cases it clarifies the relation between~$j,j'$.
It~is clear from a glance that the modular equation $\Phi_N(j,j')=0$
has solvable Galois group for many of the listed values of~$N,$ but that
the group is not solvable if $N=5,10,25$.

The table lists only $j(\tau)$ and~$j(N\tau)$ as a function of~$t_N(\tau),$
but in~fact $j(d\tau)$~is rational in~$t_N(\tau)$ for all~$d\divides{N},$
the rational map being the composition $t_N(\tau)\mapsto t_d(\tau)\mapsto
j(d\tau)$.  The formulas for the cases $1<d<N$ are easily worked~out, and
include two that throw light on Weber's functions.  Recall that his
functions
\begin{equation}
  \gamma_2(\tau):=j^{1/3}(\tau),\qquad \gamma_3(\tau):=(j-12^3)^{1/2}(\tau)
\end{equation}
are useful in the computation of singular moduli.  By group theory one can
prove that $\gamma_2(3\tau),\gamma_3(2\tau)$ are automorphic under
$\Gamma_0(9),\Gamma_0(4),$ respectively \cite[\S12A]{Cox89}.  But by direct
computation one obtains from Tables \ref{tab:intermediate}
and~\ref{tab:jfactored},
{\small
\begin{equation}
\gamma_2(3\tau)=\frac{(t_9+3)(t_9+9)(t_9^2+27)}{t_9\,(t_9^2+9t_9+27)},\qquad \gamma_3(2\tau)=\frac{(t_4-16)(t_4+8)(t_4+32)}{t_4(t_4+16)}.
\label{eq:verytired}
\end{equation}
}%
These are more precise statements about modularity, since they indicate the
locations of zeroes and poles.

The first equation in~(\ref{eq:verytired}) also clarifies the modular
setting of the Hesse--Dixon family of elliptic curves.  Any elliptic curve
$E/\mathbb{C}$ has a cubic Hesse model, the function field of which was
studied by Dixon~\cite{Dixon1890}; namely,
\begin{equation}
  x^3+y^3 + 1 - (\gamma+3)xy = 0.
\label{eq:hessedixonmodel}
\end{equation}
Here
$\gamma\in\mathbb{C}\setminus\{0,3(\zeta_3-1),3(\zeta^2_3-1)\}=\mathbb{C}\setminus\{0,\frac{-9\pm3\sqrt{-3}}{2}\}$
is the Hesse--Dixon parameter.  This $\gamma$-invariant may be chosen to be
a single-valued function of $\tau\in\mathcal{H},$ and the $j$-invariant
expressed in~terms of it by
{\small
\begin{equation}
j(\tau)=\frac{(\gamma+3)^3(\gamma+9)^3(\gamma^2+27)^3}{\gamma^3\,(\gamma^2+9\gamma+27)^3}(\tau)
=
\frac{(\gamma+3)^3(\gamma^3+9\gamma^2+27\gamma+3)^3}{\gamma(\gamma^2+9\gamma+27)}(3\tau).
\label{eq:gammaformula}
\end{equation}
}%
In fact, one may take $\gamma(\tau)=t_9(\tau/3)$.  (To~see all this,
compute the $j$-invariant of the model~(\ref{eq:hessedixonmodel}),
obtaining the first equality in~(\ref{eq:gammaformula}), and notice the
equivalence to~(\ref{eq:verytired}a).)  Since $t_9(\tau/3)$ is a Hauptmodul
for a subgroup conjugated to~$\Gamma_0(9)$ in ${\it PSL}(2,\mathbb{R})$ by
a $3$-isogeny, namely~$\Gamma(3),$ the Hesse--Dixon model is associated
to~$\Gamma(3),$ just as the Legendre and Jacobi models are to $\Gamma(2)$
and~$\Gamma_0(4)$.  This was shown by Beauville~\cite{Beauville82}, but the
present derivation is more concrete than his.  Since
$t_9=3^3\cdot[9]^3/\,[1]^3,$ one sees that the $\gamma$-invariant equals
$3^3\cdot[3]^3/\,[\frac13]^3$.  Its $q$-expansion is
$3^3\cdot[q_3+3q_3^2+9q_3^3+22q_3^4+\cdots],$ with
$q_3:=q^{1/3}=e^{2\pi{\rm i}\tau/3}$.

The parameters $\gamma,\gamma'$ of Hesse--Dixon elliptic curves $E,E'$ with
period ratios $\tau,\tau'$ satisfying $\tau'=2\tau$ are linked
parametrically by
\begin{equation}
  \gamma=\frac{t(t+3)^2}{t+2},\qquad \gamma'=\frac{t^2(t+3)}{(t+2)^2}.
\end{equation}
This follows from the degree-$2$ $t_9$\textendash\nobreak$t_9'$ modular
equation in Table~\ref{tab:intermediate}; `$t$'~is $t_{18}(\tau/3)$.  This
degree-$2$ $\gamma$\textendash\nobreak$\gamma'$ relation agrees with the
unparametrized one derived by Dixon \cite[\S79]{Dixon1890}.  He also
worked~out the multivalued function $\gamma'=\gamma'(\gamma),$ expressible
in~terms of radicals, for the case $\tau'=3\tau$~\cite[\S85]{Dixon1890}.
But the cubic $\gamma$\textendash\nobreak$\gamma'$ modular relation cannot
be rationally parametrized, since $X_0(27)$ unlike~$X_0(18)$ is not of
genus zero.

\section{The Canonical Weight-$1$ Modular Form $\mathfrak{h}_N(\tau)$ for~$\Gamma_0(N)$}
\label{sec:modular}

The preceding sections focused on the canonical Hauptmoduln~$t_N$ for the
genus-zero congruence subgroups $\Gamma_0(N)$ of~$\Gamma(1),$ and on
$t_N$\textendash\nobreak$t_N'$ modular equations, which are Hauptmodul
relations.  We now begin the derivation of modular equations (including
`multipliers') for the corresponding elliptic families.  These are really
identities satisfied by Gauss--Manin connections, or in classical language,
transformation laws for Picard--Fuchs equations.  But because of our
interest in ${}_2F_1$ and other special function identities, we shall
develop them in a very concrete way, by defining certain weight-$1$ modular
forms $\mathfrak{h}_N(\tau)$ that are (multivalued) functions of the
corresponding Hauptmoduln, according to
$\mathfrak{h}_N(\tau)=h_N(t_N(\tau))$.  Each~$h_N$ is a solution of a
normal-form Picard--Fuchs equation.  Each modular form~$\mathfrak{h}_N,$ or
equivalently the multivalued function~$h_N,$ is defined by a pullback along
$X_0(N)\to X(1)$.

\begin{definition} 
If\/ $\Gamma<\Gamma(1)$ is of genus zero with Hauptmodul~$t,$ and
  $j$~equals $P(t)/Q(t)$ with $P,Q\in\mathbb{C}[t]$ having no~factor of
  positive degree in~common, and $t$~equals zero at the cusp $\tau={\rm
  i}\infty$ {\rm(}so that $Q(0)=0${\rm)}, then in a neighborhood of the
  point $t=0$ on the quotient\/ $\Gamma\setminus\mathcal{H}^*,$ the
  holomorphic function~$h_{\Gamma,t}$ is defined by
  \begin{equation*}
    h_{\Gamma,t}(t) = \left[P(t)/P(0)\right]^{-1/12}
{}_2F_1\left(\tfrac1{12},\tfrac{5}{12};\,1;\,12^3Q(t)/P(t)\right).
  \end{equation*}
  This definition of~$h_{\Gamma,t}$ is unaffected by the Hauptmodul
  $t$~being replaced by any nonzero scalar multiple, i.e.,
  $h_{\Gamma,t}=h_{\Gamma,\alpha t}$ for any nonzero~$\alpha$.
\label{def:1}
\end{definition}

For useful facts on the Gauss\ function ${}_2F_1,$ see the appendix.  Here
it suffices to know that ${}_2F_1(\frac1{12},\frac{5}{12};1;z)$ equals
unity at~$z=0$ and is holomorphic on the unit disk.  It has a branch point
of square-root type at~$z=1,$ but can be holomorphically extended to
$\mathbb{P}^1(\mathbb{C})_z,$ slit along the positive real axis from $z=1$
to~$z=\infty$.

\begin{definition}
\label{def:hN}
  For notational simplicity, let $\hat h_1:=h_{\Gamma(1),\hat J}$ where
  $\hat J:=1/J=12^3/j$ is the abovementioned alternative Hauptmodul
  for~$\Gamma(1),$ which equals zero at the cusp $\tau={\rm i}\infty$.
  Similarly, let $h_N:=h_{\Gamma_0(N),t_N}$.  So if $j=P_N(t_N)/Q_N(t_N),$
  then
  \begin{align*}
    &\hat h_1(\hat J(\tau)) = {}_2F_1\left(\tfrac1{12},\tfrac5{12};\,1;\,\hat J(\tau)\right),\\
    &h_N(t_N(\tau)) = \left[P_N(t_N(\tau))/P_N(0)\right]^{-1/12}\hat h_1(\hat J(\tau)),
  \end{align*}
  on a neighborhood of the point\/ $\tau={\rm i}\infty$ in~$\mathcal{H}^*$.
\label{def:2}
\end{definition}

By examining Table~\ref{tab:coverings}, one sees that for each~$N,$ the
function $j=j(t_N)$ maps the positive real axis $t_N>0$ into the interval
$12^3\le j<\infty$.  That~is, if ${0<t_N<\infty}$ then $\hat J=12^3/j$
satisfies $0<\hat J\le1$.  Moreover, at each point $t_N=t_N^*>0$ at~which
$j=12^3,$ i.e., $\hat J=1,$ the behavior of $j-12^3$ is quadratic.  (The
point $t_2=8$ on~$X_0(2)$ is an example.)  It follows that $h_N=h_N(t_N),$
which as defined above is holomorphic at~$t_N=0$ and equal to unity there,
has a real holomorphic continuation along the positive real axis $t_N>0$.

A connection with differential equations is made by Theorems \ref{thm:ford}
and~\ref{thm:2} below, which are standard~\cite{Ford51,Stiller88}.  First,
recall some facts on Fuchsian differential operators and equations.
A~second-order operator $L=D_t^2+\mathcal{A}\cdot D_t+\mathcal{B}$
on~$\mathbb{P}^1(\mathbb{C})_t,$ where
$\mathcal{A},\mathcal{B}\in\mathbb{C}(t),$ is said to be Fuchsian if all
its singular points are regular, i.e., if it has two characteristic
exponents $\alpha_{i,1},\alpha_{i,2}\in\mathbb{C}$ (which may be the same)
at each singular point~$\mathfrak{s}_i\in \mathbb{P}^1(\mathbb{C})$.
(These include the poles of~$\mathcal{A},\mathcal{B},$ and perhaps the
infinite point~$t=\infty$.)  In the nondegenerate case
$\alpha_{i,1}-\alpha_{i,2}\not\in\mathbb{Z},$ the existence of two
exponents means the differential equation $L u=0$ has local solutions
$u_{i,j},$ $j=1,2,$ at~$\mathfrak{s}_i$ of the form
$t^{\alpha_{i,j}}$~times an invertible function of~$t,$ where $t$~is a
local uniformizing parameter (if $\alpha_{i,1}-\alpha_{i,2}\in\mathbb{Z},$
one solution may be logarithmic).  The definition of characteristic
exponents extends trivially from singular to ordinary (i.e., non-singular)
points.  Any finite point~$t$ that is ordinary has exponents~$0,1$.

\begin{definition}
  A second-order Fuchsian differential operator $L$
  on\/~$\mathbb{P}^1(\mathbb{C})$ is in normal form if it has a zero
  exponent at each of its finite singular points.
\end{definition}
An example of a normal-form operator is the Gauss\ hypergeometric operator
$L_{a,b;c},$ which has the function ${}_2F_1(a,b;c;\cdot)$ in its kernel.
In~general, any monic second-order Fuchsian operator in normal form will be
of the form~\cite{Poole36}
\begin{equation}
\label{eq:Poole}
\frac{d^2}{dt^2} +
\left[\sum_{i=1}^{n-1}\frac{1-\rho_i}{t-a_i}\right]\cdot\frac{d}{dt}
+\left[\frac{\Pi_{n-3}(t)}{\prod_{i=1}^{n-1}(t-a_i)}\right],
\end{equation}
where $\{a_i\}_{i=1}^{n-1}$ are the finite singular points, with
exponents~$0,\rho_i,$ and $\Pi_{n-3}(t)$ is a degree-$(n-3)$ polynomial.
Its leading coefficient determines the exponents at $t=\infty,$ and its
$n-3$~trailing coefficients are so-called accessory parameters.  The local
monodromy of~$L,$ i.e., the monodromy of the differential equation ${Lu=0}$
around each of its singular points, is determined by the exponents, and its
global monodromy is determined, additionally, by the $n-3$ accessory
parameters.  Since the operator $L_{a,b;c}$ has $n=3,$ it contains
no~accessory parameters: its global monodromy is determined uniquely by its
local monodromy.

Any substitution $u=f^{-\alpha} \hat u,$ where $f\in\mathbb{C}(t)$ and
$\alpha\in\mathbb{C},$ will transform the equation $Lu=0$ to $\hat L\hat
u=0,$ where $\hat L$~has transformed coefficients
$\hat{\mathcal{A}},\hat{\mathcal{B}}\in\mathbb{C}(t)$.  Any such `index
transformation' leaves exponent {\em differences\/} invariant: the
exponents of~$\hat L$ at any finite point $t=t_0$ at~which $f$~has order of
vanishing~$m$ will be those of~$L,$ shifted down by~$\alpha m$.  By index
transformations, any Fuchsian operator on~$\mathbb{P}^1(\mathbb{C})$ may be
reduced to one in normal form.

\begin{theorem}
\label{thm:ford}
  Let\/ $\Gamma$ be a Fuchsian subgroup\/ {\rm(}of the first kind\/{\rm)}
  of the automorphism group\/ ${\it PSL}(2,\mathbb{R})$
  of\/~$\mathcal{H}\ni\tau$ for which the quotient curve\/
  $X:=\Gamma\setminus\mathcal{H}^*$ is of genus zero, and let\/ $t$~denote
  a Hauptmodul.  Then in a neighborhood of any point on\/~$X,$ any branch
  of~$\tau,$ which can be viewed as a multivalued function on\/~$X,$ will
  equal the ratio of two independent solutions of some second-order
  Fuchsian differential equation\/ $Lu:=(D_t^2+\mathcal{A}\cdot D_t +
  \mathcal{B})u=0,$ where\/ $\mathcal{A},\mathcal{B}\in\mathbb{C}(t)${\rm.}
  The space of local solutions of this Picard--Fuchs equation will be\/
  $\left[\mathbb{C}\tau(\cdot)+\mathbb{C}\right]H(\cdot),$ where\/ $H$~is
  some particular local solution.  One can choose\/~$L$ so~that its
  singular points are the fixed points of\/~$\Gamma$ on\/~$X,$ with the
  difference of characteristic exponents equaling~$1/k$ at each fixed
  point of order\/~$k,$ and zero at each parabolic fixed point\/
  {\rm(}i.e., cusp\/{\rm)}{\rm.}
\end{theorem}

This is a special (e.g., genus zero) case of a classical theorem dealing
with Fuchsian automorphic functions of the first kind~\cite[\S44,\
Thm.~15]{Ford51}.  It does not require that $\Gamma$~be a subgroup of
$\Gamma(1)={\it PSL}(2,\mathbb{Z})$.  The following theorem is also
classical.

\begin{theorem}
\label{thm:2}
  For any first-kind Fuchsian subgroup\/~$\Gamma$ and Hauptmodul\/~$t,$ if
  two Picard--Fuchs equations of the form\/ $(D_t^2+\mathcal{A}\cdot D_t +
  \mathcal{B})u=0$ have the same characteristic exponents\/ {\rm(}not
  merely exponent differences{\rm)} at each singular point, then they must
  be equal.
\end{theorem}
\begin{corollary}
  Requiring the Picard--Fuchs equation mentioned in the last sentence of
  Theorem\/~{\rm\ref{thm:ford}} to be in normal form determines it
  uniquely.
\end{corollary}
\begin{proof}
  If any finite point
  $t\in\Gamma\setminus\mathcal{H}^*\cong\mathbb{P}^1(\mathbb{C})$ is a
  cusp, its exponents will be~$0,0$; if it is a quadratic (resp.\ cubic)
  elliptic fixed point, they will be $0,1/2$ (resp.\ $0,1/3$){\rm.}  The
  exponents at~$t=\infty$ are uniquely determined by Fuchs's relation: the
  sum of of all $2k$~characteristic exponents of any second-order Fuchsian
  differential equation with $k$~singular points
  on~$\mathbb{P}^1(\mathbb{C})$ must equal~${k-2}$.
\end{proof}

\begin{theorem}
  \label{thm:3}
  If\/ $\Gamma=\Gamma(1),$ so that\/ $X=X(1),$ and the Hauptmodul\/~$t$
  equals\/~$\hat J,$ then in a neighborhood of the cusp\/ $\hat J=0$
  {\rm(}i.e., $\tau={\rm i}\infty${\rm)}, the unique normal-form
  Picard--Fuchs equation is the Gauss hypergeometric equation with
  parameters\/ $a=\frac1{12},\allowbreak{b}=\frac{5}{12},\allowbreak{c}=1$
  and independent variable\/~$\hat J,$ i.e.,
  \begin{equation*}
    \hat{\mathcal{L}}_1\, u:= L_{\frac1{12},\frac5{12};1}u=
    \left\{D_{\!\hat J}^2
    + \left[ \tfrac1{\hat J} + \tfrac1{2(\hat J-1)} \right] D_{\!\hat J}
    + \tfrac{5/144}{\hat J(\hat J-1)}\right\} u = 0.
  \end{equation*}
  The fundamental local solution\/~$H=H(\hat J)$ can be taken to be\/~$\hat
  h_1(\hat J)$ as defined above, i.e.,
  ${}_2F_1(\frac1{12},\frac5{12};1;\hat J),$ the unique local solution of\/
  $\hat{\mathcal{L}}_1u=0$ which is holomorphic at\/~$\hat J=0,$ up~to
  normalization.  Also, the following connection to the theory of modular
  forms exists: in a neighborhood of the cusp\/ $\tau={\rm i}\infty,$
  \begin{align*}
    &E_4(\tau)=\hat h_1^4(\hat J(\tau)),\\
    &E_6(\tau)=[1-\hat J(\tau)]^{1/2}\,\hat h_1^6(\hat J(\tau)),\\
    &\Delta(\tau)=(2\pi)^{12}12^{-3}(E_4^3-E_6^2)(\tau)=(2\pi)^{12}12^{-3}\hat J(\tau)\,\hat h_1^{12}(\hat J(\tau)),
  \end{align*}
  where\/ $E_4,E_6,\Delta$ are the classical Eisenstein sums and modular
  discriminant.
\end{theorem}
\begin{proof}
  The Fuchsian equation $\hat{\mathcal{L}}_1u=0$ is in normal form, with
  exponents $0,0$ at~$\hat J=0,$ $0,\frac12$ at~$\hat J=1,$ and\/
  $\frac{1}{12},\frac{5}{12}$ at~$\hat J=\infty$; so its exponent
  differences are $0$~at the cusp $\tau={\rm i}\infty,$ $\frac12$~at the
  quadratic elliptic point~$\tau={\rm i},$ and $\frac13$~at the cubic
  elliptic point~$\tau=\zeta_3,$ in agreement with Theorem~\ref{thm:ford}.
  Up~to shifts of exponents, any hypergeometric equation is uniquely
  determined by its exponent differences, since it contains no~accessory
  parameters; hence the first sentence of the theorem follows.  The
  remaining two are due to Stiller~\cite{Stiller88} (another derivation of
  the expression for~$\Delta$ in~terms of $\hat J$ and~$\hat h_1$ will be
  mentioned below).
\end{proof}

\begin{theorem}
\label{thm:4}
  If\/ $\Gamma=\Gamma_0(N),$ so that\/ $X=X_0(N),$ and the Hauptmodul\/~$t$
  equals\/~$t_N,$ then in a neighborhood of the cusp\/ $t_N=0$ {\rm(}i.e.,
  $\tau={\rm i}\infty${\rm)}, the unique normal-form Picard--Fuchs
  equation\/ $\mathcal{L}_Nu=0$ has
  \begin{enumerate}
    \item one singular point with characteristic exponents\/
    $\frac1{12}\psi(N),\frac1{12}\psi(N),$ viz., the cusp\/ $t_N=\infty$
    {\rm(}i.e., $\tau=0${\rm);}
    \item $\sigma_\infty(N)-1$ singular points with exponents\/ $0,0,$
    viz., the remaining cusps, including\/ $t_N=0$ {\rm(}i.e., $\tau={\rm
    i}\infty${\rm);}
    \item $\epsilon_2(N)$ singular points with exponents\/ $0,\frac12,$ viz.,
    the order-$2$ elliptic fixed points;
    \item $\epsilon_3(N)$ singular points with exponents\/ $0,\frac13,$
    viz., the order-$3$ elliptic fixed points;
  \end{enumerate}
  and is the equation on\/ $X_0(N)$ obtained by\/ {\rm(i)}~pulling back\/
  $\hat{\mathcal{L}}_1u=0$ to $X_0(N)$ along the covering map\/
  $j=P_N(t_N)/Q_N(t_N)$ of $X_0(N)/X(1),$ and\/ {\rm(ii)}~performing the
  substitution\/ $\hat u=P_N(t_N)^{-1/12}\, u$.  The local solution\/
  $H=H(t_N)$ can be taken to be\/ $h_N=h_N(t_N)$ as defined above.  It is
  the unique local solution of\/ $\mathcal{L}_N u=0$ that is holomorphic at
  the cusp\/ $t_N=0$ and equals unity there.
\end{theorem}
\begin{proof}
  Pulling back from $X(1)$ to $X_0(N)$ and performing the indicated
  substitution will not remove the property that the Picard--Fuchs equation
  should have; namely, that any branch of~$\tau$ should equal the ratio of
  two of its solutions.  So, all that needs to be proved are the statements
  about the exponents of the resulting operator~$\mathcal{L}_N$; and also
  the final two sentences of the theorem.

  At any point $t_N\in X_0(N)$ at which $\hat J=\hat J(t_N) =
  12^3Q_N(t_N)/P_N(t_N),$ the covering map, has ramification index~$k,$ the
  pullback $(\hat{\mathcal{L}}_1)^*$ of~$\hat{\mathcal{L}}_1$ has exponents
  equal to $k$~times those of $\hat{\mathcal{L}}_1$ at~$\hat J(t_N)$.  So
  at each cusp of~$X_0(N),$ i.e., at each point on the fibre above $\hat
  J=0,$ the pulled-back operator will have exponents~$0,0$.

  The points on $X_0(N)$ above $\hat J=1$ are partitioned into the
  $\epsilon_2(N)$ order-$2$ elliptic fixed points of~$\Gamma_0(N)$
  (at~which $k=1$), and non-fixed points at which
  ${k=2}$~\cite{Schoeneberg74}.  The corresponding exponents will be
  $0,\frac12$ and~$0,1,$ so the latter will be ordinary (non-singular)
  points of~$(\hat{\mathcal{L}}_1)^*$.  Similarly, the points on~$X_0(N)$ above $\hat
  J=0$ are partitioned into the $\epsilon_3(N)$ order-$3$ elliptic fixed
  points of~$\Gamma_0(N)$ (at~which $k=1$) and non-fixed points
  at~which $k=3$.  The corresponding exponents will be
  $\frac1{12},\frac5{12}$ and~$\frac14,\frac54$.  Since the
  polynomial~$P_N$ has a simple root at each of the former and a triple
  root at each of the latter, performing the substitution $\hat
  u=P_N(t_N)^{-1/12}u$ will shift these exponents to $0,\frac13$ and~$0,1$
  respectively; so the non-fixed points will become ordinary.  It will
  also, since $\deg P_N=\psi(N),$ shift the exponents at~$t_N=\infty$ from
  $0,0$ to~$\frac1{12}\psi(N),\frac1{12}\psi(N)$.

  If the operator~$\mathcal{L}_N$ resulting from the substitution is taken
  to be monic, then $\mathcal{L}_N=P_N(\cdot)^{-1/12}
  (\hat{\mathcal{L}}_1)^* P_N(\cdot)^{1/12}$.  But by the definition of a
  pullback,
  \begin{equation}
    (\hat{\mathcal{L}}_1)^*\left[{}_2F_1(\tfrac1{12},\tfrac5{12};\,1;\,\hat J(\cdot))\right]= 0.
  \end{equation}
  So
  \begin{equation}
    \mathcal{L}_N[h_N] = \mathcal{L}_N\left[ P_N^{-1/12}(\cdot)\,\, {}_2F_1(\tfrac1{12},\tfrac5{12};\,1;\,\hat
      J(\cdot)) \right]=0.
  \end{equation}
The final sentence of the theorem now follows by the general theory of
Fuchsian differential equations: since $\mathcal{L}_N$~has exponents $0,0$
at~$t_N=0,$ its space of holomorphic solutions there is one-dimensional,
and must accordingly be~$\mathbb{C}h_N$.  The penultimate sentence of the
theorem follows by examination.
\end{proof}

Theorem~\ref{thm:3} reveals why the prefactor $[P(t)/P(0)]^{-1/12}$ was
included in Definition~\ref{def:1}.  If it were absent, then the
Picard--Fuchs equation satisfied by~$h_N$ would not be in normal form.

In the theory of conformal mapping, there is another type of normal-form
Picard--Fuchs equation in widespread use~\cite{Ford51,Harnad2000}.  If
$t$~is a Hauptmodul for a first-kind Fuchsian subgroup of ${\it
PSL}(2,\mathbb{R})$ as in Theorem~\ref{thm:ford}, any branch of~$\tau$ will
equal the ratio of two independent solutions of the second-order Fuchsian
differential equation $\left[D_t^2+\mathcal{Q}(t)\right]v=0,$ where
$\mathcal{Q}\in\mathbb{C}(t)$ is defined by
$\mathcal{Q}:=\frac12\{\tau,t\},$ with $\{\cdot,\cdot\}$ the Schwarzian
derivative.  This is a {\em self-adjoint\/} Picard--Fuchs equation, and its
space of local solutions is
$\left[\mathbb{C}\tau(\cdot)+\mathbb{C}\right](dt/d\tau)^{1/2}(\cdot)$.
The reader may wonder why we introduced a non-self-adjoint normal form,
instead.  It is because a Picard--Fuchs equation in our normal form is more
convenient for deriving special function identities.  Though asymmetric, it
does permit an elegant modular interpretation of its fundamental
solution~$h_N,$ as the following theorem and corollaries indicate.

\begin{theorem}
\label{thm:hNformula}
  For each\/ $N$ with\/ $\Gamma_0(N)$ of genus zero, if\/ $h_N=h_N(t_N)$ is
  the holomorphic function defined above on a neighborhood of the point\/
  $t_N=0,$ and\/ $\mathfrak{h}_N=\mathfrak{h}_N(\tau)$ is defined on a
  neighborhood of the infinite cusp\/ $\tau={\rm i}\infty$ by\/
  $\mathfrak{h}_N(\tau):=h_N(t_N(\tau)),$ then\/ $\mathfrak{h}_N$~extends
  to\/~${\mathcal{H}}^*\ni \tau$ by continuation, yielding a weight\/-$1$
  modular form for\/~$\Gamma_0(N),$ with some multiplier system.  In
  particular,
  \begin{equation}
    \label{eq:hN}
    \mathfrak{h}_N(\tau):= h_N(t_N(\tau)) = P_N(0)^{1/12} Q_N(t_N(\tau))^{-1/12}\,\eta^2(\tau).
  \end{equation}
  This modular form is regular and non-vanishing at each cusp in\/
  $\mathbb{P}^1(\mathbb{Q})$ other than those in the class\/
  $\bigl[\frac11\bigr]_N\ni0,$ at each of which its order of vanishing
  is\/~$\psi(N)/12N$.
\end{theorem}
\begin{proof}
  In a neighborhood of the infinite cusp, it is a striking fact that the
  holomorphic function $\hat h_1={}_2F_1(\frac1{12},\frac5{12};1;\hat J)$
  of~$\hat J$ can be expressed in~terms of the $j$-invariant and the eta
  function, as $12^{1/4}\hat J^{-1/12}\,\eta^2$.  This representation was
  known to Dedekind~\cite[p.~137]{Chandrasekharan85}, and was rediscovered
  by Stiller~\cite{Stiller88}.  Formula (\ref{eq:hN}) follows from~it,
  since $\hat J(\tau)=12^3 Q_N(t_N(\tau))/P_N(t_N(\tau))$.  (Incidentally,
  the representation for the discriminant $\Delta=\Delta(\tau)$ in~terms of
  $\hat J$ and~$\hat h_1$ given in Theorem~\ref{thm:3} also follows
  from~it, since $\Delta=(2\pi)^{12}\,\eta^{24}$.)  Since the factors
  $Q_N(t_N(\tau))^{-1/12}$ and~$\eta^2(\tau)$ are holomorphic
  on~$\mathcal{H}\ni\tau,$ $h_N(t_N(\tau))$ can be continued from the
  neighborhood of~$\tau={\rm i}\infty$ on~which it was originally defined,
  to all of~$\mathcal{H}$.  As a modular form it has weight~$1,$ since
  $\eta$~has weight~$1/2$.

  The roots of the polynomial~$Q_N$ correspond to the cusps of
  $\Gamma_0(N)$ on $X_0(N)$ other than $t_N=\infty,$ i.e., to the cusp
  classes $\bigl[\frac{a}d\bigr]_N\subset\mathbb{P}^1(\mathbb{Q})$ other
  than~$\bigl[\frac11\bigr]_N$.  Any root $t_N=t^*$ corresponding
  to~$\bigl[\frac{a}d\bigr]_N$ appears in~$Q_N$ with multiplicity equal to
  the cusp width $e_{d,N}=N/d(N,N/d)$.  But this is also the multiplicity
  with which $t_N=t^*$ is mapped to~$X(1)\cong\mathbb{P}^1(\mathbb{C})_j$
  by the covering map $j=P_N(t_N)/Q_N(t_N)$.  So the order of vanishing of
  $Q_N(t_N(\tau))$ at~any cusp not in~$\bigl[\frac11\bigr]_N$ is unity, and
  (as the order of~$\eta(\tau)$ at any cusp
  $\tau\in\mathbb{P}^1(\mathbb{Q})$ is~$\frac1{24}$), the order of
  $Q_N(t_N(\tau))^{-1/12}\,\eta^2(\tau)$ will be
  $-\frac1{12}+2\cdot\frac1{24}=0$.  Definition~\ref{def:hN} implies that
  at any cusp in~$\bigl[\frac11\bigr]_N,$ the order of~$h_N(t_N(\tau))$
  will be $\frac1{12}(\deg P)/e_{1,N}=\psi(N)/12N$.
\end{proof}

\begin{corollary}
\label{cor:hNformula1}
$t_N,\mathfrak{h}_N$ have divisors\/
$\bigl(\bigl[\frac1N\bigr]_N\bigr)-\frac1N\cdot\bigl(\bigl[\frac11\bigr]_N\bigr),$
$\frac{\psi(N)}{12N}\cdot\bigl(\bigl[\frac11\bigr]_N\bigr),$ if viewed as
functions on\/ $\mathcal{H}^*\ni\tau$.  If viewed respectively as a
univalent and a multivalued function on\/~$X_0(N),$ they have divisors\/
$(t_N=0)-(t_N=\infty)$ and $\frac{\psi(N)}{12}\cdot(t_N=\infty)$.
\end{corollary}
\begin{remark*}
The term `divisor' is used in a generalized sense here, referring to an
element of a free $\mathbb{Q}$-module, rather than a free
$\mathbb{Z}$-module.  The coefficients of a generalized divisor are orders
of vanishing, which may not be integers.
\end{remark*}
\begin{proof}
  That the divisor of $\mathfrak{h}_N(\tau)=h_N(t_N(\tau))$ is
  $\frac{\psi(N)}{12N}\cdot\bigl(\bigl[\frac11\bigr]_N\bigr)$ follows from
  the last sentence of the theorem, since $h_N$~has no zeroes or poles
  on~$\mathcal{H}$.  That $t_N,$ as a function on~$X_0(N),$ has the stated
  divisor is trivial.  The remaining statements follow by taking the cusp
  widths $e_{1,N}=N$ and $e_{N,N}=1$ into account.
\end{proof}

\begin{corollary}
\label{cor:hNformula2}
  For each\/ $N$ with\/ $\Gamma_0(N)$ of genus zero, the weight\/-$1$
  modular form\/ $\mathfrak{h}_N(\tau)=h_N(t_N(\tau))$ for\/~$\Gamma_0(N)$
  has the alternative representation
  \begin{multline*}
    \prod_{i=1}^{\epsilon_2(N)} \Biggl[1-\frac{t_N(\tau)}{t_{N,i}^{(2)}}\Biggr]^{-1/4}
    \,
    \prod_{i=1}^{\epsilon_3(N)} \Biggl[1-\frac{t_N(\tau)}{t_{N,i}^{(3)}}\Biggr]^{-1/3}
    \,
    \prod_{i=1}^{\sigma_\infty(N)-2} \Biggl[1-\frac{t_N(\tau)}{t_{N,i}^{(\infty)}}\Biggr]^{-1/2}\\
    {}\times t_N^{-1/2}(\tau)\left[\frac{1}{2\pi{\rm i}}\,\frac{dt_N}{d\tau}(\tau)\right]^{1/2},
  \end{multline*}
  in which the three products run over\/ {\rm(i)}~the quadratic fixed
  points $t_N=t_{N,i}^{(2)},$ $i=1,\dots,\epsilon_2(N),$ {\rm(ii)}~the
  cubic fixed points $t_N=t_{N,i}^{(3)},$ $i=1,\dots,\epsilon_3(N),$ and\/
  {\rm(iii)}~the cusps $t_N=t_{N,i}^{(\infty)},$
  $i=1,\dots,\sigma_\infty(N)-2,$ other than the two distinguished cusps
  $t_N=0$ {\rm(}i.e., $\tau={\rm i}\infty${\rm)} and $t_N=\infty$
  {\rm(}i.e., $\tau=0${\rm)}.  These fixed points on~$X_0(N)$ are given in
  Table\/~{\rm\ref{tab:hauptmoduln}}.
\end{corollary}
\begin{remark*}
  This result clarifies the connection between two Picard--Fuchs equations
  for $\Gamma_0(N)$: the self-adjoint one noted above, of the form
  $[D_{t_N}^2+\mathcal{Q}_N(t_N)]v=0,$ and the normal-form one
  $\mathcal{L}_Nu=0$.  They are related by a similarity transformation, and
  in~fact are projectively equivalent, since their respective solution
  spaces are ${[\mathbb{C}\tau(\cdot)+\mathbb{C}]}(dt_N/d\tau)^{1/2}(\cdot)$
  and $[\mathbb{C}\tau(\cdot)+\mathbb{C}]h_N(\cdot)$.  The given
  representation for $h_N(t_N(\tau))$ yields a formula for~$u/v,$ the
  quotient of their solutions.
\end{remark*}
\begin{proof}
  Despite the disconcerting presence of fractional powers, the given
  expression extends from a neighborhood of $\tau={\rm i}\infty$ to a
  single-valued function of~$\tau\in\mathcal{H},$ without zeroes or poles.
  To see this, begin by examining $\tau=\zeta_3,$ at which
  $t_N=t_{N,i}^{(3)},$ some cubic fixed point on~$X_0(N)$.  To leading
  order one has $j\sim C(\tau-\zeta_3)^3$ near $\tau=\zeta_3,$ and also
  $j\sim C'(t_N-t_{N,i}^{(3)})$.  So $[dt_N/dt]^{1/2}\sim
  C''(\tau-\zeta_3)$ and $[1-t_N/t_{N,i}^{(3)}]^{-1/3}\sim
  C'''(\tau-\zeta_3)^{-1},$ implying finiteness and single-valuedness of
  their product near~$\tau=\zeta_3$.  This extends to points
  $\tau\in\mathcal{H}$ congruent to~$\zeta_3$.  Quadratic fixed points can
  be handled similarly.

  The expression also has a nonzero, non-infinite limit as any cusp
  $\tau\in\mathbb{P}^1(\mathbb{Q})$ not congruent to~$\tau=0$
  under~$\Gamma_0(N)$ is approached.  The limit as~$\tau\to{\rm i}\infty$
  is easily seen to be unity, since $t_N\sim C^{({\rm iv})}q$ to leading
  order.  Cusps $\tau\in\mathbb{P}^1(\mathbb{Q})$ `above' the
  $\sigma_\infty(N)-2$ parabolic fixed points $t_{N,i}^{(\infty)},$ i.e.,
  cusps not congruent to $\tau={\rm i}\infty$ or~$\tau=0,$ are almost as
  easy to handle: one expands to leading order, as in the last paragraph,
  and in each case finds a nonzero, non-infinite limit.

  When $\tau\to0,$ $t_N\to\infty,$ yielding leading-order behavior $\sim
  C^{({\rm v})} t_N^{-a},$ where the exponent~$a$ equals
  $\frac14\epsilon_2(N)+\frac13\epsilon_3(N)+\frac12\sigma_\infty(N)-1$.
  From the genus formula~(\ref{eq:genus}), one deduces
  $a=\frac1{12}\psi(N)$.  By the preceding corollary, the quotient of the
  given expression and $h_N(t_N(\tau))$ must tend to a nonzero,
  non-infinite limit as any cusp congruent to~$\tau=0$ is approached; and
  by the last paragraph, the same is true of any other cusp.  The quotient
  must therefore be a nonzero constant.  By considering the limit
  $\tau\to{\rm i}\infty,$ one sees that it equals unity.
\end{proof}

It is natural to wonder whether the function $h_N=h_N(t_N)$ alone, rather
than the two-dimensional space of solutions
$[\mathbb{C}\tau(\cdot)+\mathbb{C}]h_N(\cdot)$ of the Picard--Fuchs
equation, supplies a means of computing the multivalued function~$\tau$
on~$X_0(N)$.  It is easiest to consider the half-line $t_N>0,$ since
$h_N$~extends holomorphically to~it.  For each of the 14~values of~$N,$
$t_N(\tau)\sim\kappa_N\cdot q$ to leading order as~$\tau\to{\rm i}\infty,$ and
accordingly $\tau\sim\frac{{\rm i}}{2\pi}\log(1/t_N)$ as~$t_N\to0^+$.  By
choosing the principal branch of the logarithm one obtains a unique
function $\tau=\tau(t_N),$ $t_N>0$.  It~is evident that $\tau/\,{\rm i}$~is
positive and monotone decreasing on~$t_N>0$.  In~fact it decreases to zero
as~$t_N\to\infty,$ since $\tau$~must satisfy
$\tau(\kappa_N/t_N)=-1/[N\tau(t_N)],$ which follows from the identity
$t_N(\tau)t_N(-1/N\tau)=\kappa_N$ introduced
in~\S\ref{sec:parametrized1}.

Theorem~\ref{thm:clever} below reveals how the function $\tau=\tau(t_N),$
$t_N>0,$ can be computed from the corresponding function $h_N=h_N(t_N)$ on
the half-line.

\begin{lemma}
  \label{lem:lone}
  For each\/ $N,$ the Picard--Fuchs equation\/ $\mathcal{L}_N u=0$ is
  invariant under the substitution\/ $\tilde u(t_N)=t_N^{-\psi(N)/12}
  u_N(\kappa_N/t_N)$.  That~is, $\mathcal{L}_N\tilde u=0$.
\end{lemma}
\begin{proof}
  Suppose that applying the substitution to $\mathcal{L}_N u=0$ yields a
  transformed equation $\tilde{\mathcal{L}}_N \tilde u=0$.  Since
  $t_N(\tau)t_N(-1/N\tau)=\kappa_N,$ the transformed equation has the
  defining property of a Picard--Fuchs equation; namely, that any branch of
  the multivalued function~$\tau$ is locally a quotient of two of its
  solutions.  But the map $t_N\mapsto \kappa_N/t_N$ (i.e.,
  $\tau\mapsto-1/N\tau$) is the Fricke involution, which separately
  involutes cusps, quadratic elliptic points, and cubic elliptic points.
  So the transformed operator $\tilde{\mathcal{L}}_N=D_{t_N}^2 +
  \tilde{\mathcal{A}}(t_N)D_{t_N} + \tilde{\mathcal{B}}(t_N)$ must have,
  like~$\mathcal{L}_N,$ exponents~$0,0$ at each cusp other
  than~$t_N=\infty,$ exponents $0,\frac12$~at each quadratic point, and
  $0,\frac13$ at each cubic point.  (The factor $t_N^{-\psi(N)/12}$ adjusts
  the exponents at~$t_N=0$ to be~$0,0,$ and those at~$t_N=\infty$ to be
  $\frac1{12}\psi(N),\frac1{12}\psi(N)$.)  Since the exponents
  of~$\tilde{\mathcal L}_N$ and~$\mathcal{L}_N$ agree,
  $\tilde{\mathcal{L}}_N=\mathcal{L}$ by Theorem~\ref{thm:2}.
\end{proof}

\begin{theorem}
  \label{thm:clever}
  For each\/ $N$ with\/ $\Gamma_0(N)$ of genus zero, on\/ $t_N>0$ one has
  that
  \begin{equation*}
    \tau(t_N)/\,{\rm i} = A_N\,
    t_N^{-\psi(N)/12}\,\,\frac{h_N(\kappa_N/t_N)}{h_N(t_N)},
  \end{equation*}
  for some\/ $A_N>0$.  Here\/ $h_N,\tau/\,{\rm i}$ are the single-valued
  holomorphic continuations introduced above, which are real and positive
  on the half-line\/~$t_N>0$.
\end{theorem}
\begin{proof}
  Let $\tau_1:=\tau h_N$ and $\tau_2:=h_N$.  These are solutions of the
  Picard--Fuchs equation $\mathcal{L}_N u=0,$ with $\tau=\tau_1/\tau_2$.
  Since $\mathcal{L}_N$~has exponents $0,0$ at~$t_N=0,$ solutions of
  $\mathcal{L}_N u=0$ must be asymptotic to ${\rm const} \times \log
  (1/t_N)$ or to~${\rm const},$ as~$t_N\to0^+$.  The solutions
  $\tau_1,\tau_2$ are of these two types, respectively.

  Similarly, since $\mathcal{L}_N$~has exponents
  $\frac1{12}\psi(N),\frac1{12}\psi(N)$ at~$t_N=\infty,$ solutions must be
  asymptotic to ${\rm const}\times t_N^{-\psi(N)/12}\log t_N$ or to ${\rm
  const}\times t_N^{-\psi(N)/12}$ as $t_N\to\infty$.  Consider
  in~particular $\tilde\tau_i(t_N):=t_N^{-\psi(N)/12}\tau_i(\kappa_N/t_N),$
  for $i=1,2,$ which by Lemma~\ref{lem:lone} are solutions.  Their ratio
  $[\tilde\tau_1/\tilde\tau_2](t_N)=\tau(\kappa_N/t_N)=-1/[N\tau(t_N)]$
  converges to zero as~$t_N\to0^+,$ so it must be the case that
  $\tilde\tau_1(t_N)\sim{\rm const}$ and $\tilde\tau_2(t_N)\sim{\rm
  const}\times \log(1/t_N),$ as~$t_N\to0^+$.  Hence $\tau_1(t_N)\sim{\rm
  const}\times t_N^{-\psi(N)/12}$ and $\tau_2(t_N)\sim {\rm const}\times
  t_N^{-\psi(N)/12}\log t_N$ as~$t_N\to\infty$.

  For $\tau(\kappa_N/t_N)=-1/[N\tau(t_N)]$ to be true, it must in~fact be
  the case that
  \begin{align}
    \tau_1(t_N)&={\rm i}A_N\,\tilde\tau_2(t_N) = {\rm i}A_N\, t_N^{-\psi(t_N)/2}\,\tau_2(\kappa_N/t_N),\\
    \tau_2(t_N)&=B_N\,\tilde\tau_1(t_N) = B_N\, t_N^{-\psi(t_N)/2}\,\tau_1(\kappa_N/t_N),
  \end{align}
  with $A_N=B_N/N,$ where by examination $A_N,B_N$ are real and positive.
  Since $\tau=\tau_1/\tau_2$ and $\tau_2=h_N,$ the claimed identity
  follows.
\end{proof}

\section{The Modular Form $\mathfrak{h}_N(\tau)$: Explicit Computations~(I)}
\label{sec:modular2}

In \S\ref{sec:modular}, the canonical weight-$1$ modular form
$\mathfrak{h}_N(\tau)=h_N(t_N(\tau))$ was defined for each genus-zero
group~$\Gamma_0(N)$.  In this section, we make $\mathfrak{h}_N(\tau)$ much
more concrete by (i)~expressing it as an eta product, (ii)~determining
(with much labor!)\ its multiplier system, and (iii)~working~out its
$q$-expansion.  As a byproduct of the last, we obtain interesting
combinatorial identities that extend those of Fine~\cite{Fine88}.

For each
$\left(\begin{smallmatrix}a&b\\c&d\end{smallmatrix}\right)\in\bar\Gamma_0(N),$
the inhomogeneous version of~$\Gamma_0(N),$ one must have ($h_N\circ
t_N)\left(\frac{a\tau+b}{c\tau+d}\right)=\hat\chi_N(a,b,c,d)\cdot(c\tau+d)\,(h_N\circ
t_N)(\tau),$ where the function $\hat \chi_N:\bar\Gamma_0(N)\to U(1)$ is
the multiplier system.  The simplest case is when $\hat\chi_N$ depends only
on~$d,$ i.e., $\hat\chi_N(a,b,c,d)=\chi_N(d),$ and $\chi_N$~has period~$N$.
In this case $\chi_N$~is a homomorphism of $(\mathbb{Z}/N\mathbb{Z})^*$
to~$U(1),$ i.e., is a Dirichlet character to the modulus~$N$.  It~must be
an odd function, since the M\"obius transformation
$\tau\mapsto\frac{a\tau+b}{c\tau+d}$ is unaltered by negating the matrix
$\left(\begin{smallmatrix}a&b\\c&d\end{smallmatrix}\right)$.  So, e.g.,
$\chi_N(-1)$ must equal~$-1$.

\begin{table}
\caption{The canonical Hauptmodul $t_N$ and weight-$1$ modular
form~$\mathfrak{h}_N$ for~$\Gamma_0(N)$; and if the multiplier system
of~$\mathfrak{h}_N$ is given by an (odd, real) Dirichlet character mod~$N,$
its conductor.}
\begin{center}
{\small
\begin{tabular}{cllcc}
\hline\noalign{\smallskip}
$N$ & \hfil $t_N(\tau)=\kappa_N\cdot\hat t_N(\tau)$ \hfil \hfil \hfil \hfil & $\mathfrak{h}_N(\tau)=h_N(t_N(\tau))$ & ${\rm cond}(\chi_N)$ \\
\noalign{\smallskip}\hline\noalign{\smallskip}
$2$ & $2^{12}\cdot[2]^{24}/\,[1]^{24}$ & $[1]^4/\,[2]^2\vphantom{\bigl\{[1]^5/\,[5]\bigr\}^{1/2}}$ & --- \\
$3$ & $3^{6}\cdot[3]^{12}/\,[1]^{12}$ & $[1]^3/\,[3]\vphantom{\bigl\{[1]^5/\,[5]\bigr\}^{1/2}}$ & --- \\
$4$ & $2^{8}\cdot[4]^{8}/\,[1]^{8}$ & \ \ \ \ [$=h_2(t_2(\tau))\vphantom{\bigl\{[1]^5/\,[5]\bigr\}^{1/2}}$] & --- \\
$5$ & $5^{3}\cdot[5]^{6}/\,[1]^{6}$ & $\bigl\{[1]^5/\,[5]\bigr\}^{1/2}$ & --- \\
$6$ & $2^{3} 3^{2}\cdot [2][6]^{5}/\,[1]^{5}[3]$ & $[1]^6[6]\,/\,[2]^3[3]^2\vphantom{\bigl\{[1]^5/\,[5]\bigr\}^{1/2}}$ & $3$ \\
$7$ & $7^{2}\cdot[7]^{4}/\,[1]^{4}$ & $\bigl\{[1]^7/\,[7]\bigr\}^{1/3}$ & --- \\
$8$ & $2^{5}\cdot[2]^{2}[8]^{4}/\,[1]^{4}[4]^{2}$ & \ \ \ \ [$=h_2(t_2(\tau))\vphantom{\bigl\{[1]^5/\,[5]\bigr\}^{1/2}}$] & $4$ & \\
$9$ & $3^{3}\cdot[9]^{3}/\,[1]^{3}$ & \ \ \ \ [$=h_3(t_3(\tau))\vphantom{\bigl\{[1]^5/\,[5]\bigr\}^{1/2}}$] & $3$ \\
$10$ & $2^{2}5\cdot[2]\,[10]^{3}/\,[1]^{3}[5]$ & $\bigl\{[1]^{10}[10]\,/\,[2]^5[5]^2\bigr\}^{1/2}$ & --- \\
$12$ & $2^{2}3\cdot[2]^2[3]\,[12]^3/\,[1]^{3}[4]\,[6]^2$ & \ \ \ \ [$=h_6(t_6(\tau))\vphantom{\bigl\{[1]^5/\,[5]\bigr\}^{1/2}}$] & $3$ \\
$13$ & $13\cdot[13]^{2}/\,[1]^{2}$ & $\bigl\{[1]^{13}/\,[13]\bigr\}^{1/6}$ & --- \\
$16$ & $2^{3}\cdot[2]\,[16]^{2}/\,[1]^{2}[8]$ & \ \ \ \ [$=h_2(t_2(\tau))\vphantom{\bigl\{[1]^5/\,[5]\bigr\}^{1/2}}$] & $4$ \\
$18$ & $2\cdot3\cdot[2]\,[3]\,[18]^2/\,[1]^{2}[6]\,[9]$ & \ \ \ \ [$=h_6(t_6(\tau))\vphantom{\bigl\{[1]^5/\,[5]\bigr\}^{1/2}}$] & $3$ \\
$25$ & $5\cdot[25]\,/\,[1]$ & \ \ \ \ [$=h_5(t_5(\tau))\vphantom{\bigl\{[1]^5/\,[5]\bigr\}^{1/2}}$] & --- \\
\noalign{\smallskip}\hline
\end{tabular}
}%
\end{center}
\label{tab:xNhN}
\end{table}

What remains to be determined is whether the multiplier system of
each~$\mathfrak{h}_N$ is in~fact specified by a Dirichlet
character~$\chi_N$; which may be imprimitive, with its fundamental period,
or conductor, equal to a proper divisor of~$N$.

\begin{theorem}
  For each\/~$N$ with\/ $\Gamma_0(N)$ of genus zero, the eta product
  representations of the canonical Hauptmodul\/~$t_N$ and weight\/-$1$
  modular form\/~$\mathfrak{h}_N$ are as shown in
  Table\/~{\rm\ref{tab:xNhN}}{\rm.}  When\/ $N=6,8,9,12,16,18,$ the
  multiplier system of\/~$\mathfrak{h}_N$ is given by the unique odd,
  $\pm1$-valued Dirichlet character with the stated conductor.  When\/
  $N$~is one of the other eight values, it is not of Dirichlet type.
\end{theorem}
\begin{remark*}
  If $N=6,9,12,18,$ then
  $\chi_N(d)=\bigl(\frac{-3}d\bigr)=\bigl(\frac{d}3\bigr)$; and if
  $N=8,16,$ then
  $\chi_N(d)=\bigl(\frac{-4}d\bigr)=\bigl(\frac{-1}d\bigr)$.  The
  non-Dirichlet cases are $N=2,3,4$; and
  $N=5,7,\allowbreak10,\allowbreak13,25,$ when the eta product for~$h_N$
  contains fractional powers.
\end{remark*}
\begin{proof}
The eta products for the~$t_N$ are taken from Table~\ref{tab:hauptmoduln},
but those for the~$\mathfrak{h}_N$ are new.  They come from
Theorem~\ref{thm:hNformula}.  The polynomial $Q_N(t_N)$ in that theorem is
simply the denominator of the rational expression for~$j$ in~terms
of~$t_N,$ given in Table~\ref{tab:coverings}.  If $\Gamma_0(N)$ has no
irrational cusps then $Q_N$~factors over~$\mathbb{Q}$ into linear factors,
eta products for which follow from Tables \ref{tab:hauptmoduln}
and~\ref{tab:alt}.  For example, $Q_4(t_4)={t_4(t_4+16)},$ and by Tables
\ref{tab:hauptmoduln} and~\ref{tab:alt}, on~$\mathcal{H}^*$ the functions
$t_4,t_4+16$ have respective eta-representations $2^8\cdot[4]^8/\,[1]^8,$
$2^4\cdot[2]^{24}/\,[1]^{16}[4]^8,$ yielding $h_4(t_4(\tau))=[1]^4/\,[2]^2$
as claimed.

When $\Gamma_0(N)$ has irrational cusps, i.e., when $N=9,16,18,25,$
one must work harder, since $Q_N(t_N)$ contains at~least one higher-degree
factor.  For all but $N=25$ these are bivalent functions on $X_0(N)$ with
zeroes at a conjugate pair of irrational cusps.  Fortunately, an eta
product for each factor follows from Table~\ref{tab:intermediate}.  For
example, the factor $t_9^2+9t_9+27$ of~$Q_9(t_9)$ equals~$t_3/t_9,$ i.e.,
$3^3\cdot[3]^{12}/\,[1]^9[9]^3$.  In~all, the multivalent functions that
appear are:
\begin{gather}
t_9^2+9t_9+27 = 3^3\cdot[3]^{12}/\,[1]^9[9]^3,\label{eq:gatherstart}\\
t_{16}^2+4t_{16}+8 = 2^3\cdot [4]^{10}/\,[1]^4[2]^2[8]^4,\label{eq:special16}\\
t_{18}^2+6t_{18}+12 = 2^23\cdot [6]^6[9]\,/\,[1]^3[3]^2[18]^2,\label{eq:special18a}\\
t_{18}^2+3t_{18}+3 = 3\cdot [2]^3[3]^8[18]\,/\,[1]^6[6]^4[9]^2,\label{eq:special18b}\\
t_{25}^4+5t_{25}^3+15t_{25}^2+25t_{25}+25 = 5^2\cdot [5]^6/\,[1]^5[25].\label{eq:gatherend}
\end{gather}
To continue with the $N=9$ example, $Q_9(t_9)=t_9(t_9^2+9t_9+27),$ which
equals $3^6\cdot[3]^{12}/\,[1]^{12},$ yielding $h_9(t_9(\tau))=[1]^3/\,[3]$
as claimed.  The remaining three `irrational' values of~$N$ are handled
similarly.

To determine if the multiplier system of~$\mathfrak{h}_N$ is given by an
(odd) Dirichlet character mod~$N,$ reason as~follows.  Let
$g_1,\dots,g_l\in\bar\Gamma_0(N)$ be a generating set
for~$\bar\Gamma_0(N),$ and let $a(g_i),b(g_i),\allowbreak c(g_i),d(g_i)$
denote the matrix elements of~$g_i$.  For each $1\le i\le l,$ the character
$\hat\chi_N(a(g_i),b(g_i),c(g_i),d(g_i))$ can be computed from the eta
product for~$\mathfrak{h}_N$ given in Table~\ref{tab:xNhN}, using the
transformation formula for~$\eta,$ which is~\cite{Rademacher73}
\begin{equation*}
  \eta(\tfrac{a\tau+b}{c\tau+d}) =\left\{
  \begin{array}{ll}
    \left(\frac{d}c\right)\zeta_{24}^{3(1-c)+bd(1-c^2)+c(a+d)}[-{\rm i}(c\tau+d)]^{1/2}\,\eta(\tau), & \quad c{\rm\ odd},\\
    \left(\frac{c}d\right)\zeta_{24}^{3d+ac(1-d^2)+d(b-c)}[-{\rm i}(c\tau+d)]^{1/2}\,\eta(\tau), & \quad d{\rm\ odd},\\
  \end{array}
\right.
\end{equation*}
if $c>0,$ with $\bigl(\tfrac{c}{-d}\bigr):=\bigl(\tfrac{c}{d}\bigr)$ by
convention.  If there is a unique Dirichlet character mod~$N,$ say~$\chi,$
with $\chi(d(g_i)) = \hat\chi_N(a(g_i),b(g_i),c(g_i),d(g_i))$ for $1\le
i\le l,$ then one must have $\chi_N=\chi$.  If on the other hand no~such
Dirichlet character exists, then the multiplier system of~$\mathfrak{h}_N$
cannot be so~encoded.

This procedure requires, for each genus-zero $\Gamma_0(N),$ an explicit
generating set.  The cardinality of any {\em minimal\/} generating set is
$(\sigma_\infty+\epsilon_2+\epsilon_3)(N)-1,$ the number of fixed points
minus~$1,$ but the generating set need not be minimal.  For prime~$N,$
a~standard result is available.  Rademacher~\cite{Rademacher29} showed that
the set of matrices comprising
$\pm\bigl(\begin{smallmatrix}1&1\\0&1\end{smallmatrix}\bigr)$ and
$\pm\bigl(\begin{smallmatrix}k'&1\\-(kk'+1)&-k\end{smallmatrix}\bigr),$
$k=1,\dots,{N-1},$ where $k'=k'(k)$ is determined by the condition
$kk'\equiv-1\pmod N,$ generates $\Gamma_0(N)$; and also worked~out a
minimal generating set for each prime $N$ up~to~$31$ (reproduced by
Apostol~\cite{Apostol90}).  For general~$N,$ one can use the algorithm of
Kulkarni~\cite{Kulkarni91}, which is provided by the \textsc{Magma} system,
or that of Lascurain Orive~\cite{Lascurain2002}.

Actually, we shall need generating sets for~$\Gamma_0(N)$ only in the cases
$N=3$ and $N=4,6,8,9$.  For the latter small composite values, we shall use
the (minimal) sets given by Harnad and McKay~\cite{Harnad2000}, who
obtained them by {\em ad~hoc\/} methods.

For the $14$ values of~$N,$ whether or not the multiplier system
of~$\mathfrak{h}_N$ is of Dirichlet type is determined as follows.

\renewcommand\labelitemi{$\circ$}
\begin{itemize}
  \item $N=2$.  There is no odd Dirichlet character mod~$2$.
  \item $N=3$.  A generating set for~$\Gamma_0(3),$
    from~\cite{Rademacher29}, is
    $\left(\begin{smallmatrix}1&1\\0&1\end{smallmatrix}\right),$
    $\left(\begin{smallmatrix}1&1\\-3&-2\end{smallmatrix}\right),$ up~to
    sign.  These have $d\equiv1,1\pmod3$ respectively.  The power
    of~$\zeta_{24}$ appearing in the transformation formula
    for~$\mathfrak{h}_3=[1]^3/\,[3]$ is computed to be $0,20$ respectively.
    Since $0\neq20,$ the multiplier system is non-Dirichlet.
  \item $N=4$.  A generating set for~$\Gamma_0(4),$ from~\cite{Harnad2000},
    is $\left(\begin{smallmatrix}1&1\\0&1\end{smallmatrix}\right),$
    $\left(\begin{smallmatrix}1&-1\\4&-3\end{smallmatrix}\right),$ up~to
    sign.  These have $d\equiv1,1\pmod4$ respectively.  The power
    of~$\zeta_{24}$ appearing in the transformation formula for
    $\mathfrak{h}_4=\mathfrak{h}_2=[1]^4/\,[2]^2$ is computed to be $0,12$
    respectively.  Since $0\neq12,$ the multiplier system is non-Dirichlet.
  \item $N=6,12,18$.  A generating set for~$\Gamma_0(6),$
    from~\cite{Harnad2000}, is
    $\left(\begin{smallmatrix}1&1\\0&1\end{smallmatrix}\right),$
    $\left(\begin{smallmatrix}5&-3\\12&-7\end{smallmatrix}\right),$
    $\left(\begin{smallmatrix}5&-2\\18&-7\end{smallmatrix}\right),$ up~to
    sign.  These have $d\equiv1,5,5\pmod{6}$.  The power of~$\zeta_{24}$
    appearing in the transformation formula
    for~$\mathfrak{h}_6=[1]^6[6]\,/\,[2]^3[3]^2$ is computed to be
    $0,12,12$ respectively.  Hence $\chi_{6}:(\mathbb{Z}/6\mathbb{Z})^*\to
    U(1)$ is the odd Dirichlet character
    $d\mapsto\bigl(\frac{-3}{d}\bigr),$ which takes $1,5$ to $1,-1$.  Since
    $\Gamma_0(12),\allowbreak\Gamma_0(18)$ are subgroups of~$\Gamma_0(6),$
    and $\mathfrak{h}_{12},\mathfrak{h}_{18}$ equal~$\mathfrak{h}_6,$ the
    multiplier systems when $N=12,18$ must also be of Dirichlet type, with
    the same character.
  \item $N=8,16$.  A generating set for~$\Gamma_0(8),$
    from~\cite{Harnad2000}, is
    $\left(\begin{smallmatrix}1&1\\0&1\end{smallmatrix}\right),$
    $\left(\begin{smallmatrix}3&-2\\8&-5\end{smallmatrix}\right),$
    $\left(\begin{smallmatrix}3&-1\\16&-5\end{smallmatrix}\right),$ up~to
    sign.  These have $d\equiv1,3,3\pmod{9}$.  The power of $\zeta_{24}$
    appearing in the transformation formula
    for~$\mathfrak{h}_8=\mathfrak{h}_2=[1]^4/\,[2]^2$ is computed to be
    $0,12,12$ respectively.  Hence $\chi_{8}:(\mathbb{Z}/8\mathbb{Z})^*\to
    U(1)$ is the odd Dirichlet character $d\mapsto\bigl(\frac{-4}d\bigr),$
    which takes $1,3,\allowbreak5,7$ to $1,-1,\allowbreak1,-1$.  Since
    $\Gamma_0(16)$ is a subgroup of~$\Gamma_0(8),$ and $\mathfrak{h}_{16}$
    equals~$\mathfrak{h}_8,$ the multiplier system when $N=16$ must also be
    of Dirichlet type, with the same character.
  \item $N=9$.  A generating set for~$\Gamma_0(9),$ from~\cite{Harnad2000},
    is $\left(\begin{smallmatrix}1&1\\0&1\end{smallmatrix}\right),$
    $\left(\begin{smallmatrix}5&-4\\9&-7\end{smallmatrix}\right),$
    $\left(\begin{smallmatrix}2&-1\\9&-4\end{smallmatrix}\right),$ up~to
    sign.  These have $d\equiv1,2,5\pmod{9}$.  The power of~$\zeta_{24}$
    appearing in the transformation formula
    for~$\mathfrak{h}_9=\mathfrak{h}_3=[1]^3/\,[3]$ is computed to be
    $0,12,12$ respectively.  Hence $\chi_{9}:(\mathbb{Z}/9\mathbb{Z})^*\to
    U(1)$ is the odd Dirichlet character $d\mapsto\bigl(\frac{-3}d\bigr),$
    which takes $1,2,\allowbreak4,\allowbreak5,\allowbreak7,8$ to
    $1,-1,1,-1,1,-1$.
  \item $N=5,7,10,13,25$.  Careful numerical computation
    of~$\mathfrak{h}_N,$ using the infinite product representation for the
    eta function, reveals that in these cases, the character
    $\hat\chi_N(N-1,1,-N,-1)$ is respectively
    $\zeta_{24}^0,\zeta_{24}^{20},\allowbreak\zeta_{24}^0,\allowbreak\zeta_{24}^8,\zeta_{24}^0$.
    In each case this is inconsistent with $\chi_N(-1)=\zeta_{24}^{12}=-1,$
    so in none is the multiplier system of~$\mathfrak{h}_N$ given by a
    Dirichlet character mod~$N$. \qedhere
\end{itemize}
\end{proof}

\begin{table}
\caption{Eta product representations for $\hat t_N$ and $(q\,d/dq)\hat
  t_N,$ in the cases when $\Gamma_0(N)$ has no~elliptic fixed points.}
\begin{center}
{\small
\begin{tabular}{cll}
\hline\noalign{\smallskip}
$N$ & \hfil $\hat t_N$ \hfil \hfil \hfil \hfil & $(q\,d/dq)\,\hat t_N$ \\
\noalign{\smallskip}\hline\noalign{\smallskip}
$4$ & $[4]^{8}/\,[1]^{8}$ & $[2]^{20}/\,[1]^{16}$ \\
$6$ & $[2][6]^{5}/\,[1]^{5}[3]$ & $[2]^8[3]^6/\,[1]^{10}$ \\
$8$ & $[2]^{2}[8]^{4}/\,[1]^{4}[4]^{2}$ & $[2]^8[4]^4/\,[1]^{8}$\\
$9$ & $[9]^{3}/\,[1]^{3}$ & $[3]^{10}/\,[1]^6$\\
$12$ & $[2]^2[3]\,[12]^3/\,[1]^{3}[4]\,[6]^2$ & $[2]^6[3]^2[6]^2/\,[1]^6$\\
$16$ & $[2]\,[16]^{2}/\,[1]^{2}[8]$ & $[2]^2[4]^6/\,[1]^4$ \\
$18$ & $[2]\,[3]\,[18]^2/\,[1]^{2}[6]\,[9]$ & $[2]^2[3]^4[6]^2/\,[1]^4$ \\
\noalign{\smallskip}\hline
\end{tabular}
}%
\end{center}
\label{tab:shoehorned}
\end{table}

\begin{corollary}
  For each\/ $N$ for which\/ $\Gamma_0(N)$~is of genus zero and has
  no~elliptic fixed points, there is an eta product representation not~only
  for the canonical Hauptmodul\/ $t_N=\kappa_N\cdot\hat t_N,$ but also one
  for the weight\/-$2$ modular form\/ $[(2\pi{\rm i})^{-1}\,d/d\tau]t_N,$
  i.e., for\/ $(q\,d/dq)t_N,$ as given in
  Table\/~{\rm\ref{tab:shoehorned}.}
\end{corollary}
\begin{remark*}
  An eta product for the derivative $(q\,d/dq)\hat t_N,$ when $\hat t_N$~is
  similarly expressed, constitutes a {\em combinatorial identity\/}.
  Several of the identities in Table~\ref{tab:shoehorned} are new, though
  Fine~\cite[\S33]{Fine88} gives the eta products for $(q\,d/dq)\hat t_4,$
  $(q\,d/dq)\hat t_6$ (in~effect), and $(q\,d/dq)\hat t_8,$ and
  Cooper~\cite{Cooper2005} gives the one for $(q\,d/dq)\hat t_9$.  By using
  the chain rule, one can also work~out an eta product for the derivative
  of each non-canonical Hauptmodul in Table~\ref{tab:alt}.  For instance,
  consider the invariant $\alpha=t_4/(t_4+16),$ which equals
  $2^4\cdot[1]^8[4]^{16}/\,[2]^{24}.$ One easily deduces that
  $(q\,d/dq)\alpha = 2^4\cdot[1]^{16}[4]^{16}/\,[2]^{28}$.
\end{remark*}
\begin{proof}
  If $\Gamma_0(N)$ has no elliptic fixed points, the alternative
  representation for the modular form $\mathfrak{h}_N(\tau)=h_N(t_N(\tau))$
  given in Corollary~\ref{cor:hNformula2} implies
  \begin{equation}
    \frac1{2\pi{\rm i}}\, \frac{dt_N}{d\tau}(\tau) =
    t_N(\tau)\,h_N^2(t_N(\tau))\!\!\prod_{i=1}^{\sigma_\infty(N)-2} 
    \,
    \Biggl[1-\frac{t_N(\tau)}{t_{N,i}^{(\infty)}}\Biggr].
  \end{equation}
  The formulas in Table~\ref{tab:shoehorned} follow from this, with the aid
  of the eta product representations given in Tables \ref{tab:alt}
  and~\ref{tab:xNhN}.  When $N=9,16,$ or~$18,$ not~all the
  $\sigma_\infty(N)-2$ cusps $t_N=t_{N,i}^{(\infty)}$ are rational, and for
  the right-hand side one needs also special eta product formulas:
  (\ref{eq:gatherstart}), (\ref{eq:special16}), or
  (\ref{eq:special18a})\textendash\nobreak(\ref{eq:special18b}).
\end{proof}

One can derive alternative (non-canonical) modular
forms~$\tilde{\mathfrak{h}}_N$ from~$\mathfrak{h}_N,$ just as one derives
the alternative Hauptmoduln~$\tilde t_N$ of Table~\ref{tab:alt} from the
canonical Hauptmodul~$t_N$.  For each~$N,$ if $\mathfrak{c}_0\not={\rm
i}\infty$ is a cusp of~$\Gamma_0(N),$ there is a Hauptmodul~$\tilde
t_N:=t_N/(t_N-t^*),$ where $t^*$~is the value of~$t_N$ at~$\mathfrak{c}_0,$
as listed in Table~\ref{tab:hauptmoduln}.  It has divisor
$(t_N=0)-(t_N=t^*)$ on~$X_0(N)\cong\mathbb{P}^1(\mathbb{C})_{t_N}$.  If the
cusp is rational, then $\tilde t_N$~is a rational (over~$\mathbb{Q}$)
function of~$t_N$.  The associated alternative $\tilde{\mathfrak{h}}_N$
to~$\mathfrak{h}_N$ is naturally taken to be
$\left[(t_N-t^*)/t^*\right]^{\psi(N)/12}\,{\mathfrak{h}}_N$.  By
Corollary~\ref{cor:hNformula1}, this~$\tilde{\mathfrak{h}}_N,$ regarded as
a (multivalued) function on~$X_0(N),$ will have divisor
$\frac{\psi(N)}{12}\cdot (t_N=t^*)$.

\begin{table}
\caption{Weight-$1$ modular forms
  $\tilde{\mathfrak{h}}_N(\tau)=\tilde{h}_N(t_N(\tau))$ on~$\Gamma_0(N)$
  with integer-coefficient $q$-series, including the canonical one,
  ${\mathfrak{h}}_N(\tau)={h}_N(t_N(\tau))$.  Here, $\sum$ signifies
  $\sum_{n=1}^\infty$.}
\begin{center}
{\small
\begin{tabular}{cllp{1.4in}}
\hline\noalign{\smallskip}
$N$ & $\tilde h_N=\tilde h_N(t_N)$ &$\eta$-prod.\ for $\tilde h_N(t_N(\tau))$&$q$-series for $\tilde h_N(t_N(\tau))$\\
\noalign{\smallskip}\hline\noalign{\vskip1pt}\hline\noalign{\vskip1.5pt}
$2$ & $h_2$ & $[1]^4/\,[2]^2$ & $1 + 4\sum E_1(n;4) (-q)^n$ \\
\noalign{\vskip 1pt}\hline\noalign{\vskip 1pt}
$3$ & $h_3$ & $[1]^3/\,[3]$ & $1 - 3\sum \bigl[E_1(n;3)$ \\
    &       &               & $\hfill{}-3E_1(n/3;3)\bigl]\, q^n$ \\
\noalign{\vskip 1pt}\hline\noalign{\vskip 1pt}
$4$ & $h_4$ & \quad[$=h_2$] &  \\
$4$ & $[(t_4+16)/16]^{1/2}\,h_4$ & $[2]^{10}/\,[1]^4[4]^4$ & $1 + 4\sum E_1(n;4)\,q^n$ \\
\noalign{\vskip 1pt}\hline\noalign{\vskip 1pt}
$6$ & $h_6$ & $[1]^6[6]\,/\,[2]^3[3]^2$ & $1-6\sum\bigl[E_1(n;3)$ \\
    &       &                           & $\hfill{}-2E_1(n/2;3)\bigr]\,q^n$ \\
$6$ & $[(t_6+8)/8]\,h_6$ & $[2]^6[3]\,/\,[1]^3[6]^2$ & $1+3\sum E_1(n;6)\,q^n$ \\
$6$ & $[(t_6+9)/9]\,h_6$ & $[2][3]^6/\,[1]^2[6]^3$ & $1+2\sum E_{1,2}(n;6)\,q^n$ \\
\noalign{\vskip 1pt}\hline\noalign{\vskip 1pt}
$8$ & $h_8$ & \quad[$=h_2$] &  \\
$8$ & $[(t_8+4)/4]\,h_8$ & $[2]^{10}/\,[1]^4[4]^4$ & $1 + 4\sum E_1(n;4)\, q^n$ \\
$8$ & $[(t_8+8)/8]\,h_8$ & $[4]^{10}/\,[2]^4[8]^4$ & $1 + 4\sum E_1(n;4)\, q^{2n}$ \\
\noalign{\vskip 1pt}\hline\noalign{\vskip 1pt}
$9$ & $h_9$ & \quad[$=h_3$]  &  \\
\noalign{\vskip 1pt}\hline\noalign{\vskip 1pt}
$12$ & $h_{12}$ & \quad[$=h_6$]  &  \\
$12$ & $[(t_{12}+2)/2]^2\,h_{12}$ & $[2]^{15}[3]^2[12]^2$ & $1-6\sum\bigl[E_1(n;3)$ \\
     &                            & $\quad{}/\,[1]^6[4]^6[6]^5$  & $\hfill{}-2E_1(n/2;3)\bigr](-q)^n$ \\
$12$ & $[(t_{12}+3)/3]^2\,h_{12}$ & $[2][3]^6/\,[1]^2[6]^3$ & $1+2\sum E_{1,2}(n;6)\,q^n$ \\
$12$ & $[(t_{12}+4)/4]^2\,h_{12}$ & $[4]^6[6]\,/\,[2]^3[12]^2$ & $1+3\sum E_1(n;6)\,q^{2n}$ \\
$12$ & $[(t_{12}+6)/6]^2\,h_{12}$ & $[1]^2[4]^2[6]^{15}$ & $1+2\sum E_{1,2}(n;6)(-q)^n$ \\
   &                              & $\quad{}/\,[2]^5[3]^6[12]^6$ &  \\
\noalign{\vskip 1pt}\hline\noalign{\vskip 1pt}
$16$ & $h_{16}$ & \quad[$=h_2$] &  \\
$16$ & $[(t_{16}+2)/2]^2\,h_{16}$ & $[2]^{10}/\,[1]^4[4]^4$ & $1 + 4\sum E_1(n;4)\,q^n$ \\
$16$ & $[(t_{16}+4)/4]^2\,h_{16}$ & $[8]^{10}/\,[4]^4[16]^4$ & $1 + 4\sum E_1(n;4)\,q^{4n}$ \\
\noalign{\vskip 1pt}\hline\noalign{\vskip 1pt}
$18$ & $h_{18}$ & \quad[$=h_6$] &  \\
$18$ & $[(t_{18}+2)/2]^3\,h_{18}$ & $[2]^6[3]\,/\,[1]^3[6]^2$ & $1+3\sum E_1(n;6)\,q^n$ \\
$18$ & $[(t_{18}+3)/3]^3\,h_{18}$ & $[6][9]^6/\,[3]^2[18]^3$ & $1+2\sum E_{1,2}(n;6)\,q^{3n}$ \\
\noalign{\smallskip}\hline
\end{tabular}
}%
\end{center}
\label{tab:windup}
\end{table}

One example will suffice.  The $\alpha$-invariant, i.e., the alternative
Hauptmodul $\tilde t_4=t_4/(t_4+16)$ for~$\Gamma_0(4),$ equals
$2^4\cdot[1]^8[4]^{16}/\,[2]^{24}$ and corresponds via the $2$-isogeny
$\tau\mapsto\tau/2$ to the $\lambda$-invariant
$2^4\cdot[\frac12]^8[2]^{16}/\,[1]^{24},$ a~Hauptmodul for~$\Gamma(2)$.
The associated weight-$1$ modular form
$\tilde{\mathfrak{h}}_4=\left[(t_4+16)/16\right]^{1/2} \mathfrak{h}_4$
equals, by direct computation, the eta product $[2]^{10}/\,[1]^4[4]^4$.
Under the $2$-isogeny this corresponds to $[1]^{10}/\,[\frac12]^4[2]^4,$
a~weight-$1$ modular form for~$\Gamma(2)$.

If the eta product for~$\mathfrak{h}_N(\tau)$ contains fractional powers
of~$\eta,$ i.e., if $N=5,7,10,\allowbreak 13,25,$ its $q$-expansion about
the infinite cusp $q=0,$ say $\sum_{n=0}^\infty a_n^{(N)}q^n,$ turns~out to
contain coefficients $a_n^{(N)}$ that are not integers (but of~course
$a_n^{(N)}\in\mathbb{Q}$ for all~$n$).  The $q$-series is nonetheless
`almost integral': one can show that the associated scaled sequence
$\alpha_N^na_n^{(N)},$ $n\geq0,$ where
\begin{equation}
\label{eq:alphaN}
  \alpha_N:=2^2,3^2,2,2^23^2,2^2,\qquad{\rm for}\quad N=5,7,10,13,25,
\end{equation}
is an integer sequence.  The number-theoretic interpretation of the
integers $\alpha_N^na_n^{(N)}$ is unclear.  The integral $q$-series for the
nine remaining~$\mathfrak{h}_N,$ and for the
alternatives~$\tilde{\mathfrak{h}}_N$ derived from them, can more easily be
expressed in closed form.

\begin{theorem}
The canonical and alternative weight\/-$1$ modular forms\/
$\mathfrak{h}_N,\tilde{\mathfrak{h}}_N$ for\/~$\Gamma_0(N),$ when\/
$N=2,3,4,6,8,9,12,16,18,$ have the eta product and\/ $q$-series
representations given in Table\/~{\rm\ref{tab:windup}}{\rm.}  In the\/
$q$-series, $E_{r,s,\dots}(n;k)$ denotes the excess of the number of
divisors of\/~$n$ congruent to\/ $r,s,\dots\!\pmod k$ over the number
congruent to\/~$-r,-s,\dots\!\pmod k$; or zero, if\/ $n$ is not an integer.
\end{theorem}
\begin{proof}
The eta products for the $\tilde{\mathfrak{h}}_N$ follow from those for
$t_N,\mathfrak{h}_N$ given in Table~\ref{tab:xNhN}.  Several of these eta
products were expanded in multiplicative $q$-series by
Fine~\cite[\S32]{Fine88}, and $q$-expansions of the remainder follow by
applying such transformations as $q\mapsto q^2$ and~$q\mapsto-q$.  Under
$q\mapsto-q,$ i.e., $\tau\mapsto\tau+\frac12,$ the function $[m]$
on~$\mathcal{H}\ni\tau$ is taken to itself if the integer $m$~is even, and
to $[2m]^3/\,[m][4m]$ if it is odd.
\end{proof}

\begin{table}
\caption{Alternative weight-$1$ modular forms
  $\bar{\mathfrak{h}}_N(\tau)=\bar h_N(t_N(\tau))$ for~$\Gamma_0(N)$ that
  are zero at the infinite cusp, and have integer-coefficient $q$-series.
  Here, $\sum$ signifies $\sum_{n=1}^\infty$.}
\begin{center}
{\small
\begin{tabular}{clll}
\hline\noalign{\smallskip}
$N$ & $\bar h_N=\bar h_N(t_N)$ & $\eta$-prod.\ for $\bar h_N(t_N(\tau))$ & $q$-series  for $\bar h_N(t_N(\tau))$  \\
\noalign{\smallskip}\hline\noalign{\smallskip}
$2$ & $(t_2/2^{12})^{1/4}\,h_2$ & $[2]^4/\,[1]^2$ & $q^{1/4}\left[1 + \sum E_1(4n+1;4)\,q^n\right]$ \\
$3$ & $(t_3/3^{6})^{1/3}\,h_3$ & $[3]^3/\,[1]$ & $q^{1/3}\left[1 + \sum E_1(3n+1;3)\,q^n\right]$ \\
$4$ & $(t_4/2^{8})^{1/2}\,h_4$ &  \quad[$=\bar h_2(t_2(2\tau)$] &  \\
$6$ & $(t_6/2^33^2)\,h_6$ &  $[1][6]^6/\,[2]^2[3]^3$ & $q\bigl\{1+\sum\bigl[E_1(n;6)$ \\
    &  &   & $\qquad\qquad\;{}-2E_1(n/2;3)\bigr]\,q^n\bigr\}$ \\
$8$ & $(t_8/2^5)\,h_8$ &  \quad[$=\bar h_2(t_2(4\tau)$] &  \\
$9$ & $(t_9/3^3)\,h_9$ &  \quad[$=\bar h_3(t_3(3\tau)$]  &  \\
$12$ & $(t_{12}/2^23)^2\,h_{12} $ &  \quad[$=\bar h_6(t_6(2\tau))$]  &  \\
$16$ & $(t_{16}/6)^2\,h_{16}$ &  \quad[$=\bar h_2(t_2(8\tau))$]  &  \\
$18$ & $(t_{18}/2^33^3)^2\,h_{18}$ &  \quad[$=\bar h_6(t_6(3\tau))$]  &  \\
\noalign{\smallskip}\hline
\end{tabular}
}%
\end{center}
\label{tab:insert}
\end{table}

One can also derive from each~$\mathfrak{h}_N$ a weight-$1$ modular form
$\bar{\mathfrak{h}}_N$ for~$\Gamma_0(N)$ that has a zero at the infinite
cusp, rather than at one of the finite ones.  This can be accomplished by
applying the Fricke involution~$W_N,$ which on the
half plane~$\mathcal{H}$ is the map $\tau\mapsto-1/N\tau$.  Equivalently
(up~to a constant factor), one can let
$\bar{\mathfrak{h}}_N:=(t_N/\kappa_N)^{\psi(N)/12}\,\mathfrak{h}_N$.
Defined thus, $\bar{\mathfrak{h}}_N$~will have divisor
$\frac{\psi(N)}{12}\cdot(t_N=0)$ on~$X_0(N),$ and hence will equal zero
at~$\tau={\rm i}\infty$.  The nine modular forms
$\bar{\mathfrak{h}}_N(\tau)$ that have integral $q$-series are listed in
Table~\ref{tab:insert}.  There are only three fundamental ones, the others
being obtained by modular substitutions $\tau\mapsto \ell\tau,$ i.e.,
$q\mapsto q^\ell$.  The $q$-expansion of~$\bar{\mathfrak{h}}_6(\tau)$ is
due~to M.~Somos (unpublished); the others, to Fine.

The reader will recall that each of the canonical modular forms
$\mathfrak{h}_N(\tau)=h_N(t_N(\tau))$ was originally pulled back from $\hat
h_1(\hat J(\tau))=E_4(\tau)^{1/4},$ the fourth root of a weight-$4$ modular
form for~$\Gamma(1)$.  Due~to multivaluedness this is not a true modular
form, but it does have an integral $q$-expansion about $\tau={\rm
i}\infty$.  By a useful result of Heninger et~al.~\cite{Heninger2006}, a
$q$-series~$f$ in $1+q\mathbb{Z}[[q]]$ has the property that it
equals~$g^k,$ for some $g$~in $1+q\mathbb{Z}[[q]],$ iff the reduction
of~$f$ mod~${\mu_k}$ has the same property, where $\mu_k:=k\prod_{p|k}p$.
Since the Eisenstein sum $E_4$ equals $1-240\sum\sigma_3(n)q^n,$ it follows
that its fourth and eighth roots must be in $1+q\mathbb{Z}[[q]]$.  Hence,
not~only does $\hat h_1(\hat J(\tau))$ have an integral $q$-expansion, but
so does its square root.  The $q$-expansion of $\hat h_1(\hat J(\tau))$ is
$1-60[q+99q^2+14248q^3+\cdots]$.

\section{The Modular Form $\mathfrak{h}_N(\tau)$: Explicit Computations~(II)}
\label{sec:modular3}

We now return to treating each canonical weight-$1$ form
$\mathfrak{h}_N(\tau)$ as a function of the corresponding Hauptmodul~$t_N$
for~$\Gamma_0(N),$ i.e., as a (multivalued) function~$h_N$ on the genus-$0$
curve $X_0(N)=\Gamma_0(N)\setminus\mathcal{H}^*\cong
\mathbb{P}^1(\mathbb{C})_{t_N}$.  The Picard--Fuchs equations that the
functions~$h_N$ satisfy, and their $t_N$-expansions, will be computed.
These two pieces of information will place the~$h_N$ firmly in the world of
special functions.  We pay particular attention to closed-form expressions
and recurrences for the coefficients of their $t_N$-expansions, since a
functional equation for~$h_N$ can be viewed as a {\em series identity\/}.
It turns~out that many of the coefficient sequences have cropped~up in
other contexts, and are listed in Sloane's {\em Encyclopedia\/}.

\begin{table}
\caption{The (multivalued) functions $h_N=h_N(t_N)$ that define the
canonical modular forms $\mathfrak{h}_N(\tau)=h_N(t_N(\tau)),$ expressed
in~terms of the special function ${}_2F_1(\frac1{12},\frac5{12};1;\cdot)$.
Here, $t:=t_N$.}
\begin{center}
{\small
\begin{tabular}{cl}
\hline\noalign{\smallskip}
$N$ & $\quad h_N(t_N)$\\
\noalign{\smallskip}\hline\noalign{\smallskip}
$2$ & $[16^{-1}(t+16)]^{-1/4}\,{}_2F_1\left(\frac1{12},\frac5{12};\,1;\,\frac{12^3t}{(t+16)^3}\right)$ \\
$3$ & $[3^{-6}(t+27)(t+3)^3]^{-1/12}\,{}_2F_1\left(\frac1{12},\frac5{12};\,1;\,\frac{12^3t}{(t+27)(t+3)^3}\right)$ \\
$4$ & $[16^{-1}(t+16)\,\circ\,t(t+16)]^{-1/4}\,{}_2F_1\left(\frac1{12},\frac5{12};\,1;\,\frac{12^3t}{(t+16)^3} \circ\, t(t+16)\right)$ \\
$5$ & $[5^{-1}(t^2+10t+5)]^{-1/4}\,{}_2F_1\left(\frac1{12},\frac5{12};\,1;\,\frac{12^3t}{(t^2+10t+5)^3}\right)$ \\
$6$ & $[144^{-1}(t+6)(t^3+18t^2+84t+24)]^{-1/4}\,$ \\
    & $\qquad{}\times\,_2F_1\left(\frac1{12},\frac5{12};\,1;\,\frac{12^3t}{(t+16)^3} \circ \frac{t(t+8)^3}{t+9}=\frac{12^3t}{(t+27)(t+3)^3} \circ \frac{t(t+9)^2}{t+8}\right)$ \\
$7$ & $[49^{-1}{(t^2+13t+49)(t^2+5t+1)^3}]^{-1/12}$ \\
    & $\qquad{}\times\,{}_2F_1\left(\frac1{12},\frac5{12};\,1;\,\frac{12^3t}{(t^2+13t+49)(t^2+5t+1)^3}\right)$ \\
$8$ & $[16^{-1}(t+16)\,\circ\,t(t+16)\,\circ\,t(t+8)]^{-1/4}$ \\
    & $\qquad{}\times\,{}_2F_1\left(\frac1{12},\frac5{12};\,1;\,\frac{12^3t}{(t+16)^3}\circ\,t(t+16)\,\circ\,t(t+8)\right)$ \\
$9$ & $[3^{-6}(t+27)(t+3)^3\circ\,t(t^2+9t+27)]^{-1/12}$ \\
    & $\qquad{}\times\,{}_2F_1\left(\frac1{12},\frac5{12};\,1;\,\frac{12^3t}{(t+27)(t+3)^3}\circ\, t(t^2+9t+27)\right)$ \\
$10$ & $[80^{-1}(t^6+20t^5+160t^4+640t^3+1280t^2+1040t+80)]^{-1/4}\,$ \\
     & $\qquad{}\times\,_2F_1\left(\frac1{12},\frac5{12};\,1;\,\frac{12^3t}{(t+16)^3}\circ\frac{t(t+4)^5}{t+5}=\frac{12^3t}{(t^2+10t+5)^3}\circ\frac{t(t+5)^2}{t+4}\right)$ \\
$12$ & $[144^{-1}(t+6)(t^3+18t^2+84t+24)\,\circ\,t(t+6)]^{-1/4}\,$ times any of: \\
     & $\qquad{}_2F_1\left(\frac1{12},\frac5{12};\,1;\,\frac{12^3t}{(t+16)^3}\circ\frac{t(t+8)^3}{t+9}\circ\,t(t+6)\right),$ \\
     & $\qquad{}_2F_1\left(\frac1{12},\frac5{12};\,1;\,\frac{12^3t}{(t+27)(t+3)^3}\circ\frac{t(t+9)^2}{t+8}\circ\, t(t+6)\right),$ \\
     & $\qquad{}_2F_1\left(\frac1{12},\frac5{12};\,1;\,\frac{12^3t}{(t+16)^3}\circ\,t(t+16)\, \circ\frac{t(t+4)^3}{t+3}\right)$ \\
$13$ & $[13^{-1}(t^2+5t+13)(t^4+7t^3+20t^2+19t+1)^3]^{-1/12}$ \\
     & $\qquad{}\times\,{}_2F_1\left(\frac1{12},\frac5{12};\,1;\,\frac{12^3t}{(t^2+5t+13)(t^4+7t^3+20t^2+19t+1)^3}\right)$ \\
$16$ & $[16^{-1}(t+16)\,\circ\,t(t+16)\,\circ\,t(t+8)\,\circ\,t(t+4)]^{-1/4}$ \\
     & $\qquad{}\times\,{}_2F_1\left(\frac1{12},\frac5{12};\,1;\,\frac{12^3t}{(t+16)^3}\circ\, t(t+16)\,\circ\,t(t+8)\,\circ\,t(t+4)\right)$ \\
$18$ & $[144^{-1}(t+6)(t^3+18t^2+84t+24)\,\circ\,t(t^2+6t+12)]^{-1/4}\,$ times any of: \\
     & $\qquad{}_2F_1\left(\frac1{12},\frac5{12};\,1;\,\frac{12^3t}{(t+16)^3}\circ\frac{t(t+8)^3}{t+9}\circ\,t(t^2+6t+12)\right),$ \\
     & $\qquad{}_2F_1\left(\frac1{12},\frac5{12};\,1;\,\frac{12^3t}{(t+27)(t+3)^3}\circ\frac{t(t+9)^2}{t+8}\circ\,t(t^2+6t+12)\right),$ \\
     & $\qquad{}_2F_1\left(\frac1{12},\frac5{12};\,1;\,\frac{12^3t}{(t+27)(t+3)^3}\circ\,t(t^2+9t+27)\,\circ\frac{t(t+3)^2}{t+2}\right)$ \\
$25$ & $[5^{-1}(t^2+10t+5)\,\circ\,t(t^4+5t^3+15t^2+25t+25)]^{-1/4}$ \\
     & $\qquad{}\times\,{}_2F_1\left(\frac1{12},\frac5{12};\,1;\,\frac{12^3t}{(t^2+10t+5)^3}\circ\, t(t^4+5t^3+15t^2+25t+25)\right)$ \\
\noalign{\smallskip}\hline
\end{tabular}
}%
\end{center}
\label{tab:hNasfunctionofxN}
\end{table}

\begin{theorem}
  In a neighborhood of the point\/ $t_N=0$ {\rm(}i.e., the infinite cusp\/
  $\tau={\rm i}\infty${\rm),} each function\/ $h_N=h_N(t_N)$ can be
  expressed in~terms of the Gauss\ hypergeometric function
  ${}_2F_1(\frac1{12},\frac{5}{12};1;\cdot),$ as given in
  Table\/~{\rm\ref{tab:hNasfunctionofxN}.}
\end{theorem}
\begin{proof}
  This follows from Definition~\ref{def:2}, the many composite rational
  expressions $P_N(t_N)/Q_N(t_N)$ for the $j$-invariant being taken from
  Table~\ref{tab:jfactored}.
\end{proof}
\begin{remark*}
  The hypergeometric recurrence~(\ref{eq:hyperrecurrence}) yields a $\hat
  J$-expansion of the underlying function $\hat h_1=\hat h_1(\hat
  J)={}_2F_1(\frac1{12},\frac5{12};1;\hat J)$.  If $\hat h_1(\hat
  J)=\sum_{n=0}^\infty \hat c_n^{(1)}\hat J^n$ then the sequence $\hat
  d_n^{(1)} = 12^{3n}\,\hat c_n^{(1)},$ $n\geq0,$ is integral, the first
  few terms being $1,60,\allowbreak39780,\allowbreak38454000,43751038500$.
  This is Sloane's sequence \texttt{A092870}.
\end{remark*}

\begin{table}
\caption{Differential operators~$\mathcal{L}_N$ in the normal-form
Picard--Fuchs equations $\mathcal{L}_N\,u=0$ satisfied by the modular
forms~$\mathfrak{h}_N$ for~$\Gamma_0(N),$ when viewed as functions of
$t:=t_N$.  In the four sections of the table, $\mathcal{L}_N$~has
respectively $3,4,6,8$ singular points on the curve
$X_0(N)\protect\cong\mathbb{P}^1(\mathbb{C})_{t_N}$.}
\begin{center}
{\small
\begin{tabular}{cl}
\hline\noalign{\smallskip}
$N$ & Differential operator~$\mathcal{L}_N$\\
\noalign{\smallskip}\hline\noalign{\vskip1pt}\hline\noalign{\vskip1.5pt}
$2$ & $D_t^2+\left[\frac{1}{t}+\frac{1}{2(t+64)}\right]D_t+\frac{1}{16t(t+64)}$ \\
$3$ & $D_t^2+\left[\frac{1}{t}+\frac{2}{3(t+27)}\right]D_t+\frac{1}{9t(t+27)}$ \\
$4$ & $D_t^2+\left[\frac{1}{t}+\frac{1}{t+16}\right]D_t+\frac{1}{4t(t+16)}$ \\
\noalign{\vskip 1pt}\hline\noalign{\vskip 1pt}
$5$ & $D_t^2+\left[\frac{1}{t}+\frac{t+11}{t^2+22t+125}\right]D_t+\frac{t+10}{4t(t^2+22t+125)}$ \\
$6$ & $D_t^2+\left[\frac{1}{t}+\frac{1}{t+8}+\frac{1}{t+9}\right]D_t+\frac{t+6}{t(t+8)(t+9)}$ \\
$7$ & $D_t^2+\left[\frac{1}{t}+\frac{4t+26}{3(t^2+13t+49)}\right]D_t+\frac{4t+21}{9t(t^2+13t+49)}$ \\
$8$ & $D_t^2+\left[\frac{1}{t}+\frac{1}{t+4}+\frac{1}{t+8}\right]D_t+\frac{1}{t(t+8)}$ \\
$9$ & $D_t^2+\left[\frac{1}{t}+\frac{2t+9}{t^2+9t+27}\right]D_t+\frac{t+3}{t(t^2+9t+27)}$ \\
\noalign{\vskip 1pt}\hline\noalign{\vskip 1pt}
$10$ & $D_t^2+\left[\frac{1}{t}+\frac{1}{t+4}+\frac{1}{t+5}+\frac{t+4}{t^2+8t+20}\right]D_t+\frac{9t^3+95t^2+340t+400}{4t(t+4)(t+5)(t^2+8t+20)}$ \\
$12$ & $D_t^2+\left[\frac{1}{t}+\frac{1}{t+2}+\frac{1}{t+3}+\frac{1}{t+4}+\frac{1}{t+6}\right]D_t+\frac{4(t^2+6t+6)}{t(t+2)(t+4)(t+6)}$ \\
$13$ & $D_t^2+\left[\frac{1}{t}+\frac{t+3}{t^2+6t+13}+\frac{4t+10}{3(t^2+5t+13)}\right]D_t+\frac{49t^3+351t^2+1027t+1014}{36t(t^2+5t+13)(t^2+6t+13)}$ \\
$16$ & $D_t^2+\left[\frac{1}{t}+\frac{1}{t+2}+\frac{1}{t+4}+\frac{2t+4}{t^2+4t+8}\right]D_t+\frac{4(t+2)^2}{t(t+4)(t^2+4t+8)}$ \\
\noalign{\vskip 1pt}\hline\noalign{\vskip 1pt}
$18$ & $D_t^2+\left[\frac{1}{t}+\frac{1}{t+2}+\frac{1}{t+3}+\frac{2t+3}{t^2+3t+3}+\frac{2t+6}{t^2+6t+12}\right]D_t$ \\
   & $\quad\qquad\qquad\qquad\qquad\qquad\qquad\qquad\qquad{}+\frac{9(t+2)(t^3+6t^2+12t+6)}{t(t+3)(t^2+3t+3)(t^2+6t+12)}$ \\
$25$ & $D_t^2+\left[\frac{1}{t}+\frac{4t^3+15t^2+30t+25}{t^4+5t^3+15t^2+25t+25}+\frac{t+1}{t^2+2t+5}\right]D_t$ \\
 & $\qquad\qquad\qquad\qquad\qquad\qquad\qquad{}+\frac{25(t^5+5t^4+15t^3+25t^2+25t+10)}{4t(t^2+2t+5)(t^4+5t^3+15t^2+25t+25)}$ \\
\noalign{\vskip2pt}\hline
\end{tabular}
}%
\end{center}
\label{tab:fuchsianops}
\end{table}

For most purposes the explicit formulas of Table~\ref{tab:hNasfunctionofxN}
are far less useful than the Fuchsian differential equations that the
functions~$h_N$ satisfy, and representations that can be deduced from them.
Recall that by Theorem~\ref{thm:3}, the function~$h_N$ is the unique
holomorphic solution (up~to normalization) of a normal-form Picard--Fuchs
equation $\mathcal{L}_N u=0,$ in a neighborhood of the point $t_N=0$.  The
full space of local solutions is
$h_N(\cdot)\left[\mathbb{C}\tau(\cdot)+\mathbb{C}\right],$ where
$\tau=\tau(t_N)$~means any branch of what is, in~reality, an
infinite-valued function on~$X_0(N)\cong\mathbb{P}^1(\mathbb{C})_{t_N}$.

\begin{theorem}
  For each\/~$N,$ the normal-form Fuchsian differential operator\/
  $\mathcal{L}_N=D_{t_N}^2+\mathcal{A_N}(t_N) D_{t_N}+\mathcal{B}(t_N)$ in
  the Picard--Fuchs equation\/ $\mathcal{L}_Nu=0$ of
  Theorem\/~{\rm\ref{thm:3}}, which has holomorphic local solution\/
  $u=h_N,$ is as given in Table\/~{\rm\ref{tab:fuchsianops}}{\rm.}
\end{theorem}
\begin{proof}
  Each operator~$\mathcal{L}_N$ is computed as in Theorem~\ref{thm:4}.
  This is the procedure: (i)~Pull back the Gauss\ hypergeometric operator
  $\hat{\mathcal L}_1=L_{\frac1{12},\frac{5}{12};1}$ along $X_0(N)\to
  X(1),$ or equivalently, along $\mathbb{P}^1(\mathbb{C})_{t_N}\to
  \mathbb{P}^1(\mathbb{C})_{\hat J},$ using the formula for the
  degree-$\psi(N)$ covering map $j=P_N(t_N)/Q_N(t_N)$ given in
  Table~\ref{tab:coverings}, and (ii)~Make the substitution $\hat
  u=P_N(t_N)^{-1/12} u,$ to reduce the resulting operator to normal form.
  The computations are lengthy but rewarding.
\end{proof}
\begin{remark*}
By comparing each~$\mathcal{L}_N$ of Table~\ref{tab:fuchsianops} with the
normal form~(\ref{eq:Poole}), one sees that its singular points and
exponents are in agreement with Theorem~\ref{thm:3}.  In
each~$\mathcal{L}_N$ the coefficient $\mathcal{A_N}$ of~$D_{t_N}$ has only
simple poles, and they are located at the $t_N$-values of the fixed points
of~$\Gamma_0(N)$ on~$X_0(N)$ other than the cusp $t_N=\infty$ (i.e.,
$\tau=0$).  The locations agree with those given in
Table~\ref{tab:hauptmoduln}.  For example, the three (of~a total of four)
singular points on~$X_0(5)$ evident in the rational coefficient
$\frac1{t_5} + \frac{t_5+11}{t_5^2+22t_5+125}$ of~$D_{t_5}$
in~$\mathcal{L}_5$ are located at (i)~the infinite cusp $t_5=0,$ and
(ii)~the pair of quadratic elliptic points $t_5=-11\pm2{\rm i}$.  (As was
noted in~\S\ref{subsec:parametrized1b}, the latter are naturally bijective
with the representations $1^2+2^2,\allowbreak2^2+1^2$ of $N=5$ as a sum of
two squares.)  For each~$N,$ the coefficient $\mathcal{A}_N$ of~$D_{t_N}$
has residue at each cusp equal to~$1,$ at each quadratic elliptic point
to~$1/2,$ and at each cubic one to~$2/3$.
\end{remark*}
\begin{remark*}
  The equation $\mathcal{L}_9u=0$ is equivalent to the Picard--Fuchs
  equation derived (in~matrix form) by Dwork~\cite[\S8]{Dwork64} as the
  equation satisfied by the periods of the Hesse--Dixon elliptic
  family~(\ref{eq:hessedixonmodel}), parametrized by $\gamma=\gamma(\tau)=
  t_9(\tau/3)$.
\end{remark*}
\begin{remark*}
  For each of the 14~values of~$N,$ a Picard--Fuchs equation in the
  self-adjoint form $\left[D_{t_N}^2+\mathcal{Q}_N(t_N)\right]v=0,$
  projectively equivalent to the normal-form equation $\mathcal{L}_N u=0$
  of Table~\ref{tab:fuchsianops}, has been derived by Lian and
  Wiczer~\cite{Lian2006}, as part of a larger computational project on the
  $175$ genus-zero Conway--Norton subgroups of ${\it PSL}(2,\mathbb{R})$.
  They derive each $\mathcal{Q}_N\in\mathbb{Q}(t_N)$ in a heuristic way
  from the formula $\mathcal{Q}_N:=\frac12\{\tau,t_N\},$ starting with a
  $q$-expansion for each~$t_N$.  (For the first two dozen terms in each
  expansion, see~\cite[Table~3]{Ford94}.)  Each of their Hauptmoduln is of
  the McKay--Thompson type, with $q$-expansion $q^{-1}+0\cdot q^0 +O(q^1)$;
  so it is not identical to our Hauptmodul~$t_N,$ but rather related to it
  by a linear fractional transformation over~$\mathbb{Q}$.  Their
  Picard--Fuchs equations are indexed by Conway--Norton classes; for the
  correspondence to~$N,$ see, e.g.,~\cite[Table~4]{Ford94}.
\end{remark*}
\begin{remark*}
  From the Picard--Fuchs equations $\mathcal{L}_Nh_N=0,$ one can deduce
  that
\begin{subequations}
\begin{align}
\label{eq:agma}
h_2(t_2)&=1/\AG_2\bigl(1,\sqrt{1+\sqrt{1+t_2/64}}\bigr)\\
\label{eq:agmb}
h_4(t_4)&=1/\AG_2\bigl(1,\sqrt{1+t_4/16}\bigr)\\[3pt]
\label{eq:agmc}
h_8(t_8)&=1/\AG_2(1,1+t_8/4),\\[2pt]
\label{eq:agmd}
h_{16}(t_{16})&=1/\AG_2\bigl(1,(1+t_{16}/2)^2\bigr),  
\end{align}
\end{subequations}
where $\AG_2(\cdot,\cdot)$ is the quadratic arithmetic--geometric mean
(AGM) function of Gauss~\cite{Borwein87,Borwein91}, and that
\begin{subequations}
\begin{align}
\label{eq:agm3a}
h_3(t_3)&=1/\AG_3\bigl(1,\sqrt[3]{1+t_3/27}\bigr),\\[\jot]
\label{eq:agm3b}
h_9(t_9)&=1/\AG_3(1,1+t_9/3),
\end{align}
\end{subequations}
where $\AG_3(\cdot,\cdot)$ is the cubic AGM function of the
Borweins~\cite{Borwein91}.
\end{remark*}

\begin{table}
\caption{Recurrences satisfied by the coefficients $\{c_n^{(N)}\}$ of the
  series expansion $h_N(t_N)=\sum_{n=0}^\infty c_n^{(N)}t_N^n,$ in the
  cases when\/ $\Gamma_0(N)$ has $3,4,6,8$ fixed points on~$X_0(N)$.}
\begin{center}
{\small
\begin{tabular}{cp{3.85in}}
\hline\noalign{\smallskip}
$N$ & \hfil recurrence \\
\noalign{\smallskip}\hline\noalign{\vskip1pt}\hline\noalign{\vskip1.5pt}
$2$ &  $(4n-3)^2\,c_{n-1} + 1024n^2\,c_n=0$\\
$3$ &  $(3n-2)^2\,c_{n-1} + 243n^2\,c_n=0$\\
$4$ &  $(2n-1)^2\,c_{n-1} + 64n^2\,c_n=0$\\
\noalign{\vskip 1pt}\hline\noalign{\vskip 1pt}
$5$ &  $(2n-1)^2\,c_{n-1} + 2(44n^2+22n+5)\,c_n + 500(n+1)^2\,c_{n+1}=0$\\
$6$ &  $n^2\,c_{n-1} + (17n^2+17n+6)\,c_n + 72(n+1)^2\,c_{n+1}=0$\\
$7$ &  $(3n-1)^2\,c_{n-1} + 3(39n^2+26n+7)\,c_n + 441(n+1)^2\,c_{n+1}=0$\\
$8$ &  $n^2\,c_{n-1} + 4(3n^2+3n+1)\,c_n + 32(n+1)^2\,c_{n+1}=0$\\
$9$ &  $n^2\,c_{n-1} + 3(3n^2+3n+1)\,c_n + 27(n+1)^2\,c_{n+1}=0$\\
\noalign{\vskip 1pt}\hline\noalign{\vskip 1pt}
$10$ & $(2n+1)^2\,c_{n-1} + (68n^2+152n+95)\,c_n + 4(112n^2+392n+365)\,c_{n+1}$\\
     & \hfill ${}+ 80(17n^2+81n+99)\,c_{n+2} + 1600(n+3)^2\,c_{n+3}=0$\\
$12$ & $(n+1)^2\,c_{n-1} + 3(5n^2+15n+12)\,c_n + 16(5n^2+20n+21)\,c_{n+1}$\\
     & \hfill ${}+ 36(5n^2+25n+32)\,c_{n+2} + 144(n+3)^2\,c_{n+3}=0$\\
$13$ & $(6n+1)^2\,c_{n-1} + 3(132n^2+232n+117)\,c_n + (2016n^2+6384n+5395)\,c_{n+1}$\\
     & \hfill ${}+ 78(66n^2+302n+353)\,c_{n+2} + 6084(n+3)^2\,c_{n+2}=0$\\
$16$ & $(n+1)^2\,c_{n-1} + 2(5n^2+15n+12)\,c_n + 8(5n^2+20n+21)\,c_{n+1}$\\
     & \hfill ${}+ 16(5n^2+25n+32)\,c_{n+2} + 64(n+3)^2\,c_{n+3}=0$\\
\noalign{\vskip 1pt}\hline\noalign{\vskip 1pt}
$18$ & $(n+2)^2\,c_{n-1} + 2(7n^2+35n+45)\,c_n + 12(7n^2+42n+65)\,c_{n+1}$\\
     & \qquad\qquad\quad ${}+ 39(7n^2+49n+88)\,c_{n+2} + 72(7n^2+56n+114)\,c_{n+3}$\\
     & \hfill ${}+ 72(7n^2+63n+143)\,c_{n+4} + 216(n+5)^2\,c_{n+5}=0$\\
$25$ & $(2n+3)^2\,c_{n-1} + (28n^2+116n+125)\,c_n + 5(24n^2+128n+179)\,c_{n+1}$\\
     & \qquad\qquad\quad ${}+ 5(64n^2+416n+701)\,c_{n+2} + 25(24n^2+184n+361)\,c_{n+3}$\\
     & \hfill ${}+ 50(14n^2+124n+277)\,c_{n+4} + 500(n+5)^2\,c_{n+5}=0$\\
\noalign{\vskip2pt}\hline
\end{tabular}
}%
\end{center}
\label{tab:recurrences}
\end{table}

\begin{theorem}
  For each\/~$N,$ the sequence\/ $\{c_n^{(N)}\}_{n=0}^\infty$ of
  coefficients of the\/ $t_N$-expansion of the modular
  form\/~$\mathfrak{h}_N$ about the infinite cusp\/ $t_N=0$ satisfies the
  recurrence in Table\/~{\rm\ref{tab:recurrences},} initialized by\/
  $c^{(N)}_0=1$ and $c^{(N)}_n=0,$ $n<0${\rm.}  The number of terms is the
  number of fixed points of\/~$\Gamma_0(N),$ minus\/~$1${\rm.}  The scaled
  sequence\/ $d^{(N)}_n:=(\alpha_N\kappa_N)^n c^{(N)}_n,$ $n\geq0,$ is an\/
  {\em integral\/} sequence.  Here\/ $\kappa_N$~is as in
  Table\/~{\rm\ref{tab:xNhN}}, and\/ $\alpha_N$ equals\/ $2^2,3^2,2,2^23^2$
  if $N=5,7,10,13,$ and\/ $1$ otherwise.
\end{theorem}
\begin{remark*}
  That a factor~$\kappa_N^n$ is needed for integrality is due~to the
  definition of each Hauptmodul $t_N$ as $\kappa_N\cdot \hat t_N,$ where
  $\hat t_N$~is an eta product.  The factor $\alpha_N^n$ appeared
  in~\S\ref{sec:modular2}, as the superimposed geometric growth that forces
  the $q$-series of $\mathfrak{h}_N(\tau)$ to become integral.  (See
  Eq.~(\ref{eq:alphaN}).)
\end{remark*}
\begin{proof}
  Substitute $u=h_N=\sum_{n=0}^\infty c_n^{(N)}t_N^n$ into the equation
  $\mathcal{L}_Nu=0,$ and extract and analyse the recurrence equation
  satisfied by the coefficients.
\end{proof}

The integral sequences $d_n^{(N)},$ $n\geq0,$ for $N=2,3,$ do not currently
appear in Sloane's {\em Encyclopedia\/}~\cite{Sloane2005}, perhaps because
they are not expressible in binomial form.  Their first few terms are
$1,-4,\allowbreak100,\allowbreak-3600,\allowbreak152100,\allowbreak-7033104,\allowbreak344622096$;
and
$1,-3,\allowbreak36,\allowbreak-588,\allowbreak11025,\allowbreak-223587,\allowbreak4769856,$
respectively.  The integral sequences $d_n^{(N)},$ $n\geq0,$ for
$N=4,6,8,9,$ appear in the {\em Encyclopedia\/} as \texttt{A002894},
\texttt{A093388}, \texttt{A081085}, \texttt{A006077}.  One easily sees that
$d_n^{(4)}=2^{8n}c_n^{(4)}=(-)^n{\binom{2n}{n}}^2$.  For $N=6,8,9,$ the
recurrence satisfied by $d_n^{(N)}$ was derived by Coster~\cite{Coster83}
in an investigation of modular curves $\Gamma\setminus\mathcal{H}^*,$ where
$\Gamma<{\it PSL}(2,\mathbb{R})$ has four cusps, no~elliptic points, and is
of genus~$0$.  (The possible~$\Gamma$ include $\Gamma_0(6),$ $\Gamma_0(8),$
$\Gamma_0(9),$ and three others~\cite{Beauville82}.)
It was later shown that
\begin{alignat}{2}
d^{(6)}_n&=72^nc^{(6)}_n&&=\sum_{k=0}^n \binom{n}{k}(-8)^k\sum_{j=0}^{n-k}{\binom{n-k}{j}}^3, \\
d_n^{(8)}&=2^{5n}c_n^{(8)}&&=(-)^n\sum_{k=0}^n\binom{n}{k}\binom{2n-2k}{n-k}\binom{2k}{k},\\
d_n^{(9)}&=3^{3n}c_n^{(9)}&&=\sum_{k=0}^{\lfloor n/3\rfloor} (-)^k 3^{n-3k}\binom{n}{3k}\binom{3k}{k}\binom{2k}{k},
\end{alignat}
respectively by Verrill~\cite{Verrill99}; by Larcombe and
French~\cite[Thm.~3]{Larcombe2004} and V.~Jovovi\'c; and by Zagier and
Verrill.

\begin{table}
\caption{Alternative expressions for the modular forms~$\mathfrak{h}_N$ as
functions of $t:=t_N,$ in~terms of~${}_2F_1$ (when~$\Gamma_0(N)$~has three
fixed points) or the special function~$\Hl$ (when it has four).}
\begin{center}
{\small
\begin{tabular}{cl}
\hline\noalign{\smallskip}
$N$ & \hfil $h_N(t)$ \\
\noalign{\smallskip}\hline\noalign{\vskip1pt}\hline\noalign{\vskip1.5pt}
$2$ & ${}_2F_1(\frac14,\frac14;\,1;\,-t/64)$ \\
$3$ & ${}_2F_1(\frac13,\frac13;\,1;\,-t/27)$ \\
$4$ & ${}_2F_1(\frac12,\frac12;\,1;\,-t/16)$ \\
\noalign{\vskip 1pt}\hline\noalign{\vskip 1pt}
$5$ & $\Hl \left(\frac{-11\mp2{\rm i}}{-11\pm2{\rm i}},-\frac1{50}(-11\mp2{\rm i});\,\frac12,\frac12,1,\frac12;\,t/[-11\pm2{\rm i}\,]\right)\vphantom{\Bigl(\Bigr)}$ \\
$6$ & $\Hl \left(\frac98,\frac34;\,1,1,1,1;\,-t/8\right)=\Hl \left(\frac89,\frac23;\,1,1,1,1;\,-t/9\right)\vphantom{\Bigl(\Bigr)}$ \\
$7$ & $\Hl \left(\frac{-13\mp3\sqrt{-3}}{-13\pm3\sqrt{-3}},-\frac1{21}\bigl(\frac{-13\mp3\sqrt{-3}}{2}\bigr);\,\frac23,\frac23,1,\frac23;\,t/\bigl[\frac{-13\pm3\sqrt{-3}}{2}\bigr]\right)\vphantom{\Bigl(\Bigr)}$ \\
$8$ & $\Hl \left(2,1;\,1,1,1,1;\,-t/4\right)=\Hl \left(\frac12,\frac12;\,1,1,1,1;\,-t/8\right)\vphantom{\Bigl(\Bigr)}$ \\
$9$ & $\Hl \left(\frac{-9\mp3\sqrt{-3}}{-9\pm3\sqrt{-3}},-\frac1{9}\bigl(\frac{-9\mp3\sqrt{-3}}{2}\bigr);\,1,1,1,1;\,t/\bigl[\frac{-9\pm3\sqrt{-3}}{2}\bigr]\right)\vphantom{\Bigl(\Bigr)}$ \\
\noalign{\vskip2pt}\hline
\end{tabular}
}%
\end{center}
\label{tab:specfuncs}
\end{table}

\begin{theorem}
In the cases when\/ $\Gamma_0(N)$ has exactly\/ $3$ fixed points
on\/~$X_0(N)$ {\rm(}i.e., $N=2,3,4${\rm)}, and\/ $4$ fixed points\/
{\rm(}i.e., $N=6,7,8,9${\rm)}, the modular form\/ $\mathfrak{h}_N$ as a
function of the Hauptmodul\/ $t_N$ can be expressed in terms of the special
functions ${}_2F_1$ and\/~$\Hl,$ as shown in
Table\/~{\rm\ref{tab:specfuncs}}{\rm.}
\label{thm:nameless}
\end{theorem}
\begin{remark*}
  When $N=2,3,4,$ one can write, compactly and elegantly,
  \begin{equation}
    h_N(t_N) = {}_2F_1\bigl(\tfrac1{12}\psi(N),\tfrac1{12}\psi(N);\,1;\,-t_N/\kappa_N^{1/2}\bigr),
  \end{equation}
  since $\kappa_N=2^{12},3^{6},2^{8},$ respectively.
\end{remark*}
\begin{remark*}
Using the equivalence of the expressions of Table~\ref{tab:specfuncs} to
those of Table~\ref{tab:hNasfunctionofxN}, in the cases $N=2,3,4$ one
obtains transformation formulas for~${}_2F_1$.  They are respectively
cubic, quartic, and sextic, and are special cases of Goursat's
transformations.  In the cases $N=5,6,7,8,9$ one obtains reductions of
${\it Hl}$ to~${}_2F_1,$ which are new (such reductions have never been
classified).
\end{remark*}
\begin{proof}
Compare the recurrences of Table~\ref{tab:recurrences} with those for
${}_2F_1,\Hl$ given in the Appendix; or more directly, compare the
differential operators of Table~\ref{tab:fuchsianops} with the differential
equations satisfied by ${}_2F_1,\Hl,$ also given in the Appendix.  The
awkwardness of the expressions for $h_5,h_7,h_9$ (which are not
over~$\mathbb{Q}$) comes from the three singular points of each
of~$\mathcal{L}_N,$ $N=5,7,9,$ other than the cusp $t_N=\infty,$ not being
collinear on~$\mathbb{C}\ni t_N$. (Cf.\ Table~\ref{tab:hauptmoduln}.)
\end{proof}

\begin{corollary}
\label{cor:nameless}
When\/ $N=2,3,4,$ on the half-line\/ $t_N>0$ one has the identity
  \begin{equation*}
    \tau(t_N)/\,{\rm i} = A_N\,
    t_N^{-\psi(N)/12}\,\,\frac{h_N(\kappa_N/t_N)}{h_N(t_N)},
  \end{equation*}
with\/ $A_N=2,\sqrt3,2,$ respectively.  Here\/ $h_N,\tau$ are the
single-valued holomorphic continuations to\/~$t_N>0$ discussed
in\/~\S{\rm\ref{sec:modular},} which are real and positive.
\end{corollary}
\begin{proof}
  This continues the proof of Theorem~\ref{thm:clever}.  It is
  known~\cite[\S2.10]{Erdelyi53} that
  \begin{equation}
    \label{eq:nameless2}
    {}_2F_1(a,a;\,1;\,z)\sim\frac{\sin (\pi a)}{\pi}\,(-z)^{-1}  \left[\log(-z) + h_0(a)\right]
  \end{equation}
  to leading order as~$z\to-\infty,$ with
  $h_0(a):=2\Psi(1)-\Psi(a)-\Psi(1-a),$ where $\Psi$~is the Euler
  digamma function, i.e., $\Psi(a)=d(\log\Gamma(a))/da$.  Substituting the
  values of $\Psi(a),\Psi(1-a)$ for $a=\frac12,\frac13,\frac14,$ which are
  known~\cite[vol.~2,\ app.~II.3]{Prudnikov86}, reveals that the
  expressions in Table~\ref{tab:specfuncs} for~$h_N(t_N),$ $N=2,3,4,$ have
  leading-order behavior
  \begin{equation}
    h_N(t_N)\sim B_N\,t_N^{-\psi(N)/12}\,\log t_N
  \end{equation}
  as~$t_N\to\infty,$ where $B_N=4,3\sqrt3,8$ respectively.  One then uses
  $A_N=B_N/N,$ shown in the proof of Theorem~\ref{thm:clever}.
\end{proof}

The Picard--Fuchs equations $\mathcal{L}_N\,u=0$ of
Table~\ref{tab:fuchsianops} and the recurrences of
Table~\ref{tab:recurrences}, though `canonical', are not unique.  As was
explained in~\S\ref{sec:modular2}, to each alternative Hauptmodul~$\tilde
t_N$ of Table~\ref{tab:alt} there is associated an alternative
weight\nobreakdash-$1$ modular form~$\tilde h_N\circ t_N$.  For
convenience, write $t_N=\tilde\phi(\tilde t_N),$ where $\tilde\phi$~is the
appropriate M\"obius transformation.  The alternative modular form will
then be $(\tilde h_N\circ\tilde\phi)(\tilde t_N(\tau))$.  The function
$\tilde h_N\circ \tilde \phi$ will satisfy its own normal-form
Picard--Fuchs equation, with independent variable~$\tilde t_N$.  When
$N\le9,$ the unique holomorphic solution equal to unity at~$\tilde t_N=0$
will be expressible in~terms of ${}_2F_1$ or~$\Hl,$ providing an easy route
to the recurrence satisfied by the coefficients of its $\tilde
t_N$\nobreakdash-expansion.  But unlike the recurrences of
Table~\ref{tab:recurrences}, this recurrence will contain at~least one
minus sign.

The following pair of examples shows how in the just-described way, one may
expand $\tilde h_N(t_N(\tau))$ as a function of~$\tilde t_N(\tau),$
sometimes obtaining a combinatorially interesting series.  In the case
$N=6,$ consider the Hauptmoduln $\tilde t_6,\tilde{\tilde t}_6$ defined by
\begin{subequations}
\begin{align}
  \tilde t_6 &= t_6/(t_6+8),\qquad t_6=\tilde\phi(\tilde t_6):=8\,\tilde
  t_6/(1-\tilde t_6)\\ 
  \tilde{\tilde t}_6 &= t_6/(t_6+9),\qquad
  t_6=\tilde{\tilde\phi}(\tilde{\tilde t}_6):=9\,\tilde{\tilde
  t}_6/(1-\tilde{\tilde t}_6),
\label{eq:franelHlhauptmodul}
\end{align}
\end{subequations}
with a pole at $\bigl[\frac12\bigr]_6,$ resp.\ $\bigl[\frac13\bigr]_6$.
The associated modular forms for~$\Gamma_0(6)$ are
\begin{subequations}
\begin{align}
  \tilde h_6 &=\tilde h_6(t_6)=  [(t_6+8)/8]\,h_6(t_6),   \\
  \tilde{\tilde h}_6 &=\tilde{\tilde h}_6(t_6)= [(t_6+9)/9]\,h_6(t_6),
\label{eq:franelHlform}
\end{align}
\end{subequations}
with a zero at $\bigl[\frac12\bigr]_6,$ resp.\ $\bigl[\frac13\bigr]_6$.
Expressing $h_6$ in~terms of~$\Hl$ as in Table~\ref{tab:specfuncs}, and
employing~(\ref{eq:genPfaff}), the generalized Pfaff transformation
of~$\Hl,$ yields
\begin{subequations}
{\small
\begin{align}
(\tilde h_6\circ\tilde\phi)(\tilde t_6) &= \Hl (9,3;\,1,1,1,1;\,\tilde t_6) = \Hl (\tfrac19,\tfrac13;\,1,1,1,1;\,\tilde t_6/9),\\
(\tilde{\tilde h}_6\circ\tilde{\tilde\phi})(\tilde {\tilde t}_6) &= \Hl (-8,-2;\,1,1,1,1;\,\tilde{\tilde t}_6) = \Hl (-\tfrac18,\tfrac14;\,1,1,1,1;\,-\tilde{\tilde t}_6/8).
\label{eq:franelHl}
\end{align}
}%
\end{subequations}
If $\tilde h_6\circ\tilde\phi$ and $\tilde{\tilde
h}_6\circ\tilde{\tilde\phi}$ are expanded about the infinite cusp, as
$\sum_{n=0}^\infty \tilde c_n^{(6)} \tilde t^n$ and $\sum_{n=0}^\infty
\tilde {\tilde c}_n^{(6)} \tilde{\tilde t}^n,$ the coefficients will
satisfy the three-term recurrences
\begin{subequations}
\begin{align}
& n^2\,\tilde c^{(6)}_{n-1} - (10n^2+10n+3)\,\tilde c^{(6)}_n +  9(n+1)^2\,\tilde c^{(6)}_{n+1} = 0,\\
& n^2\,\tilde{\tilde c}^{(6)}_{n-1} + (7n^2+7n+2)\,\tilde{\tilde c}^{(6)}_n - 8(n+1)^2\,\tilde{\tilde c}^{(6)}_{n+1} = 0.
\end{align}
\end{subequations}
The scaled coefficients $\tilde d_n^{(6)}:=9^n\tilde c_n^{(6)},$
$\tilde{\tilde d}_n^{(6)}:=8^n\tilde{\tilde c}_n^{(6)}$ will be integral,
and one can show that $\tilde d_n^{(6)}=
\sum_{k=0}^n{\binom{n}{k}}^2\binom{2k}{k}$; and that the integers
$\tilde{\tilde d}_n^{(6)},$ $n\geq0,$ are the so-called Franel numbers,
i.e., $\tilde{\tilde d}_n^{(6)}= \sum_{k=0}^n{\binom{n}{k}}^3$.  These
recurrences and closed-form solutions were first derived by Stienstra and
Beukers~\cite{Stienstra85}, in a somewhat different framework (see also
Verrill~\cite{Verrill99}).  In Sloane's {\em Encyclopedia\/} the sequence
$\tilde d_n^{(6)},$ $n\geq0,$ is listed as {\tt A002893}, and $\tilde
{\tilde d}_n^{(6)},$ $n\geq0,$ as {\tt A000172}.

\section{Modular Equations for Elliptic Families}
\label{sec:families}

Now that the weight-$1$ modular forms $\mathfrak{h}_M= h_M\circ
t_M:\mathcal{H}^*\to\mathbb{P}^1(\mathbb{C})$ for the genus-zero modular
subgroups $\Gamma_0(M)<\Gamma(1)$ have been thoroughly examined, we proceed
to derive rationally parametrized modular equations for them.  When
$\Gamma_0(NM)$ is of genus zero, there is such a modular equation of
degree~$N,$ which links $h_M(t_M(\tau)),\allowbreak h_M(t_M(N\tau))$.
That~is, it links $h_M(t_M),h_M(t_M'),$ where the relation between
$t_M,t_M'\in X_0(M)\cong\mathbb{P}^1(\mathbb{C})$ is uniformized by the
Hauptmodul $t_{NM}\in X_0(NM),$ according to rational maps
$t_M=t_M(t_{NM}),$ $t_M'=t_M'(t_{NM})$ that can be found in
Table~\ref{tab:intermediate}.  The modular equation is of the form
$h_M(t_M)=\mathcal{M}\cdot h_M(t_M'),$ where the multiplier~$\mathcal{M}$
is (a~root~of) a function in~$\mathbb{Q}(t_{NM})$.

Such modular equations can be derived mechanically, with the aid of
$q$-series representations for $h_{M}(t_M(\tau))$ and~$t_{NM}(\tau)$.  That
is not our approach, since it ignores the algebraic--geometric context.
Our key Theorem~\ref{thm:key}, which efficiently produces such modular
equations (as transformations of special functions), is really a theorem on
pullbacks of Gauss--Manin connections along maps between rational elliptic
modular surfaces.  The following remarks place it in context.

An elliptic surface is a flat morphism
$\mathfrak{E}\stackrel{\pi}{\to}\mathbb{P}^1(\mathbb{C}),$ the generic
fibre of which is an elliptic curve $E/\mathbb{C}$.  An example is the
universal family of elliptic curves over~$\mathbb{C},$ denoted
$\hat{\mathfrak{E}}_1\stackrel{\pi_1}{\to}X(1)$ here, with base equal to
the $j$-line.  It has one singular fibre, over $j=\infty$.  For any
genus-zero $\Gamma<\Gamma(1),$ a rational elliptic surface
${\mathfrak{E}_\Gamma\to\Gamma\setminus\mathcal{H}^*},$ by definition a
rational elliptic {\em modular\/} surface, can be constructed from the
quotiented half-plane $\Gamma\setminus\mathcal{H}$ and elliptic fibres
over~it, by resolving singularities~\cite{Shioda72}.  If $\Gamma$ has
Hauptmodul~$t,$ with
\begin{equation}
j=\frac{P^3(t)S(t)}{R(t)}=12^3+\frac{Q^2(t)T(t)}{R(t)}
\end{equation}
(in~lowest terms,
as in~Table~\ref{tab:coverings}), then this surface, regarded as an
elliptic family parametrized by~$t,$ will have Weierstrass presentation
\begin{equation}
y^2=4x^3-3P(t)S(t)T(t)x-Q(t)S(t)T^2(t).
\end{equation}
There will be a singular fibre
over every fixed point of~$\Gamma$ on~$\Gamma\setminus\mathcal{H}^*$.

The surfaces $\mathfrak{E}_M\stackrel{\pi_M}{\to}X_0(M)$ attached to
genus-zero $\Gamma_0(M)$ are of primary interest.  For instance,
$\mathfrak{E}_5\stackrel{\pi_5}{\to}X_0(5),$ with parameter $t=t_5$ on the
base, has presentation
\begin{equation*}
y^2=4x^3 - 3(t^2+10t+5)(t^2+22t+125)x - (t^2+4t-1)(t^2+22t+125)^2.
\end{equation*}
This rational elliptic family, with four singular fibres, can be found in
the classification of Herfurtner~\cite{Herfurtner91}.  (The singular ones
are of Kodaira types $I_1,I_5,{\it III},{\it III},$ located above the two
cusps $t=0,\infty$ and the two quadratic elliptic fixed points
$t=-11\pm2{\rm i}$; cf.\ Table~\ref{tab:hauptmoduln}.)  The Hesse--Dixon
family~(\ref{eq:hessedixonmodel}) yields an elliptic surface with base
parameter $\gamma=\gamma(\tau)= t_9(\tau/3),$ and four singular fibres.  It
is attached to~$\Gamma(3),$ which is conjugated to~$\Gamma_0(9)$ by a
$3$-isogeny.  Verrill~\cite{Verrill2001} gives additional examples.

For any $M,N$ for which $\Gamma_0(M)$ and $\Gamma_0(NM)$ are of genus zero,
there is a commutative diagram of rational elliptic surfaces,
\begin{equation}
\begin{CD}
  \mathfrak{E}_{NM} @>>> \mathfrak{E}_{M} @>>> \hat{\mathfrak{E}}_{1} \\
  @V{\pi_{NM}}VV	  @V{\pi_M}VV	  @V{{\pi}_1}VV \\
  X_0(NM) @>>> X_0(M) @>>> X(1)
\end{CD}
\end{equation}
where the maps $X_0(NM)\to X_0(M)$ and $X_0(M)\to X(1)$ are given by the
rational functions $t_{NM}\mapsto t_M$ and $t_{M}\mapsto j,$ listed in
Tables \ref{tab:intermediate} and~\ref{tab:jfactored} respectively, and the
fibre maps are induced by them.  (If $t_{NM}\mapsto t_M$ is replaced by
$t_{NM}\mapsto t_M',$ i.e., $t_{NM}(\tau)\mapsto t_M(N\tau),$ or
$t_{M}\mapsto j$ by $t_{M}\mapsto j',$ i.e., $t_{M}(\tau)\mapsto j(M\tau),$
then the resulting diagram will still commute.)

Now consider the Picard--Fuchs equations attached to the genus-zero
subgroups $\Gamma_0(M)$.  At any generic point~$t_M$ on the base of
$\mathfrak{E}_M\stackrel{\pi_M}{\to}X_0(M),$ the fibre (an~elliptic curve)
will have period module $\mathbb{Z}\tau_1\oplus\mathbb{Z}\tau_2,$ where
$\tau_1/\tau_2$ and~$\tau_2$ are what we have been calling $\tau$
and~$h_M(t_M(\tau))$.  The space of local solutions of the Picard--Fuchs
equation $\mathcal{L}_M u = 0$
on~$X_0(M)\cong\mathbb{P}^1(\mathbb{C})_{t_M}$ is
$h_M(\cdot)\left[\mathbb{C}\tau(\cdot)+\mathbb{C}\right],$ i.e.,
$\mathbb{C}\tau_1(\cdot)+\mathbb{C}\tau_2(\cdot)$.

Any element $\omega\in H^1({\pi_M}^{-1}(t_M))$ of the de~Rham cohomology
group of a fibre of~$\mathfrak{E}_M$ above a generic point $t_M\in X_0(M)$
is a period, i.e., a point in~$\mathbb{C}\tau_1+\mathbb{C}\tau_2$.  So the
second-order differential equation $\mathfrak{L}_Mu=0$ defines a flat
connection~$\nabla_M$ on a $2$-dimensional period bundle, each fibre being
a de~Rham cohomology group; viz., a~Gauss--Manin connection.

In~\S\ref{sec:modular3}, each Picard--Fuchs operator~$\mathcal{L}_M$ was
constructed explicitly, by pulling back along $X_0(M)\to X(1)$ the Gauss
hypergeometric operator~$\hat{\mathcal{L}}_1$ that defines a Gauss--Manin
connection for the universal family
$\hat{\mathfrak{E}}_1\stackrel{\pi_1}{\to} X(1)$.  (Much of the work went
into keeping $\mathcal{L}_M$ in normal form, to permit the derivation of
${}_2F_1$ and other special function identities.)  So by construction, the
connection~$\nabla_1$ on~$\hat{\mathfrak{E}}_1$ pulls back to~$\nabla_M$
on~${\mathfrak{E}}_M,$ which in~turn pulls back to~$\nabla_{NM}$
on~${\mathfrak{E}}_{NM}$.

Theorem~\ref{thm:key}, which in order to facilitate computation is phrased
in~terms of the respective modular forms $\hat h_1(\hat
J(\tau)),\allowbreak h_M(t_M(\tau)),\allowbreak h_{NM}(t_{NM}(\tau))$ for
$\Gamma(1),\allowbreak\Gamma_0(M),\allowbreak\Gamma_0(NM),$ says the
following.  (1)~Pulling back~$\nabla_1$ along the distinct maps
$t_{M}(\tau)\mapsto j(\tau)$ and $t_{M}(\tau)\mapsto j(M\tau)$ yields the
{\em same\/} connection; namely,~$\nabla_{M}$.  (2)~Pulling back~$\nabla_M$
along the distinct maps $t_{NM}(\tau)\mapsto t_M(\tau)$ and
$t_{NM}(\tau)\mapsto t_M(N\tau)$ does too, in both cases
yielding~$\nabla_{NM}$.

\begin{theorem}
  \label{thm:key}
\ \begin{enumerate}
    \item If\/ $\Gamma_0(M)$ is of genus zero, and if one has\/ {\rm(}in
      lowest terms{\rm)} that\hfil\break $j=P(t_M)/Q(t_M),$ $j'=P'(t_M)/Q'(t_M),$
      then the identity
      \begin{align*}
	h_{M}(t_{M}) &= [P(t_{M})/P(0)]^{-1/12}\, \hat h_1(P(t_{M})/Q(t_{M}))\\
	             &= [P'(t_{M})/P'(0)]^{-1/12}\, \hat h_1(P'(t_{M})/Q'(t_{N}))
      \end{align*}
      holds near the point\/ $t_{M}=0$ on\/ $X_0(M),$ i.e., near the
      infinite cusp.
    \item If\/ $\Gamma_0(M),\Gamma_0(NM)$ are of genus zero, and
      if\/ {\rm(}in lowest terms{\rm)}\hfil\break $t_M=P(t_{NM})/Q(t_{NM}),$
      $t_M'=P'(t_{NM})/Q'(t_{NM}),$ then the identity
      \begin{align*}
	h_{NM}(t_{NM}) &= [Q(t_{NM})/Q(0)]^{-\psi(M)/12}\, h_M(P(t_{NM})/Q(t_{NM}))\\
	               &= [Q'(t_{NM})/Q'(0)]^{-\psi(M)/12}\, h_M(P'(t_{NM})/Q'(t_{NM}))
      \end{align*}
      holds in a neighborhood of the point\/ $t_{NM}=0$ on\/ $X_0(NM),$
      i.e., of the infinite cusp.
  \end{enumerate}
\end{theorem}
\begin{proof}
First, recall that by Theorems~\ref{thm:3} and~\ref{thm:4}, $\hat h_1=\hat
h_1(\hat J)$ and $h_M=h_M(t_M)$ are the unique local solutions of the
canonical Picard--Fuchs equations $\hat{\mathcal{L}}_1 u=0$ for $X(1)$ and
$\mathcal{L}_M u=0$ for $X_0(M)$ that are holomorphic at the infinite cusp
(i.e., at~$\hat J=0,$ resp.\ $t_M=0$), and are normalized to equal unity
there.

The first equality in part~1 is true by definition (see
Definition~\ref{def:2}).  To~prove the first in part~2, pull back
$\mathcal{L}_M$ along $X_0(NM)\to X_0(M),$ i.e., along $t_M=t_M(t_{NM}),$
which takes $0$ to~$0$.  The pulled-back operator $(\mathcal{L}_M)^*$ will
be a Picard--Fuchs operator for~$\Gamma_0(NM),$ but it may not be the
unique normal-form one~$\mathcal{L}_{NM}$.  By Theorem~\ref{thm:4},
$\mathcal{L}_M$~has exponents $\frac1{12}\psi(M),\frac1{12}\psi(M)$ at the
cusp $t_M=\infty$ (i.e., $\tau=0$) and $0,0$ at other cusps.  If the
function $t_M=t_M(t_{NM})$ is not a polynomial, there will be a cusp other
than $t_{NM}=\infty$ in the inverse image of $t_M=\infty,$ and there will
be a disparity in exponents: at~any such cusp, the operator
$\mathcal{L}_{NM}$ will have exponents~$0,0,$ but $(\mathcal{L}_M)^*$~will
have exponents $\frac{k}{12}\psi(M),\frac{k}{12}\psi(M),$ where $k$~is the
multiplicity with which the cusp appears in the fibre.  Applying a
similarity transformation to $(\mathcal{L}_M)^*,$ by performing the change
of variable (substitution) $\hat u =[Q(t_{NM})/Q(0)]^{-\psi(M)/12}\,u,$
will remove this disparity, for all such cusps.  By Theorem~\ref{thm:2},
the resulting transformed operator must equal~$\mathcal{L}_{NM},$ which has
unique normalized holomorphic solution~$h_{NM}$ in a neighborhood
of~$t_{NM}=0$; so the left and right sides are equal as claimed.

The second equality in each of parts~$1,2$ is proved in a related way.
Consider $\hat{\mathcal L}_1'u=0$ and $\mathcal{L}_M'u=0,$ the Fuchsian
differential equations on $X(1)'$ and~$X_0(M)$ obtained by formally
substituting~$\hat J'$ (or~$j'$) and~$t_M'$ for~$\hat J$ (or~$j$)
and~$t_M,$ in $\hat{\mathcal L}_1u=0$ and $\mathcal{L}_Mu=0$ respectively.
Any ratio of independent solutions of $\hat{\mathcal{L}}_1u=0$ or
$\mathcal{L}_Mu=0$ is of the form $(a\tau+b)/(c\tau+d),$ where
$a,b,c,d\in\mathbb{C}$ with $ad\neq bc$; so the same is true of
$\hat{\mathcal{L}}_1'u=0$ and $\mathcal{L}_M'u=0$.  (In~effect, adding
primes multiplies the coefficients $a,c$ by~$N$.)  The same is necessarily
true of their respective pullbacks to $X_0(M)$ and~$X_0(NM),$ by the
definition of a pullback.  But the pullbacks may not equal $\mathcal{L}_M$
and~$\mathcal{L}_{NM}$ respectively, due~to the problem encountered in the
last paragraph: exponents may not agree, even though exponent {\em
differences\/} must.  A similarity transformation by
$[P'(t_{M})/P'(0)]^{-1/12},$ resp.\ $[Q'(t_{NM})/Q'(0)]^{-\psi(M)/12}$ will
solve this problem.  If the exponents agree then the Picard--Fuchs
operators and equations must agree, by Theorem~\ref{thm:2}; and the claimed
equalities follow from the uniqueness of the local solution $h_M=h_M(t_M),$
resp.\ $h_{NM}=h_{NM}(t_{NM})$.
\end{proof}

\begin{corollary}
  If all primes that divide\/ $N$ also divide\/~$M,$ then the modular
forms\/ $\mathfrak{h}_M = h_M\circ t_M$ and\/ $\mathfrak{h}_{NM} =
h_{NM}\circ t_{NM}$ are equal.
\end{corollary}
\begin{proof}
  Under the divisibility assumption, $t_M=P(t_{NM})/Q(t_{NM})$ reduces to a
  polynomial map $t_M=P(t_{NM})$ by Proposition~\ref{prop:short}; so by
  part~2 of the theorem, $h_{NM}(t_{NM}(\tau)) =
  h_N(P(t_{NM}(\tau)))=h_M(t_M(\tau))$.
\end{proof}
\begin{remark*}
  This corollary reveals why there are duplications in column~3 of
  Table~\ref{tab:xNhN}, which lists the eta product representations of the
  fourteen modular forms $\mathfrak{h}_M = h_M\circ t_M$.  The seven
  duplicates are due~to the seven identities
  \begin{subequations}
    \def\theequation{\thesection.\arabic{parentequation}\alph{equation}}
  \begin{align}
    \def\theequation{\arabic{parentequation}\alpha{equation}}
    \label{eq:convoluted1a}
    h_4(t)&= h_2(t(t+16)),\\
    \label{eq:convoluted1b}
    h_8(t)&= h_2(t(t+16)\circ t(t+8)),\\
    \label{eq:convoluted1c}
    h_{16}(t)&= h_2(t(t+16)\circ t(t+8)\circ t(t+4)),\displaybreak[0]\\[\jot]
    \addtocounter{parentequation}{1}
    \setcounter{equation}{0}
    \def\theequation{\arabic{parentequation}\alpha{equation}}
    \label{eq:convoluted2a}
    h_9(t)&= h_3(t(t^2+9t+27)),\displaybreak[0]\\[\jot]
    \addtocounter{parentequation}{1}
    \setcounter{equation}{0}
    \def\theequation{\arabic{parentequation}\alpha{equation}}
    h_{12}(t)&= h_6(t(t+6)),\\
    h_{18}(t)&= h_6(t(t^2+6t+12)),\\[\jot]
    \addtocounter{parentequation}{1}
    \setcounter{equation}{0}
    \def\theequation{\arabic{parentequation}\alpha{equation}}
    h_{25}(t)&= h_5(t(t^4+5t^3+15t^2+25t+25)).
  \end{align}
  \end{subequations}
  The identities
  (\ref{eq:convoluted1a})\textendash\nobreak(\ref{eq:convoluted1c})
  and~(\ref{eq:convoluted2a}) clarify the relations among the quadratic AGM
  representations, Eqs.\ (\ref{eq:agma})\textendash\nobreak(\ref{eq:agmd}),
  and the cubic ones, Eqs.\
  (\ref{eq:agm3a})\textendash\nobreak(\ref{eq:agm3b}), respectively.
\end{remark*}

\begin{table}
\caption{The function $\Delta(N\tau)/\Delta(\tau),$ modular
for~$\Gamma_0(N),$ in~terms of~$t:=t_N$.}  
\begin{center}
{\small
\begin{tabular}{cl}
\hline\noalign{\smallskip}
$N$ & \hfil $N^{12}\Delta(N\tau)/\Delta(\tau)$ \\
\noalign{\smallskip}\hline\noalign{\smallskip}
$2$ & $t\vphantom{\frac{t^{15}(t+4)^3}{(t+2)^3}}$ \\
$3$ & $t^2\vphantom{\frac{t^{15}(t+4)^3}{(t+2)^3}}$ \\
$4$ & $t^3\vphantom{\frac{t^{15}(t+4)^3}{(t+2)^3}}$ \\
$5$ & $t^4\vphantom{\frac{t^{15}(t+4)^3}{(t+2)^3}}$ \\
$6$ & $\frac{t^5(t+9)}{t+8}\vphantom{\frac{t^{15}(t+4)^3}{(t+2)^3}}$ \\
$7$ & $t^6\vphantom{\frac{t^{15}(t+4)^3}{(t+2)^3}}$ \\
$8$ & $\frac{t^7(t+8)}{t+4}\vphantom{\frac{t^{15}(t+4)^3}{(t+2)^3}}$ \\
\noalign{\smallskip}\hline
\end{tabular}
\hfil\hfil
\begin{tabular}{cl}
\hline\noalign{\smallskip}
$N$ & \hfil $N^{12}\Delta(N\tau)/\Delta(\tau)$ \\
\noalign{\smallskip}\hline\noalign{\smallskip}
$9$ & $t^8\vphantom{\frac{t^{15}(t+4)^3}{(t+2)^3}}$ \\
$10$& $\frac{t^9(t+5)^3}{(t+4)^3}\vphantom{\frac{t^{15}(t+4)^3}{(t+2)^3}}$ \\
$12$& $\frac{t^{11}(t+4)(t+6)^2}{(t+2)^2(t+3)}\vphantom{\frac{t^{15}(t+4)^3}{(t+2)^3}}$ \\
$13$& $t^{12}\vphantom{\frac{t^{15}(t+4)^3}{(t+2)^3}}\vphantom{\frac{t^{15}(t+4)^3}{(t+2)^3}}$ \\
$16$& $\frac{t^{15}(t+4)^3}{(t+2)^3}\vphantom{\frac{t^{15}(t+4)^3}{(t+2)^3}}$ \\
$18$& $\frac{t^{17}(t+3)^7(t^2+6t+12)}{(t+2)^7(t^2+3t+3)}\vphantom{\frac{t^{15}(t+4)^3}{(t+2)^3}}$ \\
$25$& $t^{24}\vphantom{\frac{t^{15}(t+4)^3}{(t+2)^3}}$ \\
\noalign{\smallskip}\hline
\end{tabular}
}%
\end{center}
\label{tab:delta}
\end{table}

\begin{corollary}
  For each\/~$N$ with\/ $\Gamma_0(N)$ of genus zero, the modular
  discriminant\/ $\Delta$ satisfies the modular equations given in
  Table\/~{\rm\ref{tab:delta},} which express\/
  $\Delta(N\tau)/\Delta(\tau)$ in~terms of the Hauptmodul\/~$t_N$
  for\/~$\Gamma_0(N)${\rm.}
\end{corollary}
\begin{proof}
  For each $N,$ the modular equation for $E_4=(\hat h_1 \circ \hat J)^4$
  follows from part~1 of Theorem~\ref{thm:key}, if one uses the formulas
  for $j=j(t_{N}),\allowbreak j'=j'(t_N)$ given in Tables
  \ref{tab:coverings} and~\ref{tab:coverings2}.  Also, by
  Theorem~\ref{thm:3}, $E_6=(1-\hat J)^{1/2}(\hat h_1 \circ \hat J)^6$ and
  $\Delta=(2\pi)^{12} 12^{-3}\hat J\, (\hat h_1\circ \hat J)^{12},$ so the
  modular equations for~$E_6,\Delta$ follow from some elementary further
  manipulations, using the fact that $\hat J=12^3/j$.  The modular
  equations for~$\Delta$ are given in the table; those for $E_4,E_6$ are
  omitted.
\end{proof}
\begin{remark*}
  The expressions for $\Delta(N\tau)/\Delta(\tau)$ in Table~\ref{tab:delta}
  seem not to have been tabulated before.  It is well known that if
  $N-1\divides{24},$ then $N^{12}\Delta(N\tau)/\Delta(\tau)$ equals our
  canonical Hauptmodul~$t_N$.  (See Apostol~\cite[Ch.~4]{Apostol90}.)  But
  the table reveals what happens in the six cases when
  $N-1\notdivides{24}$.

  Each expression for $\Delta(N\tau)/\Delta(\tau)$ in the table
  can alternatively be verified by rewriting it as an eta product, since
  $\Delta(N\tau)/\Delta(\tau) = [N]^{24}/\,[1]^{24}$.  To do
  this, one would exploit the eta product representations given in Tables
  \ref{tab:hauptmoduln} and~\ref{tab:alt}.  To~handle the case $N=18,$ one
  would also need the eta product representations (\ref{eq:special18a})
  and~(\ref{eq:special18b}) for the bivalent functions
  $t_{18}^2+6t_{18}+12$ and $t_{18}^2+3t_{18}+3$ on~$X_0(18),$ each of
  which is zero on a conjugate pair of irrational cusps.
\end{remark*}

The following theorem is a slight extension of Theorem~\ref{thm:key}.  It
is proved in the same way.
\addtocounter{theorem}{-1}
\def\thetheorem{\arabic{section}.\arabic{theorem}${}'$}
\begin{theorem}
\ \begin{enumerate}
    \item If\/ $\Gamma_0(N)$ is of genus zero, and if for some\/
      $d\divides{N},$ the function\/ $j(d\tau)$ can be expressed rationally
      in~terms of\/ $t_N(\tau)$ {\rm(}in lowest terms{\rm)} by the
      formula\/ $j(d\tau)=P(t_N(\tau))/Q(t_N(\tau)),$ then
      \begin{equation*}
	h_{N}(t_{N}) = [Q(t_{N})/Q(0)]^{-1/12}\, \hat h_1(P(t_{N})/Q(t_{N}))
      \end{equation*}
      holds in a neighborhood of the point\/ $t_{N}=0${\rm.}
  \item If\/ $\Gamma_0(M),\Gamma_0(NM)$ are of genus zero, and if for
    some\/~$d\divides{N},$ the function\/ $t_M(d\tau)$ can be expressed
    rationally in~terms of\/ $t_{NM}(\tau)$ by the formula\/
    $t_M(d\tau)=P(t_{NM}(\tau))/Q(t_{NM}(\tau))$ {\rm(}in lowest
    terms{\rm)}, then
    \begin{equation*}
      h_{NM}(t_{NM}) = [Q(t_{NM})/Q(0)]^{-\psi(M)/12}\, h_N(P(t_{NM})/Q(t_{NM}))
      \end{equation*}
      holds in a neighborhood of the point\/ $t_{NM}=0${\rm.}
  \end{enumerate}
\label{thm:extendedkey}
\end{theorem}
\def\thetheorem{\arabic{section}.\arabic{theorem}}

\begin{table}
\caption{Rationally parametrized modular equations for the weight-$1$
modular forms~$\mathfrak{h}_M(\tau)$ for~$\Gamma_0(M),$ viewed as
functions~$h_M=h_M(t_M)$ of the Hauptmoduln~$t_M$.  In each, $t:=t_{NM}$.}
\begin{center}
{\small
\begin{tabular}{ccl}
\hline\noalign{\smallskip}
$M$ & $N$ & Modular equation for~$h_M,$ of degree~$N$\\
\noalign{\smallskip}\hline\noalign{\vskip1pt}\hline\noalign{\vskip1.5pt}
$2$ & $2$ & $h_2\left(t(t+16)\right)=2\,(t+16)^{-1/4}\,h_2\left(\frac{t^2}{t+16}\right)$\\
$2$ & $3$ & $h_2\left(\frac{t(t+8)^3}{t+9}\right) = 3\,(t+9)^{-1/2} \,h_2\left(\frac{t^3(t+8)}{(t+9)^3}\right)$ \\
$2$ & $4$ & $h_2\left(t(t+16)\,\circ\,t(t+8)\right) = 4\,\left[(t+4)(t+8)\right]^{-1/4}\,h_2\left(\frac{t^2}{t+16} \circ \frac{t^2}{t+4}\right)$ \\
$2$ & $5$ & $h_2\left(\frac{t(t+4)^5}{t+5}\right)=5\,(t+5)^{-1}\,h_2\left(\frac{t^5(t+4)}{(t+5)^5}\right)$ \\
$2$ & $6$ & $h_2\left(\frac{t(t+8)^3}{t+9}\circ\,t(t+6)=t(t+16)\, \circ\frac{t(t+4)^3}{t+3}\right)$\\
    &     & \quad{\scriptsize$=6\,\left[(t+2)(t+3)(t+6)^3\right]^{-1/4}\,h_2\left(\frac{t^3(t+8)}{(t+9)^3}\circ\frac{t^2}{t+2}=\frac{t^2}{t+16}\circ\frac{t^3(t+4)}{(t+3)^3}\right)$}\\
$2$ & $8$ & $h_2\left(t(t+16)\,\circ\,t(t+8)\,\circ\,t(t+4)\right)$\\
    &     & \quad $=8\,\left[(t+2)(t+4)^4(t^2+4t+8)\right]^{-1/4}\,h_2\left(\frac{t^2}{t+16}\circ\frac{t^2}{t+4}\circ\frac{t^2}{t+2}\right)$\\
$2$ & $9$ & $h_2\left(\frac{t(t+8)^3}{t+9}\circ\,t(t^2+6t+12)\right) = 9\,(t+3)^{-2}\,h_2\left(\frac{t^3(t+8)}{(t+9)^3}\circ\frac{t^3}{t^2+3t+3}\right)$ \\
\noalign{\vskip 1pt}\hline\noalign{\vskip 1pt}
$3$ & $2$ & $h_3\left(\frac{t(t+9)^2}{t+8}\right)= 2\,(t+8)^{-1/3}\,h_3\left(\frac{t^2(t+9)}{(t+8)^2}\right)$\\
$3$ & $3$ & $h_3\left(t(t^2+9t+27)\right)= 3\,(t^2+9t+27)^{-1/3}\,h_3\left(\frac{t^3}{t^2+9t+27}\right)$\\
$3$ & $4$ & $h_3\left(\frac{t(t+9)^2}{t+8}\circ\, t(t+6)\right)= 4\,(t+4)^{-1}\,h_3\left(\frac{t^2(t+9)}{(t+8)^2}\circ \frac{t^2}{t+2}\right)$\\
$3$ & $6$ & $h_3\left(t(t^2+9t+27)\,\circ\frac{t(t+3)^2}{t+2}=\frac{t(t+9)^2}{t+8}\circ\,t(t^2+6t+12)\right)$\\
    &     & \quad{\scriptsize$=6\,\left[\frac{(t^3+3t+3)(t^2+6t+12)^2}{t+2}\right]^{-1/3}h_3\left(\frac{t^3}{t^2+9t+27}\circ\frac{t^2(t+3)}{(t+2)^2}=\frac{t^2(t+9)}{(t+8)^2}\circ\frac{t^3}{t^2+3t+3}\right)$} \\
\noalign{\vskip 1pt}\hline\noalign{\vskip 1pt}
$4$ & $2$ & $h_4\left(t(t+8)\right) = 2\,(t+4)^{-1/2}\,h_4\left(\frac{t^2}{t+4}\right)$\\
$4$ & $3$ & $h_4\left(\frac{t(t+4)^3}{t+3}\right) = 3\,(t+3)^{-1}\,h_4\left(\frac{t^3(t+4)}{(t+3)^3}\right)$\\
$4$ & $4$ & $h_4\left(t(t+8)\,\circ\,t(t+4)\right)$\\
    &     & \quad{$= 4\,\left[(t+2)(t^2+4t+8)\right]^{-1/2}\,h_4\left(\frac{t^2}{t+4}\circ\frac{t^2}{t+2}\right)$}\\
\noalign{\vskip 1pt}\hline\noalign{\vskip 1pt}
$5$ & $2$ & $h_5\left(\frac{t(t+5)^2}{t+4}\right) = 2\,(t+4)^{-1/2}\,h_5\left(\frac{t^2(t+5)}{(t+4)^2}\right)$ \\
$5$ & $5$ & $h_5\left(t(t^4+5t^3+15t^2+25t+25)\right)$\\
    &     & \quad $= 5\,\left[t^4+5t^3+15t^2+25t+25\right]^{-1/2}\,h_5\left(\frac{t^5}{t^4+5t^3+15t^2+25t+25}\right)$ \\
\noalign{\vskip 1pt}\hline\noalign{\vskip 1pt}
$6$ & $2$ & $h_6\left(t(t+6)\right)=2\,(t+2)^{-1}\,h_6\left(\frac{t^2}{t+2}\right)$\\
$6$ & $3$ & $h_6\left(t(t^2+6t+12)\right)=3\,(t^2+3t+3)^{-1}\,h_6\left(\frac{t^3}{t^2+3t+3}\right)$\\
\noalign{\vskip 1pt}\hline\noalign{\vskip 1pt}
$8$ & $2$ & $h_8\left(t(t+4)\right)=2\,(t+2)^{-1}\,h_8\left(\frac{t^2}{t+2}\right)$\\
\noalign{\vskip 1pt}\hline\noalign{\vskip 1pt}
$9$ & $2$ & $h_9\left(\frac{t(t+3)^2}{t+2}\right)=2\,(t+2)^{-1}\,h_9\left(\frac{t^2(t+3)}{(t+2)^2}\right)$\\
\noalign{\vskip2pt}\hline
\end{tabular}
}
\end{center}
\label{tab:keymodular}
\end{table}

\begin{theorem}
  For each\/~$M$ with\/ $\Gamma_0(M)$ of genus zero, the function\/~$h_M$
  defined in a neighborhood of\/~$t_M=0$ on\/~$X_0(M)$ {\rm(}i.e., of the
  infinite cusp\/{\rm)} satisfies the rationally parametrized degree-$N$
  modular equations given in Table\/~{\rm\ref{tab:keymodular},} of the
  form\/ $h_M(t_M(t_{NM}))=\mathcal{M}_{M,N}(t_{NM})\cdot
  h_M(t'_M(t_{NM}))${\rm.}
  \label{thm:maintable}
\end{theorem}
\begin{proof}
  For each $M,N,$ the claimed modular equation follows from
  Theorem~\ref{thm:key} (part~2), if one takes into account the formulas
  $t_M=t_M(t_{NM}),\allowbreak t_M'=t_M'(t_{NM})$ given in
  Table~\ref{tab:intermediate}.
\end{proof}
\begin{altproof}
  For each $M,N,$ evaluate both sides as eta products, and verify that they
  are the same.  One needs the eta product representations for the modular
  forms~$\mathfrak{h}_M$ given in Table~\ref{tab:xNhN}, and those for
  Hauptmoduln and alternative Hauptmoduln given in Tables
  \ref{tab:hauptmoduln} and~\ref{tab:alt}.  To~handle the cases
  $M=9,16,18,25,$ one also needs the special eta product formulas for
  certain non-univalent functions of~$t_M$ given in
  (\ref{eq:gatherstart})\textendash\nobreak(\ref{eq:gatherend}).  The
  computations are straightforward and are left to the reader.  But this
  alternative way of proceeding is a {\em verification\/}, not a proof: it
  sheds little light on the algebraic--geometric origin of the modular
  equations.
\end{altproof}
\begin{remark*}
  At a general point $\tau\in\mathcal{H},$ each equation in
  Table~\ref{tab:keymodular} evaluates to
  \begin{equation}
    (h_M\circ t_M)(\tau) = \mathcal{M}_{M,N}(t_{NM}(\tau))\cdot(h_M\circ
    t_M)(N\tau),
  \end{equation}
  and can be called an $N:1$ modular equation.  The degree-$N$ multiplier
  $\mathcal{M}_{M,N}=\mathcal{M}_{M,N}(t_{NM})$ specifies the relationship
  between the period modules of the fibres (elliptic curves) over
  $t_M,t_M'\in X_0(M),$ for the elliptic surface
  $\mathfrak{E}_M\stackrel{\pi_M}{\to}X_0(M)$.  It~is algebraic over each
  of the related points $t_M,t_M'$ on the base.
\end{remark*}
\begin{remark*}
  Each equation in Table~\ref{tab:keymodular} holds for all $t>0,$ if one
  interprets $h_M=h_M(t_M)$ as the unique holomorphic continuation along
  the positive $t_M$-axis.  The two argument functions in each equation,
  i.e., $t=t_{NM}\mapsto t_M$ and $t=t_{NM}\mapsto t'_M,$ by examination
  take $(0,\infty)$ bijectively to~$(0,\infty)$.
\end{remark*}

There are two additional modular equations that one can derive, using not
Theorem~\ref{thm:key} but its enhancement Theorem~\ref{thm:extendedkey}.
\begin{theorem}
  The functions\/ $h_2=h_2(t_2),h_3=h_3(t_3),$ in neighborhoods of the
  points\/ $t_2=0,\allowbreak t_3=0,$ satisfy respective functional
  equations
  {\small
  \begin{align}
    & h_2\left(\tfrac{t^2}{t+16}\circ\tfrac{t(t+4)^3}{t+3} = \tfrac{t(t+8)^3}{t+9}\circ\tfrac{t^2}{t+2}\right) \nonumber\\
    & \qquad =\left(\tfrac32\right) \left[\tfrac{(t+3)^5}{(t+2)^3(t+6)}\right]^{-1/4}\! h_2\left(t(t+16)\,\circ\tfrac{t^3(t+4)}{(t+3)^3} = \tfrac{t^3(t+8)}{(t+9)^3}\circ\,t(t+6)\right), \nonumber\\
    & h_3\left(t(t^2+9t+27)\,\circ\tfrac{t^2(t+3)}{(t+2)^2}=\tfrac{t^2(t+9)}{(t+8)^2}\circ\, t(t^2+6t+12)\nonumber\right) \\
    & \qquad = \left(\tfrac32\right) \left[\tfrac{(t^2+3t+3)^2(t^2+6t+12)}{(t+2)^5}\right]^{-1/3} \nonumber \\
    & \qquad\qquad\qquad\qquad\qquad\qquad
     {}\times h_3\left(\tfrac{t^3}{t^2+9t+27}\circ\tfrac{t(t+3)^2}{t+2}=\tfrac{t(t+9)^2}{t+8}\circ\tfrac{t^3}{t^2+3t+3}\right),\nonumber
  \end{align}
  }%
  which can be called modular equations of degree\/ $3:2,$ since they
  relate\hfil\break $(h_2\circ t_2)(2\tau)$ to\/ $(h_2\circ t_2)(3\tau),$ and\/
  $(h_3\circ t_3)(2\tau)$ to\/ $(h_3\circ t_3)(3\tau)$.
\label{thm:supplement}
\end{theorem}
\begin{proof}
In these two equations, the parameter~$t$ signifies respectively the
Hauptmoduln $t_{12}$ and~$t_{18}$.  The first equation follows from the
$M=2,\allowbreak N=6$ case of Theorem~\ref{thm:extendedkey} (part~2), by
equating the expressions for $h_{12}=h_{12}(t_{12})$ obtained from the
cases $d=2,3$.  The second follows similarly from the $M=3,\allowbreak N=6$
case, by equating expressions for $h_{18}=h_{18}(t_{18})$.
\end{proof}

Each modular equation for~$h_M$ in Table~\ref{tab:keymodular} (and
Theorem~\ref{thm:supplement}) with $M=2,3,4$ can be written as an algebraic
transformation of~${}_2F_1,$ with the aid of the expressions for~$h_M$
given in Table~\ref{tab:specfuncs}.  These transformations will be
discussed in~\S\ref{sec:ramanujan}.  The modular equations with
$M=5,6,\allowbreak7,\allowbreak8,9$ can similarly be written as algebraic
transformations of the local Heun function~${\it Hl}$.  As special function
identities, they are quite novel, since there is not even a rudimentary
theory of Heun transformations (unlike the classical theory of
transformations of~${}_2F_1$).  It is possible to restate them as
functional equations for combinatorial generating functions, without
explicit reference to~${\it Hl}$.  The following striking proposition
illustrates this.
\begin{proposition}
\label{prop:franel}
Let\/ $F=F(x)=\sum_{n=0}^\infty a_nx^n$ be the generating function of the
Franel numbers\/ $a_n=\sum_{k=0}^n \binom{n}{k}^3,$ $n\geq0$.  Then\/
$F,$~which is defined on the disk\/ $\left|x\right|<1/8,$ satisfies the
quadratic and cubic functional equations
\begin{align*}
  F\left(\tfrac{t(t+6)}{8(t+3)^2}\right) &= 2\left[\tfrac{t+3}{t+6}\right]\, F\left(\tfrac{t^2}{8(t+3)(t+6)}\right), \\
  F\left(\tfrac{t(t^2+6t+12)}{8(t+3)(t^2+3t+3)}\right) &= 3\left[\tfrac{t^2+3t+3}{(t+3)^2}\right]\, F\left(\tfrac{t^3}{8(t+3)^3}\right),
\end{align*}
for\/ $\left|t\right|$ sufficiently small, and also for all\/ $t>0$.
\end{proposition}
\begin{proof}
  The Franel numbers $a_n=\tilde{\tilde d}_n^{(6)}$ were introduced
  in~\S\ref{sec:modular3}, and it follows from
  (\ref{eq:franelHlhauptmodul}), (\ref{eq:franelHlform}),
  and~(\ref{eq:franelHl}) that (changing the dummy variable~$t$ to~$t_6$)
  \begin{equation}
    \begin{split}
    F(t_6) &= (\tilde{\tilde h}_6\circ\tilde{\tilde\phi})(8t_6) = (1-8t_6)^{-1}h_6\left(\tfrac{72t_6}{1-8t_6}\right)\\
         &= \Hl (-8,-2;\,1,1,1,1;\,8t_6) = \Hl (-\tfrac18,\tfrac14;\,1,1,1,1;\,-t_6).
    \end{split}
  \end{equation}
The functional equations accordingly follow from the modular equations
for~$h_6$ of degrees $N=2,3$ in Table~\ref{tab:keymodular}.  The
parameter~$t$ is respectively $t_{12},t_{18}$.
\end{proof}

The lone modular equation for~$h_9$ in Table~\ref{tab:keymodular}, of
degree $N=2,$ implies the following proposition on the behavior of the
periods of the Hesse--Dixon family
\begin{equation}
  x^3+y^3 + 1 - (\gamma+3)xy = 0,
\end{equation}
parametrized by
$\gamma\in\mathbb{C}\setminus\{0,3(\zeta_3-1),3(\zeta^2_3-1)\},$ under the
action of the level-$2$ modular correspondence that doubles the period
ratio $\tau=\tau_1/\tau_2$.
\begin{proposition}
  Let\/ $E$ and\/~$E'$ be Hesse--Dixon elliptic curves with respective
  parameters\/ $\gamma$ and\/~$\gamma',$ which are related parametrically
  by
  \begin{equation*}
    \gamma=\frac{t(t+3)^2}{t+2},\qquad \gamma'=\frac{t^2(t+3)}{(t+2)^2}.
  \end{equation*}
  Then their respective fundamental periods\/ $\tau_1,\tau_2$ and
  $\tau_1',\tau_2'$ will be related by\hfil\break $\tau_1'=(t+2)\tau_1,$
  $\tau_2'=[(t+2)/2]\tau_2$.
\end{proposition}

\section{Ramanujan's Modular Equations in Signatures $2,3,4$}
\label{sec:ramanujan}

The modular equations derived in~\S\ref{sec:families}, for the elliptic
families $\mathfrak{E}_M\stackrel{\pi_M}{\to}X_0(M),$ $M=2,3,4,$ yield
rationally parametrized modular equations in Ramanujan's theories of
elliptic integrals to alternative bases.  These include but are not limited
to the ones proved by Berndt, Bhargava, and Garvan~\cite{Berndt95}, and
others.  In this section we derive the full set of such modular equations.
Most are shown in Table~\ref{tab:ramanujan}.  Our treatment in this section
is on the level of special function manipulations.  In~\S\ref{sec:final},
we shall give a direct modular interpretation of Ramanujan's theories.

His complete elliptic integral (of the first kind) in the theory of
signature~$r$ (or~base~$r$), where $r>1,$ is
$\K_r:(0,1)\to\mathbb{R}^+,$ a~monotone increasing function defined
by
\begin{subequations}
\begin{align}
\label{eq:altintrep0}
\K_r(\alpha_r) &= \tfrac\pi2 \,\, {}_2F_1\left(\tfrac1r,
1-\tfrac1r;\,1;\,\alpha_r\right)\\
&= \tfrac{\sin(\pi/r)}2 \int_0^1 x^{-1/r} (1-x)^{-1+1/r} (1-\alpha_r x)^{-1/r}\,dx.
\label{eq:altintrep}
\end{align}
\end{subequations}
Its limits as $\alpha_r\to0^+,1^-$ are $\pi/2,\infty,$ and in the latter
limit it diverges logarithmically.  The representation~(\ref{eq:altintrep})
comes from Euler's integral formula for~${}_2F_1$
\cite[\S2.1.3]{Erdelyi53}.  By the Schwartz--Christoffel theory,
$\K_r(\alpha_r)$ can be interpreted~\cite{Anderson2000} as the side
length of a certain parallelogram with angles $\frac\pi{r},\pi(1-\frac1r),$
its aspect ratio being such that it can be mapped conformally onto the
upper half plane, with its vertices taken to $0,1,\alpha_r^{-1},\infty$.

In Ramanujan's theory of signature~$r,$ an
$\alpha_r$\textendash\nobreak$\beta_r$ modular equation of degree~$N$ is an
explicit relation between $\alpha_r,\beta_r\in(0,1)$ induced by
\begin{equation}
  \label{eq:modulareqndef}
  \frac{\K_r'(\beta_r)}{\K_r(\beta_r)} = N\,
  \frac{\K_r'(\alpha_r)}{\K_r(\alpha_r)},
\end{equation}
where $\K_r'(\alpha_r):=\K_r(1-\alpha_r)$ is the so-called complementary
complete elliptic integral, and $N\in\mathbb{Q}^+$ is the degree.
A~degree-$N$ modular equation in signature~$r,$ interpreted more strongly,
should also include an expression for the `multiplier'
$\K_r(\alpha_r)/\K_r(\beta_r)$ from which (\ref{eq:modulareqndef})~may be
recovered.  If $N=N_1/N_2$ in lowest terms, a degree-$N$ modular equation
will optionally be referred to here as a modular equation of degree
$N_1:N_2$.  In all closed-form signature-$r$ modular equations found
to~date, $r$~is $2,3,4,$ or~$6,$ and the
$\alpha_r$\textendash\nobreak$\beta_r$ relation is algebraic.

To see that modular equations for the families of elliptic curves attached
to $\Gamma_0(M),$ $M=2,3,4,$ can be converted to modular equations in the
theories of signature $r=4,3,2,$ where $r=12/\psi(M),$ recall
from~\S\ref{sec:modular3} that when $M=2,3,4,$
\begin{equation}
    h_M(t_M) = {}_2F_1\bigl(\tfrac1{12}\psi(M),\tfrac1{12}\psi(M);\,1;\,-t_M/\kappa_M^{1/2}\bigr),
\end{equation}
with $\kappa_M=2^{12},3^6,2^8,$ respectively.  The function~$h_M$ has a
holomorphic continuation from a neighborhood of~$t_M=0$ to the half-line
$t_M>0$.  Let a bijection between $(0,\infty)\ni t_M$ and
$(0,1)\ni\alpha_r,$ where $r=12/\psi(M),$ be given by
\begin{subequations}
\begin{align}
\alpha_r &= \alpha_r(t_M) = t_M/(t_M+\kappa_M^{1/2}),\\
t_M &= t_M(\alpha_r) = \kappa_M^{1/2}\,\alpha_r/(1-\alpha_r).
\end{align}
\end{subequations}
Then one can write
\begin{subequations}
\begin{align}
\tfrac2\pi \K_r(\alpha_r) &= (1-\alpha_r)^{-\psi(M)/12}\,\, h_M(t_M(\alpha_r)),\\
h_M(t_M) &= (1+t_M/\kappa_M^{1/2})^{-\psi(M)/12}\,\, \tfrac2\pi \K_r(\alpha_r(t_M)).
\label{eq:hNKN}
\end{align}
\end{subequations}
These follow immediately from~(\ref{eq:Pfaff}), Pfaff's transformation
of~${}_2F_1$.

\begin{theorem}
  If\/ $M=2,3,4,$ on\/ $t_M>0$ the ratio\/ ${\rm
  i}\,\K_r'(\alpha_r(t_M))/\K_r(\alpha_r(t_M))$ equals\/
  $C_M\tau(t_M),$ where\/ $\tau(\cdot)$ is the\/ ${\rm i}\mathbb{R}$-valued
  branch of the period ratio\/ $\tau$ that was introduced
  in\/~{\rm\S}{\rm\ref{sec:modular}.}  In all three cases, the prefactor\/
  $C_M$ equals\/ $M^{1/2}$.
\end{theorem}
\begin{proof}
  Rewrite the representation for $\tau(t_M)$ given in
  Corollary~\ref{cor:nameless} in~terms of $\K_r(\alpha_r(t_M)),$
  with the aid of~(\ref{eq:hNKN}).
\end{proof}

\begin{remark*}
  By reviewing the proof of Corollary~\ref{cor:nameless}, one discovers
  that the simple formula $C_M=M^{1/2}$ for the prefactor comes ultimately
  from a subtle result (Eq.~(\ref{eq:nameless2})) on the asymptotic
  behavior of the Gauss hypergeometric function, together with evaluations
  of the Euler digamma function $\Psi(a)$ at the rational points
  $a=\frac12,\frac13,\frac14$ and $a=1-\frac12,1-\frac13,1-\frac14$.  This
  is a remarkably complicated proof of a seemingly simple result.
\end{remark*}

\begin{table}
\caption{Rationally parametrized modular equations in the theories of
signature $r=2,3,4$.  The parameter~$t$ signifies respectively
$t_{4N},t_{3N},t_{2N},$ i.e., the canonical Hauptmodul for
$\Gamma_0(4N),\Gamma_0(3N),\Gamma_0(2N)$.}
\begin{center}
{\small
\begin{tabular}{ccp{261pt}}
\hline\noalign{\smallskip}
$r$ & $N$ & Modular equation for~$\K_r,$ of degree~$N$\\
\noalign{\vskip2pt}\hline\noalign{\vskip1pt}\hline\noalign{\vskip1.5pt}
$2$ & $2$ & $\K_2\left(\frac{t(t+8)}{(t+4)^2}\right)=2\left[\frac{t+4}{t+8}\right]\,\K_2\left(\frac{t^2}{(t+8)^2}\right)$ \\
$2$ & $3$ & $\K_2\left(\frac{t(t+4)^3}{(t+2)^3(t+6)}\right)=3\left[\frac{t+2}{t+6}\right]\,\K_2\left(\frac{t^3(t+4)}{(t+2)(t+6)^3}\right)$ \\
$2$ & $4$ & $\K_2\left(\frac{t(t+8)}{(t+4)^2}\circ\,t(t+4)\right)=4\left[\frac{t+2}{t+4}\right]^{2}\,\K_2\left(\frac{t^2}{(t+8)^2}\circ\frac{t^2}{t+2}\right)$ \\
\noalign{\vskip 1pt}\hline\noalign{\vskip 1pt}
$3$ & $2$ & $\K_3\left(\frac{t(t+9)^2}{(t+6)^3}\right)=2\left[\frac{t+6}{t+12}\right] \,\K_3\left(\frac{t^2(t+9)}{(t+12)^3}\right)$ \\
$3$ & $3$ & $\K_3\left(\frac{t(t^2+9t+27)}{(t+3)^3}\right)=3\left[\frac{t+3}{t+9}\right]\,\K_3\left(\frac{t^3}{(t+9)^3}\right)$\\
$3$ & $4$ & $\K_3\left(\frac{t(t+9)^2}{(t+6)^3}\circ\,t(t+6)\right)=4\left[\frac{t^2+6t+6}{t^2+12t+24}\right]\,\K_3\left(\frac{t^2(t+9)}{(t+12)^3}\circ\frac{t^2}{t+2}\right)$ \\
$3$ & $6$ & $\K_3\left(\frac{t(t+9)^2}{(t+6)^3}\circ\,t(t^2+6t+12)=\frac{t(t^2+9t+27)}{(t+3)^3}\circ\frac{t(t+3)^2}{t+2}\right)$\\
    &     & $=6\left[\frac{t^3+6t^2+12t+6}{t^3+12t^2+36t+36}\right]\,\K_3\left(\frac{t^2(t+9)}{(t+12)^3}\circ\frac{t^3}{t^2+3t+3}=\frac{t^3}{(t+9)^3}\circ\frac{t^2(t+3)}{(t+2)^2}\right)$\\
\noalign{\vskip 1pt}\hline\noalign{\vskip 1pt}
$4$ & $2$ & $\K_4\left(\frac{t(t+16)}{(t+8)^2}\right)=2\left[\frac{t+8}{t+32}\right]^{1/2}\K_4\left(\frac{t^2}{(t+32)^2}\right)$\\
$4$ & $3$ & $\K_4\left(\frac{t(t+8)^3}{(t^2+12t+24)^2}\right)=3\left[\frac{t^2+12t+24}{t^2+36t+216}\right]^{1/2}\K_4\left(\frac{t^3(t+8)}{(t^2+36t+216)^2}\right)$ \\
$4$ & $4$ & $\K_4\left(\frac{t(t+16)}{(t+8)^2}\circ\,t(t+8)\right)=4\left[\frac{t^2+8t+8}{t^2+32t+128}\right]^{1/2}\K_4\left(\frac{t^2}{(t+32)^2}\circ\frac{t^2}{t+4}\right)$\\
$4$ & $5$ & $\K_4\left(\frac{t(t+4)^5}{(t^2+6t+4)^2(t^2+8t+20)}\right)$\\
    &     & $=5\,\left[\frac{t^2+6t+4}{t^2+30t+100}\right]^{1/2}\K_4\left(\frac{t^5(t+4)}{(t^2+8t+20)(t^2+30t+100)^2}\right)$ \\
$4$ & $6$ & $\K_4\left(\frac{t(t+16)}{(t+8)^2}\circ\frac{t(t+4)^3}{t+3}=\frac{t(t+8)^3}{(t^2+12t+24)^2}\circ\,t(t+6)\right)$\\
    &     & $=6\,\left[\frac{t^4+12t^3+48t^2+72t+24}{t^4+36t^3+288t^2+864t+864}\right]^{1/2}$\\
    &     & \hfill${}\times\K_4\left(\frac{t^2}{(t+32)^2}\circ\frac{t^3(t+4)}{(t+3)^3}=\frac{t^3(t+8)}{(t^2+36t+216)^2}\circ\frac{t^2}{t+2}\right)$\\
$4$ & $8$ & $\K_4\left(\frac{t(t+16)}{(t+8)^2}\circ\,t(t+8)\,\circ\,t(t+4)\right)$\\
    &     & $=8\left[\frac{t^4+8t^3+24t^2+32t+8}{t^4+32t^3+192t^2+512t+512}\right]^{1/2}\K_4\left(\frac{t^2}{(t+32)^2}\circ\frac{t^2}{t+4}\circ\frac{t^2}{t+2}\right)$\\
$4$ & $9$ & $\K_4\left(\frac{t(t+8)^3}{(t^2+12t+24)^2}\circ\,t(t^2+6t+12)\right) $\\
    &     & $=9\,\left[\frac{t^6+12t^5+60t^4+156t^3+216t^2+144t+24}{t^6+36t^5+324t^4+1404t^3+3240t^2+3888t+1944}\right]^{1/2}$\\
    &     & \hfill${}\times\K_4\left(\frac{t^3(t+8)}{(t^2+36t+216)^2}\circ\frac{t^3}{t^2+3t+3}\right)$\\
\noalign{\vskip2pt}\hline
\end{tabular}
}%
\end{center}
\label{tab:ramanujan}
\end{table}

\begin{corollary}
  In the cases\/ $M=2,3,4,$ any degree\/ $N_1:N_2$ modular equation for the
  function\/ $h_M=h_M(t_N)$ associated to the group\/ $\Gamma_0(N)$ may be
  converted to a degree\/ $N_1:N_2$ modular equation for Ramanujan's
  complete elliptic integral of the first kind,
  $\K_{r}=\K_{r}(\alpha_{r}),$ with $r=12/\psi(M)$.
\end{corollary}

\begin{theorem}
  In the theories of signature\/ $r=2,3,4,$ one has the\/ $N:1$ modular
  equations of Table\/~{\rm\ref{tab:ramanujan}.}  Each holds for all\/
  $t>0$ and incorporates a rationally parametrized\/
  $\alpha_r$\textendash\nobreak$\beta_r$ modular equation of degree\/~$N,$
  the pair\/ $\alpha_r,\beta_r$ being the arguments of\/ $\K_r$ on
  the left and right sides.  In addition, $\K_3$
  and\/ $\K_4$ satisfy

  {\small
  \begin{align}
    & \K_3\left(\tfrac{t(t^2+9t+27)}{(t+3)^3}\circ\tfrac{t^2(t+3)}{(t+2)^2}=\tfrac{t^2(t+9)}{(t+12)^3}\circ\, t(t^2+6t+12)\right)\nonumber\\
    & \qquad = \left(\tfrac32\right) \left[\tfrac{t^3+6t^2+12t+12}{t^3+6t^2+18t+18}\right]\! \K_3\left(\tfrac{t^3}{(t+9)^3}\circ\tfrac{t(t+3)^2}{t+2}=\tfrac{t(t+9)^2}{(t+6)^3}\circ\tfrac{t^3}{t^2+3t+3}\right),\nonumber
    \\
    & \K_4\left(\tfrac{t^2}{(t+32)^2}\circ\tfrac{t(t+4)^3}{t+3} = \tfrac{t(t+8)^3}{(t^2+12t+24)^2}\circ\tfrac{t^2}{t+2}\right)\nonumber\\
    & \qquad =\left(\tfrac32\right) \left[\tfrac{t^4+12t^3+48t^2+96t+96}{t^4+12t^3+72t^2+216t+216}\right]^{1/2}\nonumber\\
    &\qquad\qquad\qquad\qquad\quad{}\times \K_4\left(\tfrac{t(t+16)}{(t+8)^2}\circ\tfrac{t^3(t+4)}{(t+3)^3} = \tfrac{t^3(t+8)}{(t^2+36t+216)^2}\circ\,t(t+6)\right),\nonumber
  \end{align}
  }%
  which are modular equations of degree\/ $3:2,$ i.e., of
  degree\/~{\rm$\tfrac32$.}
\label{thm:ramafull}
\end{theorem}
\begin{proof}
  Rewrite the modular equations given in Theorems \ref{thm:maintable}
  and~\ref{thm:supplement} for~$h_M,$ $M=2,3,4,$ in~terms of
  $\K_r,$ $r=4,3,2,$ with the aid of~(\ref{eq:hNKN}).
\end{proof}

One can also derive modular equations of a certain mixed type from the
elliptic-family modular equations of~\S\ref{sec:families}.  These are
$\alpha_r$\textendash\nobreak$\beta_s$ relations, induced by
\begin{equation}
  \frac{\K_s'(\beta_r)}{\K_s(\beta_r)} = N\,
  \frac{\K_r'(\alpha_r)}{\K_r(\alpha_r)},
\end{equation}
where $r\neq s$.  If the degree $N\in\mathbb{Q}^+$ equals $N_1/N_2$ in
lowest terms, such an algebraic relation may be called an
$\alpha_r$\textendash\nobreak$\beta_s$ modular equation of degree
${N_1:N_2}$\,; or if accompanied by an explicit multiplier, a modular
transformation of degree ${N_1:N_2}$ from $\K_r$ to~$\K_s$.
The source of such mixed modular equations, when $r,s\in\{2,3,4\},$ is {\em
commensurability\/}: the fact that the intersection of any two of the
subgroups $\Gamma_0(4),\Gamma_0(3),\Gamma_0(2)$ has finite index in both.

\begin{theorem}
\label{thm:fullcollection}
  The full collection of rationally parametrized\/
  $\alpha_r$\textendash\nobreak$\beta_s$ modular equations, for\/ $r\neq s$
  with $r,s\in\{2,3,4\},$ includes
  \renewcommand\labelitemi{$\circ\ $}
  \begin{itemize}
    \item $\alpha_2$\textendash\nobreak$\beta_3$ modular equations of degrees\/ $4,2,\frac43,1,\frac23,\frac13${\rm;}
    \item $\alpha_2$\textendash\nobreak$\beta_4$ modular equations of degrees\/
    $8,6,4,3,2,1,\frac23,\frac12,\frac13,\frac14${\rm;}
    \item $\alpha_3$\textendash\nobreak$\beta_4$ modular equations of degrees\/ $9,\frac92,6,3,\frac32,2,1,\frac34,\frac12,\frac13,\frac14,\frac16${\rm.}
  \end{itemize}
  Each such\/ $\alpha_r$\textendash\nobreak$\beta_s$ relation is the basis of
  a rationally parametrized modular transformation from\/ $\K_r$
  to\/~$\K_s$.
\end{theorem}
\begin{proof}
  For any~$\mathcal{N}$ with $\Gamma_0(\mathcal{N})$ of genus zero,
  $t_M(d\tau)$ may be expressed rationally in~terms
  of~$t_{\mathcal{N}}(\tau)$ if $M\divides{\mathcal{N}}$
  and~$d\divides(\mathcal{N}/M)$.  So if $M_1\divides\mathcal{N}$ and
  $d_i\divides(\mathcal{N}/M_i),$ $i=1,2,$ there is an
  $t_{M_1}(d_1\tau)$\textendash\nobreak$t_{M_2}(d_2\tau)$ modular equation,
  i.e., an $t_{M_1}$\textendash\nobreak$t_{M_2}$ modular equation of
  degree~$d_2/d_1$.  By Theorem~\ref{thm:extendedkey} (part~2), it yields a
  modular transformation from $h_{M_1}$ to~$h_{M_2}$.  But the
  Hauptmodul~$t_{M_i}$ and function $h_{M_i}=h_{M_i}(t_{M_i}),$ when
  $M_i=2,3,4,$ correspond to the invariant~$\alpha_r$ and elliptic integral
  $\K_r=\K_r(\alpha_r),$ where $r=4,3,2$.  The given lists
  are computed by considering all relevant~$\mathcal{N}$
  (i.e.,~$\mathcal{N}=4,6,\allowbreak8,\allowbreak12,\allowbreak16,18$),
  and enumerating divisors~$M_i\in\{2,3,4\},$ and $d_i$~such that
  $d_i\divides(\mathcal{N}/M_i)$.
\end{proof}
\begin{remark*}
  The collection mentioned in the theorem consists of 28~modular equations,
  but at most~14 are independent.  As Berndt et~al.\ have observed, for any
  $r\neq s$ there is an involution of the set of
  $\alpha_r$\textendash\nobreak$\beta_s$ modular equations, called
  `reciprocation.'  It consists in applying the simultaneous substitutions
  $\alpha_r\mapsto1-\alpha_r,$ $\beta_s\mapsto1-\beta_s$.  If
  $r=12/\psi(M_1)$ and~$s=12/\psi(M_2)$ with $M_1,M_2\in\{2,3,4\},$ it is
  not difficult to see that the procedure transforms an
  $\alpha_r$\textendash\nobreak$\beta_s$ modular equation of degree~$N$ to
  one of degree $(M_1/M_2)/N$.  For the three types of modular equation
  mentioned in the theorem, the degree involutions are $N\mapsto\frac43/N,$
  $N\mapsto2/N,$ and $N\mapsto\frac32/N$.
\end{remark*}

The extent to which the modular equations derived in this section are new
can now be discussed.  To begin with, the modular equations for~$\K_2$ in
Table~\ref{tab:ramanujan}, of degrees $N=2,3,4,$ are classical.  This is
because $r=2$ is the classical base, and $\K_2=\K_2(\alpha_2)$ is identical
to $\K=\K(\alpha),$ the classical complete elliptic integral, which is a
single-valued function of the $\alpha$-invariant in a neighborhood of the
infinite cusp, at~which $\alpha=0$. The argument~$\alpha_2$ of~$\K_2$ is
now seen to have a modular interpretation: being identical to~$\alpha,$ it
is a Hauptmodul for~$\Gamma_0(4)$.  The three underlying
$\alpha_2$\textendash\nobreak$\beta_2$ modular equations, i.e.,
$\alpha$\textendash\nobreak$\beta$ modular equations, were discussed
in~\S\ref{sec:parametrized2}.  The $N=2$ equation is Landen's
transformation, in the form already given in Eq.~(\ref{eq:duo}).

Of the remaining rationally parametrized equations in the table, Ramanujan
found the degree-$2,3,4$ modular equations for $\K_3$
and~$\K_4$.  That~is, he found the underlying
$\alpha_3$\textendash\nobreak$\beta_3$ and
$\alpha_4$\textendash\nobreak$\beta_4$ modular equations, and expressions
for the multipliers.  He also found the degree-$\frac32$ modular equation
for~$\K_3,$ given in Theorem~\ref{thm:ramafull}, and several of the
transformations mentioned in Theorem~\ref{thm:fullcollection}: those from
$\K_2$ to~$\K_3$ of degrees~$1,\frac23,$ those from
$\K_2$ to~$\K_4$ of degrees~$1,2,$ and that from
$\K_3$ to~$\K_4$ of degree~$1$.  Presumably he was aware
of the reciprocation principle.

Proofs of his results have been constructed by the Borweins, Berndt,
Garvan, and others.  The Borweins derived the transformations of degree~$1$
among
$\K_2,\K_3,\allowbreak\K_4$~\cite[\S5.5]{Borwein87},
and the degree-$3$ modular equation in the theory of
signature~$3$~\cite{Borwein91}.  Berndt, Bhargava and
Garvan~\cite{Berndt95} (see also Berndt~\cite[Chap.~33]{BerndtV})
systematically derived Ramanujan's remaining modular equations and
transformations in the theories of signature $3$ and~$4$.
Garvan~\cite[Eq.~(2.34)]{Garvan95} additionally derived the transformation
from $\K_2$ to~$\K_3$ of degree~$\frac43$.

Recently, Berndt, Chan and Liaw~\cite{Berndt2001} have conducted further
investigations into the theory of signature~$4$.  However, the
signature-$4$ modular equations in Table~\ref{tab:ramanujan} of degrees
greater than~$4$ are new, as is the signature-$4$ degree-$\frac32$
equation; and also the signature-$3$ degree-$6$ equation.  Of~the $14=28/2$
independent rationally parametrized transformations of $\K_r$
to~$\K_s$ identified in Theorem~\ref{thm:fullcollection}, all of
which can readily be worked~out explicitly if needed, $10$~are new.

\section{A Modular Approach to Signatures $2,3,4$ (and~$6$)}
\label{sec:final}

The derivation in~\S\ref{sec:ramanujan} of modular equations in Ramanujan's
theories of signature $r=2,3,4$ relied on Pfaff's transformation
of~${}_2F_1,$ and took place on the level of special functions.  In this
section we give a direct modular interpretation of the complete elliptic
integral~$\K_r$ and its argument~$\alpha_r$ (the
`$\alpha_r$-invariant').  Modular interpretations were pioneered by the
Borweins~\cite{Borwein87,Borwein91}, but our new interpretation is very
concise, and lends itself to extension.  Like the classical complete
elliptic integral, each $\K_r(\alpha_r)$ is simply a weight-$1$
modular form, expressed as a function of a Hauptmodul.  We~close by
discussing the underdeveloped, and very interesting, theory of
signature~$6$.  It does not fit into the Picard--Fuchs framework of this
article, but we indicate a slightly more general Gauss--Manin framework
into which it fits.

When $M=2,3,4$ and correspondingly $r=12/\psi(M)=4,3,2,$ consider the
subgroup $\Gamma_0(M)<\Gamma(1)$.  The function field
of~$X_0(M)=\Gamma_0(M)\setminus\mathcal{H}^*$ is generated by the canonical
Hauptmodul 
\begin{subequations}
\begin{align}
&t_M=\kappa_M\cdot[M]^{24/(M-1)}/\,[1]^{24/(M-1)},\\
&\kappa_M:=M^{12/(M-1)},
\end{align}
\end{subequations}
and there are exactly three fixed points on~$X_0(M)$.  (See
Table~\ref{tab:hauptmoduln}.)  These are the cusp $t_M=0$ (i.e.,
$\tau\in\bigl[\frac1M\bigr]_M,$ including the infinite cusp $\tau={\rm
i}\infty$), the cusp $t_M=\infty$ (i.e., $\tau\in\bigl[\frac11\bigr]_M,$
including $\tau=0$), and the third fixed point $t_M=-\kappa_M^{1/2}$.  When
$M=2,3,4,$ the third point is respectively a quadratic elliptic point, a
cubic one, and a cusp (namely, $\tau\in\bigl[\frac12\bigr]_4$).  To focus
on this third fixed point, which is stabilized by the Fricke involution
$t_M\mapsto \kappa_M/t_M,$ one defines an alternative Hauptmodul
\begin{equation*}
\alpha_r:=t_M/(t_M+\kappa_M^{1/2})
\end{equation*}
that is zero at the infinite cusp, like the canonical Hauptmodul~$t_M,$ but
has its pole at the third point rather than at the cusp $\tau=0$.  At the
latter cusp, $\alpha_r=1$.  In~terms of~$\alpha_r,$ the Fricke involution
is the map $\alpha_r\mapsto1-\alpha_r$.

For each signature $r=2,3,4,$ associated to $M=4,3,2$ by $r=12/\psi(M),$
define also a triple of weight-$1$ modular forms
$\mathcal{A}_r,\mathcal{B}_r,\mathcal{C}_r$ for~$\Gamma_0(M)$ thus:
\begin{align*}
  \mathcal{A}_r&:=(1-\alpha_r)^{-1/r}\mathcal{B}_r, &
  \mathcal{B}_r&:=\mathfrak{h}_M, &
  \mathcal{C}_r&:=\left[\tfrac{1-\alpha_r}{\alpha_r}\right]^{-1/r}\mathcal{B}_r,
\end{align*}
where $\mathfrak{h}_M=h_M\circ t_M$ is the canonical weight-$1$ modular
form.  By Corollary~\ref{cor:hNformula1}, $\mathcal{B}_r$~has divisor
$\frac{\psi(M)}{12}(t_M=\infty) = \frac1r(\alpha_r=1)$ on~$X_0(M)$.  As a
consequence of their definition, $\mathcal{C}_r$~has divisor
$\frac{\psi(M)}{12}(t_M=0) = \frac1r(\alpha_r=0)$ and $\mathcal{A}_r$~has
divisor $\frac{\psi(M)}{12}(t_M=-\kappa_M^{1/2}) =
\frac1r(\alpha_r=\infty)$.  Also, the triple satisfies (by definition)
\begin{equation}
{\mathcal{A}_r}^r = {\mathcal{B}_r}^r + {\mathcal{C}_r}^r,
\end{equation}
which in each case is an equality between two weight-$r$ modular forms.

\begin{table}
\begin{center}
  \caption{For $M=2,3,4,$ the Hauptmodul~$t_M$ and weight-$1$ modular
  form~$\mathfrak{h}_M$ for~$\Gamma_0(M)$; and in the corresponding theory
  of signature $r=12/\psi(M)=4,3,2,$ the Hauptmodul~$\alpha_r$ and modular
  forms $\mathcal{A}_r,\mathcal{B}_r,\mathcal{C}_r,$ which satisfy
  ${\mathcal{A}_r}^{\!r}={\mathcal{B}_r}^{\!r}+{\mathcal{C}_r}^{\!r}$.
  Here $\hat{\K}_r$~signifies $(2/\pi)\K_r,$ and $\sum$ means
  $\sum_{n=1}^\infty$.}  {\small
  \begin{tabular}{cc|l}
    \hline\noalign{\smallskip}
    $M$ & $r$ & Hauptmoduln and modular forms\\
    \noalign{\smallskip}\hline\noalign{\vskip1pt}\hline\noalign{\vskip1.5pt}
    $2$ & $4$ & $t_2=2^{12}\cdot[2]^{24}/\,[1]^{24}$\\
        &     & $\mathfrak{h}_2=h_2\circ t_2=[1]^4/\,[2]^2=1+4\sum E_1(n;4)(-q)^n$\\
        &     & \\
        &     & $\alpha_4\, = 2^6\,[2]^{24}/\bigl(2^6\,[2]^{24}+[1]^{24}\bigr)$\\
        &     & $1-\alpha_4\, = [1]^{24}/\bigl(2^6\,[2]^{24}+[1]^{24}\bigr)$\\
        &     & $\mathcal{A}_4 = \hat{\K}_4\circ\alpha_4= (2^6\,[2]^{24}+[1]^{24})^{1/4}/\,[1]^2[2]^2$\\
        &     & $\hphantom{\mathcal{A}_4 } =
    \sqrt{1+24\sum\left(\sum_{d\divides{n},\,\textrm{$d$ odd}}d\right)q^n}$\\
        &     & $\mathcal{B}_4 = [1]^4/\,[2]^2= 1+4\sum E_1(n;4)(-q)^n$\\
        &     & $\mathcal{C}_4\, = 2^{3/2}\cdot[2]^4/\,[1]^2 = 2^{3/2}q^{1/4}\left[1+\sum E_1(4n+1;4)\,q^n\right]$\\
    \noalign{\vskip 1pt}\hline\noalign{\vskip 1pt}
    $3$ & $3$ & $t_3=3^6\cdot[3]^{12}/\,[1]^{12}$\\
        &     & $\mathfrak{h}_3=h_3\circ t_3=[1]^3/\,[3]=1-3\sum \left[E_1(n;3)-3E_1(n/3;3)\right]q^n$\\
        &     & \\
        &     & $\alpha_3\, = 3^3\,[3]^{12}/\bigl(3^3\,[3]^{12}+[1]^{12}\bigr)$\\
        &     & $1-\alpha_3\, = [1]^{12}/\bigl(3^3\,[3]^{12}+[1]^{12}\bigr)$\\
        &     & $\mathcal{A}_3 =\hat{\K}_3\circ\alpha_3 = (3^3\,[3]^{12}+[1]^{12})^{1/3}/\,[1][3]$\\
        &     & $\hphantom{\mathcal{A}_3 } = 1+6\sum E_1(n;3)\,q^n$\\
        &     & $\mathcal{B}_3 = [1]^3/\,[3] = 1-3\sum\left[E_1(n;3)-3E_1(n/3;3)\right]q^n$\\
        &     & $\mathcal{C}_3\, = 3\cdot [3]^3/\,[1] = 3\,q^{1/3}\left[1+\sum E_1(3n+1;3)\,q^n\right]$\\
    \noalign{\vskip 1pt}\hline\noalign{\vskip 1pt}
    $4$ & $2$ & $t_4=2^8\cdot[4]^8/\,[1]^8$\\
        &     & $\mathfrak{h}_4=h_4\circ t_4=[1]^4/\,[2]^2=1+4\sum E_1(n;4)(-q)^n$\\
        &     & \\
        &     & $\alpha_2\, = 2^4\,[4]^{8}/\bigl(2^4\,[4]^{8}+[1]^{8}\bigr) = 2^4\cdot [1]^8[4]^{16}/\,[2]^{24}$\\
        &     & $1-\alpha_2\, = [1]^{8}/\bigl(2^4\,[4]^{8}+[1]^{8}\bigr) =  [1]^{16}[4]^{8}/\,[2]^{24}$\\
        &     & $\mathcal{A}_2 = \hat{\K}_2\circ\alpha_2 = (2^4\,[4]^8+[1]^8)^{1/2}/\,[2]^2 = [2]^{10}/\,[1]^4[4]^4$\\
        &     & $\hphantom{\mathcal{A}_2 }= 1+4\sum E_1(n;4)\,q^n$\\
        &     & $\mathcal{B}_2=[1]^4/\,[2]^2=1+4\sum E_1(n;4)(-q)^n$\\
        &     & $\mathcal{C}_2\, = 2^2\cdot[4]^4/\,[2]^2 = 4\,q^{1/2}\left[1+\sum E_1(2n+1;4)\,q^n\right]$\\
    \noalign{\vskip2pt}\hline
  \end{tabular}
  }%
\end{center}
\label{tab:final}
\end{table}

\begin{theorem}
For\/ $r=2,3,4,$ corresponding to\/ $M=4,3,2$ by\/ $r=12/\psi(M),$ the
Hauptmodul\/ $\alpha_r$ and triple of weight\/-$1$ modular forms\/
$\mathcal{A}_r,\mathcal{B}_r,\mathcal{C}_r$ for\/~$\Gamma_0(M)$ have the
eta product representations and\/ $q$-expansions shown in
Table\/~{\rm\ref{tab:final}.}
\end{theorem}
\begin{proof}
The formulas for $t_M$ and $\mathcal{B}_r=\mathfrak{h}_M=h_M\circ t_M$ are
taken from Table~\ref{tab:xNhN}, and those for~$\mathcal{C}_r$ from
Table~\ref{tab:insert}, since $\mathcal{C}_r=M^{\psi(M)/2(M-1)}\,\bar
h_M\circ t_M,$ in that table's notation.  Since $\mathcal{A}_2$~equals
$\left[(t_4+16)/16\right]^{1/2}(h_4\circ t_4),$ it~is an alternative
weight\nobreakdash-$1$ modular form of the sort considered
in~\S\ref{sec:modular2}, with its (single) zero located at a finite cusp
other than $\tau=0$; so its eta product and $q$-series can be found in
Table~\ref{tab:windup}.  The $q$-series representations for $\mathcal{A}_3$
and~${\mathcal{A}_4}^2$ were discovered by Ramanujan.
\end{proof}
\begin{remark*}
  As is clear from its definition as $t_4/(t_4+16)$ and its eta product
  representation, the $\alpha_2$-invariant is the $\alpha$-invariant, an
  alternative Hauptmodul for~$\Gamma_0(4)$.  Ramanujan's theory of
  signature~$2$ is Jacobi's classical theory of elliptic integrals.  The
  $r=2$ identity ${\mathcal{A}_2}^2 = {\mathcal{B}_2}^2 +
  {\mathcal{C}_2}^2$ is, in~fact, the Jacobi theta identity
  ${\vartheta_3}(0)^4 = {\vartheta_4}(0)^4 + {\vartheta_2}(0)^4,$ as is
  evident from the $q$-series for
  $\mathcal{A}_2,\mathcal{B}_2,\mathcal{C}_2$.
\end{remark*}
\begin{remark*}
  It also follows by examining the $q$-series for
  $\mathcal{A}_r,\mathcal{B}_r,\mathcal{C}_r,$ that the identities
  ${\mathcal{A}_r}^r = {\mathcal{B}_r}^r + {\mathcal{C}_r}^r,$ $r=3,4,$ are
  explicitly modular restatements of $q$-series identities previously
  obtained by the Borweins~\cite{Borwein91}.
\end{remark*}
\begin{remark*}
  The modular form $\mathcal{A}_4$ for~$\Gamma_0(2)$ is anomalous:
  no~number-theoretic interpretation of the coefficients of its $q$-series
  $1+12[q-5q^2+64q^3\allowbreak{-917q^4}+\cdots]$ is known, unlike that
  of~${\mathcal{A}_4}^2,$ which is the theta series of the $D_4$ lattice
  packing.  Its coefficient sequence
  $1,12,\allowbreak-60,\allowbreak768,-11004,\dots$ is Sloane's {\tt
  A108096}.  For some apposite remarks on how roots of $q$-series with
  integer coefficients may unexpectedly turn~out to have the same property,
  see Heninger et~al.~\cite{Heninger2006}.
\end{remark*}

By applying Pfaff's transformation of~${}_2F_1,$ one can immediately
transform the known fact (see~\S\ref{sec:modular3}) that for $M=4,3,2,$
\begin{equation}
  h_M(t_M) = {}_2F_1\bigl(\tfrac1{12}\psi(M),\tfrac1{12}\psi(M);\,1;\,-t_M/\kappa_M^{1/2}\bigr),
\end{equation}
into the statement that $\mathcal{A}_r={\hat\K}_r\circ\alpha_r$ when $r=12/\psi(M)=2,3,4,$ with
\begin{equation}
{\hat\K}_r(\alpha_r) := \tfrac2{\pi}\K_r(\alpha_r) = {}_2F_1\left(\tfrac1r,
1-\tfrac1r;\,1;\,\alpha_r\right).
\end{equation}
That is, Ramanujan's complete elliptic integral~$\K_r(\alpha_r),$
$r=2,3,4,$ arises naturally in a modular context, when the noncanonical
weight-$1$ modular form~$\mathcal{A}_r$ is expressed as a function of the
alternative Hauptmodul~$\alpha_r$.  This has not previously been realized.
Ramanujan may, in~fact, have made a good choice in focusing
on~$\mathcal{A}_r,$ rather than on its canonical
counterpart~$\mathcal{B}_r,$ since there are ways in which
$\mathcal{A}_r$~is `nicer' than~$\mathcal{B}_r$.  It can be shown that the
multiplier systems of~$\mathcal{A}_2,\mathcal{A}_3$ are given by Dirichlet
characters (of conductors $4$ and~$3$ respectively), unlike those
of~$\mathcal{B}_2,\mathcal{B}_3$; cf.~Table~\ref{tab:xNhN}.  Also, it can
be shown that the squares ${\mathcal{A}_3}^2,{\mathcal{A}_4}^2$ are the
unique (normalized) weight\nobreakdash-$2$ modular forms with trivial
character for $\Gamma_0(3),\Gamma_0(2),$ respectively.

The preceding interpretation of Ramanujan's theories to alternative bases
does not cover his theory of signature~$6$.  But on a purely formal level,
one can develop that theory in parallel with those of signatures $2,3,4$.
His complete elliptic integral $\K_6=\K_6(\alpha_6)$ is
defined by~(\ref{eq:altintrep0}), i.e., by
\begin{equation}
  \K_6(\alpha_6) = \tfrac{\pi}2\,\,{}_2F_1\left(\tfrac16,\tfrac56;\,1;\,\alpha_6\right).
\end{equation}
To recover it, one introduces a `Hauptmodul'~$y_6,$ defined implicitly by
\begin{equation}
\label{eq:formalH}
  j=\frac{(y_6+432)^2}{y_6}=12^3 + \frac{(y_6-432)^2}{y_6}.
\end{equation}
By mechanically following the prescription of Definition~\ref{def:1}, one
defines an associated function $H_6=H_6(y_6)$ by
\begin{equation}
  H_6(y_6) = \left[432^{-1}(y_6+432)\right]^{-1/6}\,{}_2F_1\left(\tfrac1{12},\tfrac5{12};\,1;\,\tfrac{12^3y_6}{(y_6+432)^2}\right).
\end{equation}
One can show that $H_6$~satisfies a normal-form `Picard--Fuchs equation' of
hypergeometric type, with three singular points ($y_6=0$ and~$\infty,$
which by~(\ref{eq:formalH}) are formally over the cusp $j=\infty$
on~$X(1)$; and~$y_6=-432,$ over the cubic elliptic point $j=0$).  This
leads to the alternative hypergeometric representation
\begin{equation}
  \label{eq:altrep6}
  H_6(y_6) = {}_2F_1\left(\tfrac16,\tfrac16;\,1;\,-y_6/432\right),
\end{equation}
which closely resembles the formulas of Table~\ref{tab:specfuncs}.  As
in~\S\ref{sec:ramanujan} one introduces a `noncanonical
Hauptmodul'~$\alpha_6,$ related to~$y_6$ by
\begin{subequations}
\begin{align}
\alpha_6 &= \alpha_6(y_6) = y_6/(y_6+432),\\
y_6&= y_6(\alpha_6) = 432\,\alpha_6/(1-\alpha_6).
\end{align}
\end{subequations}
Then one readily deduces
\begin{subequations}
\begin{align}
\tfrac2{\pi}\K_6(\alpha_6) &=(1-\alpha_6)^{-1/6}\,H_6(y_6(\alpha_6)),\\
H_6(y_6)&= (1+y_6/432)^{-1/6}\, \tfrac2{\pi}\K_6(\alpha_6(y_6)),
\end{align}
\end{subequations}
from Pfaff's transformation of~${}_2F_1$.  In this way the
function~$\K_6$ is recovered.

One can go further, and work out modular transformations of $H_6$
and~$\K_6$.  Starting with~(\ref{eq:modular1}), the degree-$2$
classical modular equation $\Phi_2(j,j')=0,$ by a lengthy process of
polynomial elimination one first derives the $y_6$\textendash\nobreak$y_6'$
modular equation of degree~$2$.  This can be rationally parametrized with
the aid of the {\sc Maple} {\tt algcurves} package, or other software.
A~further analysis yields the degree-$2$ modular transformation of~$H_6$
that incorporates the parametrized $y_6$\textendash\nobreak $y_6'$ modular
equation, namely
\begin{multline}
\label{eq:rama6}
  H_6\left(
  \tfrac{s(s+60)^2(s+72)^2(s+96)}{(s+48)(s+80)(s+120)^2}
  \right) \\= 2\left[
  \tfrac{(s+60)(s+80)(s+96)^2}{(s+48)(s+120)^2}
    \right]^{-1/6}
  H_6\left(
  \tfrac{s^2(s+48)^2(s+72)(s+120)}{(s+60)(s+80)^2(s+96)^2}
  \right).
\end{multline}
Expressing $H_6$ in~terms of~$\K_6,$ and replacing $s$ by $12s,$
then yields
\begin{equation}
  \K_6\left(
  \tfrac{s(s+5)^2(s+6)^2(s+8)}{(s^2+10s+20)^3}
  \right) = 2\left[
  \tfrac{s^2+10s+20}{s^2+20s+80}
    \right]^{1/2}
  \K_6\left(
  \tfrac{s^2(s+4)^2(s+6)(s+10)}{(s^2+20s+80)^3}
  \right),
  \label{eq:intriguing}
\end{equation}
which incorporates a rational parametrization of the
$\alpha_6$\textendash\nobreak$\beta_6$ modular equation of degree~$2$.
Here $\alpha_6,\beta_6$ are of~course the arguments of the left and
right~$\K_6$'s.  This uniformization of the
$\alpha_6$\textendash\nobreak$\beta_6$ relation agrees with the explicit
expression for the algebraic function $\beta_6=\beta_6(\alpha_6)$ obtained
by Borwein, Borwein, and Garvan~\cite{Borwein93}, using radicals.  The
intriguing identity~(\ref{eq:intriguing}) is a parametrized signature-$6$
counterpart of Landen's transformation, Eq.~(\ref{eq:one}).  Ramanujan did
not discover~it, but he would surely have appreciated~it.

The just-sketched manipulations are however quite formal, valid only near
the infinite cusp $\tau={\rm i}\infty,$ where $y_6=0$ and~$\alpha_6=0$.
The `Hauptmodul'~$y_6,$ algebraic over~$j,$ does not extend to a
single-valued function on the upper half $\tau$-plane, as is clear from its
definition~(\ref{eq:formalH}).  Also, the Gauss hypergeometric equation
satisfied by $H_6=H_6(y_6)$ cannot be the normal-form Picard--Fuchs
equation attached to any first-kind Fuchsian subgroup $\Gamma<{\it
PSL}(2,\mathbb{Z})$.  By~(\ref{eq:altrep6}), the exponent differences at
its singular points $y_6=0,\infty,-432$ are $0,0,\frac23$ respectively; and
by Theorem~\ref{thm:ford}, a difference of~$\frac23$ should not be present.

A recasting of the theory of signature~$6$ in algebraic--geometric terms
should be based not on the elliptic modular surface
$\mathfrak{E}_\Gamma\to\Gamma\setminus\mathcal{H}^*$ attached to
some~$\Gamma,$ but rather on an elliptic {\em non-modular\/} surface
$\mathfrak{E}\stackrel{\pi}{\to}\mathbb{P}^1(\mathbb{C})_{y_6},$ the base
curve of which is parametrized by the pseudo-Hauptmodul~$y_6$.  Its
functional invariant $j=j(y_6),$ in Kodaira's terminology, is given
by~(\ref{eq:formalH}).  Each of its fibres is an elliptic curve
over~$\mathbb{C},$ except for those above the singular points
$y_6=0,\infty,-432$.  This elliptic surface is one of the rational ones,
with nonconstant $j$-invariant, three singular fibres, and a section, that
were classified by Schmickler-Hirzebruch~\cite{SchmicklerHirzebruch85}.  It
is her `Fall~$8$', and up~to a broad notion of equivalence, it is the {\em
only\/} such surface that does not come from $\mathcal{H}^*$ and
some~$\Gamma$ in the classical way.  Its singular fibres are of Kodaira
types $I_1,I_1^*,{\it IV}$.  The hypergeometric equation satisfied by
$H_6=H_6(y_6)$ is its normal-form `Picard--Fuchs equation', though the term
is not really appropriate, suggesting as it does close ties to the
classical theory of automorphic functions.

The parameter~$s$ in the degree-$2$ modular equation for~$H_6,$
Eq.~(\ref{eq:rama6}), is best viewed as the coordinate on the base of a
second rational elliptic surface,
$\tilde{\mathfrak{E}}\stackrel{\tilde\pi}{\to}\mathbb{P}^1(\mathbb{C})_s,$
which parametrizes the modular equation.  There is a commutative diagram
\begin{equation}
  \begin{CD}
  \tilde{\mathfrak{E}} @>>> \mathfrak{E} @>>> \hat{\mathfrak{E}}_{1} \\
  @V{\tilde\pi}VV	  @V{\pi}VV	  @V{{\pi}_1}VV \\
  \mathbb{P}^1(\mathbb{C})_{s} @>>> \mathbb{P}^1(\mathbb{C})_{y_6} @>>> X(1)
  \end{CD}
\end{equation}
in which the lower map $s\mapsto y_6$ is the degree-$6$ rational function
of~$s$ appearing on the left-hand side of~(\ref{eq:rama6}).  By pulling
back the `Picard--Fuchs equation' for
$\mathfrak{E}\stackrel{\pi}{\to}\mathbb{P}^1(\mathbb{C})_{y_6}$ along
$\mathbb{P}^1(\mathbb{C})_{s}\to\mathbb{P}^1(\mathbb{C})_{y_6},$ one
readily computes a (formal!)\ Picard--Fuchs equation
for~$\tilde{\mathfrak{E}}$.  It turns~out to have $10$~singular points, two
of which are only apparent, i.e., have trivial monodromy.  So
$\tilde{\mathfrak{E}}\stackrel{\tilde\pi}{\to}\mathbb{P}^1(\mathbb{C})_s$
has $8$~singular fibres.  They are located above
\begin{equation*}
s=0,\infty,-48,-60,-72,-80,-96,-120,
\end{equation*}
as one would expect from~(\ref{eq:rama6}).  Having so many singular fibres,
this elliptic surface is more complicated than those classified by
Schmickler-Hirzebruch or by Herfurtner~\cite{Herfurtner91}.

A theory of modular equations for {\em general\/} rational elliptic
surfaces must produce algebraic transformation laws, of arbitrarily high
degree, for Gauss--Manin connections pulled back from the universal
elliptic family.  It is clear from the above analysis that such modular
equations should be parametrized by elliptic surfaces, which may sometimes
themselves be rational.  In the absence of a general theory, whether there
are additional rationally parametrized modular equations in Ramanujan's
theory of signature~$6$ remains unclear.

\appendix
\renewcommand{\theequation}{A.\arabic{equation}}
\setcounter{equation}{0}

\section*{Appendix.  Hypergeometric and Heun Equations}

Any second-order Fuchsian differential equation
on~$\mathbb{P}^1(\mathbb{C})_t$ can be reduced to the normal form
\begin{equation}
\label{eq:hyper}
D_t^2\, u + \left[\frac{c}{t} + \frac{d}{t-1}\right]
D_t u + \left[\frac{ab}{t(t-1)}\right]u = 0,
\end{equation}
if it has exactly three singular points, and to the normal form
\begin{equation}
\label{eq:heun}
D_t^2\, u
+ \left[ \frac\gamma t + \frac\delta{t-1} + \frac\epsilon{t-a}
  \right]D_t u + \left[\frac{\alpha\beta\, t - q}{t(t-1)(t-a)}\right]u = 0,
\end{equation}
if it has exactly four.  The reduction employs a M\"obius transformation to
reposition (three~of) the singular points to $t=0,1,\infty$.  (In the
latter case the fourth singular point moves to some location
$a\in\mathbb{C}\setminus\{0,1\}$.)  It~also employs index transformations,
to shift one characteristic exponent at each of the finite singular points
to zero.  The exponents of (\ref{eq:hyper}) and~(\ref{eq:heun}) are
respectively
\begin{alignat*}{3}
&0,1-c;\,0,1-d;\,a,b &\qquad&\text{at}&\quad&t=0,1,\infty \\
&0,1-\gamma;\,0,1-\delta;\,0,1-\epsilon;\alpha,\beta &\qquad&\text{at}&\quad&t=0,1,a,\infty.
\end{alignat*}
Fuchs's relation (the sum of the $2k$~exponents of a Fuchsian differential
equation with $k$~singular points on~$\mathbb{P}^1(\mathbb{C})$ equaling
$k-2$) implies that $d=a+b-c+1,$ resp.\
$\epsilon=\alpha+\beta-\gamma-\delta+1$.  The quantity $q\in\mathbb{C}$
in~(\ref{eq:heun}) is an accessory parameter, which does not affect the
local monodromy about each singular point.

Each of (\ref{eq:hyper}),(\ref{eq:heun}) has a unique local solution
holomorphic at~$t=0$ and equal to unity there.  For~(\ref{eq:hyper}) this
solution is the Gauss\ hypergeometric function ${}_2F_1(a,b;c;\cdot),$ and
for~(\ref{eq:heun}) it is the local Heun function~\cite{Ronveaux95}, widely
denoted $\Hl(a,q;\alpha,\beta,\gamma,\delta;\cdot)$.  These have
$t$-expansions $\sum_{n=0}^\infty c_nt^n,$ which converge on
$\left|t\right|<1,$ resp.\
$\left|t\right|<\min\left(1,\left|a\right|\right)$.  The coefficients
satisfy respective recurrences
\begin{equation}
\label{eq:hyperrecurrence}
(n+a-1)(n+b-1) \,c_{n-1} - n(n+c-1)\,c_n   = 0
\end{equation}
and
\begin{multline}
(n+\alpha-1)(n+\beta-1)\, c_{n-1}\\
\quad{}- \bigl\{n \left[\,(n+\gamma+\delta-1)a + (n+\gamma+\epsilon-1)\,\right] + q\bigr\}\,c_n\\
+(n+1)(n+\gamma)a\, c_{n+1} = 0,
\end{multline}
initialized by $c_0=1$ (and $c_{-1}=0$ in the latter case).

In the spirit of Kummer, one can derive functional equations satisfied by
the parametrized functions ${}_2F_1$ and~$\Hl$ by applying M\"obius and
index transformations to their defining differential equations.  The idea
is that one should transform the equation to itself, with altered
parameters and argument~\cite{Maier10}.  If there are $k$~singular points
($k=3,4$ here), there will be a subgroup of the M\"obius group, isomorphic
to~$\mathfrak{S}_k,$ that permutes them.  Also, there will be a subgroup of
the group of index transformations, isomorphic to~$(C_2)^{k-1},$ that
negates the nonzero exponents of the $k-1$ finite singular points.  The
first group normalizes the latter, so the group of composite
transformations, i.e., the transformation group that they generate, is a
semidirect product $C_2^{k-1}\rtimes \mathfrak{S}_k$.  By examination, this
semidirect product is isomorphic to the group of {\em even-signed\/}
permutations of $k$~objects, which is an index\nobreakdash-$2$ subgroup of
the wreath product $C_2 \wr \mathfrak{S}_k,$ the group of {\em signed\/}
permutations of $k$~objects.

The action of the order-$24$ group $\left[C_2\wr\mathfrak{S}_3\right]_{\rm
even}$\! on~(\ref{eq:hyper}) yields $24$ local solutions, each expressed in
terms of~${}_2F_1$.  These are the famous $24$ solutions of Kummer.  The
action of the order-$192$ group $\left[C_2 \wr\mathfrak{S}_4\right]_{\rm
even}$\! on~(\ref{eq:heun}) yields $192$ local solutions, each expressed
in~terms of~$\Hl,$ which are less well known than Kummer's (but see
Refs.~\cite{Maier10,Ronveaux95}).  If the transformation in $\left[C_2
\wr\mathfrak{S}_k\right]_{\rm even}$\! applied to (\ref{eq:hyper}), resp.\
(\ref{eq:heun}), is to yield a function equal to ${}_2F_1,$ resp.\ ${\it
Hl},$ then it should stabilize the singular point $t=0,$ and perform
no~index transformation there.  Hence, the subgroup of allowed
transformations is isomorphic to $[C_2 \wr\mathfrak{S}_{k-1}]_{\rm even}$.

So, there is a group of transformations of~${}_2F_1$ isomorphic to
$\left[C_2 \wr\mathfrak{S}_{2}\right]_{\rm even},$ i.e., to the four-group
$C_2\times C_2$.  It is generated by classical transformations named after
Euler and Pfaff.  Pfaff's is used in this article.  It is
\begin{equation}
\label{eq:Pfaff}
{}_2F_1(a,b;\,c;\,t)=(1-t)^{-a}\: {}_2F_1(a,c-b;\,c;\,\tfrac{t}{t-1}),
\end{equation}
and arises from the M\"obius transformation that interchanges $t=1,\infty$.

The group of transformations of~$\Hl$ is isomorphic to $\left[C_2
\wr\mathfrak{S}_{3}\right]_{\rm even},$ which by examination is isomorphic
to the octahedral group~$\mathfrak{S}_4$.  Two of the transformations in
this group are used in this article.  They are the relatively trivial one
\begin{equation}
  \Hl (a,q;\,\alpha,\beta,\gamma,\delta;\,t)={\it
  Hl}(\tfrac{1}a,\tfrac{q}a;\,
  \alpha,\beta,\gamma,\alpha+\beta-\gamma-\delta+1;\,\tfrac{t}a),
\end{equation}
and the generalized Pfaff transformation
\begin{multline}
  \Hl (a,q;\,\alpha,\beta,\gamma,\delta;\,t) = \\ 
  (1-t)^{-\alpha}\, {\it
  Hl}(\tfrac{a}{a-1},\tfrac{-q+\gamma\alpha
  a}{a-1};\,\alpha,\alpha-\delta+1,\gamma,\alpha-\beta+1;\,\tfrac{t}{t-1}).
\label{eq:genPfaff}
\end{multline}
They arise respectively from the M\"obius transformation that in~effect
interchanges $t=1,a$ (i.e., takes $a$~to~$1$ and $1$~to~$a':=1/a$), and the
one that interchanges $t=1,\infty$ (and less importantly, takes
$t=a$~to~$t=a':=a/(a-1)$).

\begin{acknowledgement}
  The tables were prepared with {\sc Maxima}, the free version of the {\sc
  Macsyma} computer algebra system.  The author is indebted to Neil Sloane
  for his {\em On-Line Encyclopedia of Integer Sequences\/}, to volunteers
  Michael Somos and Vladeta Jovovi\'c, among others, for annotating~it; and
  to Mark van Hoeij for his {\tt algcurves} package.
\end{acknowledgement}


\end{document}